
\ifx\shlhetal\undefinedcontrolsequence\let\shlhetal\relax\fi

\input amstex
\expandafter\ifx\csname mathdefs.tex\endcsname\relax
  \expandafter\gdef\csname mathdefs.tex\endcsname{}
\else \message{Hey!  Apparently you were trying to
  \string\input{mathdefs.tex} twice.   This does not make sense.} 
\errmessage{Please edit your file (probably \jobname.tex) and remove
any duplicate ``\string\input'' lines}\endinput\fi




\catcode`\X=12\catcode`\@=11

\def\n@wcount{\alloc@0\count\countdef\insc@unt}
\def\n@wwrite{\alloc@7\write\chardef\sixt@@n}
\def\n@wread{\alloc@6\read\chardef\sixt@@n}
\def\r@s@t{\relax}\def\v@idline{\par}\def\@mputate#1/{#1}
\def\l@c@l#1X{\firstpart.#1}\def\gl@b@l#1X{#1}\def\t@d@l#1X{{}}

\def\crossrefs#1{\ifx\all#1\let\tr@ce=\all\else\def\tr@ce{#1,}\fi
   \n@wwrite\cit@tionsout\openout\cit@tionsout=\jobname.cit 
   \write\cit@tionsout{\tr@ce}\expandafter\setfl@gs\tr@ce,}
\def\setfl@gs#1,{\def\@{#1}\ifx\@\empty\let\next=\relax
   \else\let\next=\setfl@gs\expandafter\xdef
   \csname#1tr@cetrue\endcsname{}\fi\next}
\def\m@ketag#1#2{\expandafter\n@wcount\csname#2tagno\endcsname
     \csname#2tagno\endcsname=0\let\tail=\all\xdef\all{\tail#2,}
   \ifx#1\l@c@l\let\tail=\r@s@t\xdef\r@s@t{\csname#2tagno\endcsname=0\tail}\fi
   \expandafter\gdef\csname#2cite\endcsname##1{\expandafter
     \ifx\csname#2tag##1\endcsname\relax?\else\csname#2tag##1\endcsname\fi
     \expandafter\ifx\csname#2tr@cetrue\endcsname\relax\else
     \write\cit@tionsout{#2tag ##1 cited on page \folio.}\fi}
   \expandafter\gdef\csname#2page\endcsname##1{\expandafter
     \ifx\csname#2page##1\endcsname\relax?\else\csname#2page##1\endcsname\fi
     \expandafter\ifx\csname#2tr@cetrue\endcsname\relax\else
     \write\cit@tionsout{#2tag ##1 cited on page \folio.}\fi}
   \expandafter\gdef\csname#2tag\endcsname##1{\expandafter
      \ifx\csname#2check##1\endcsname\relax
      \expandafter\xdef\csname#2check##1\endcsname{}%
      \else\immediate\write16{Warning: #2tag ##1 used more than once.}\fi
      \multit@g{#1}{#2}##1/X%
      \write\t@gsout{#2tag ##1 assigned number \csname#2tag##1\endcsname\space
      on page \number\count0.}%
   \csname#2tag##1\endcsname}}

\def\multit@g#1#2#3/#4X{\def\t@mp{#4}\ifx\t@mp\empty%
      \global\advance\csname#2tagno\endcsname by 1 
      \expandafter\xdef\csname#2tag#3\endcsname
      {#1\number\csname#2tagno\endcsnameX}%
   \else\expandafter\ifx\csname#2last#3\endcsname\relax
      \expandafter\n@wcount\csname#2last#3\endcsname
      \global\advance\csname#2tagno\endcsname by 1 
      \expandafter\xdef\csname#2tag#3\endcsname
      {#1\number\csname#2tagno\endcsnameX}
      \write\t@gsout{#2tag #3 assigned number \csname#2tag#3\endcsname\space
      on page \number\count0.}\fi
   \global\advance\csname#2last#3\endcsname by 1
   \def\t@mp{\expandafter\xdef\csname#2tag#3/}%
   \expandafter\t@mp\@mputate#4\endcsname
   {\csname#2tag#3\endcsname\lastpart{\csname#2last#3\endcsname}}\fi}
\def\t@gs#1{\def\all{}\m@ketag#1e\m@ketag#1s\m@ketag\t@d@l p
\let\realscite\scite
\let\realstag\stag
   \m@ketag\gl@b@l r \n@wread\t@gsin
   \openin\t@gsin=\jobname.tgs \re@der \closein\t@gsin
   \n@wwrite\t@gsout\openout\t@gsout=\jobname.tgs }
\outer\def\localtags{\t@gs\l@c@l}
\outer\def\globaltags{\t@gs\gl@b@l}
\outer\def\newlocaltag#1{\m@ketag\l@c@l{#1}}
\outer\def\newglobaltag#1{\m@ketag\gl@b@l{#1}}

\newif\ifpr@ 
\def\m@kecs #1tag #2 assigned number #3 on page #4.%
   {\expandafter\gdef\csname#1tag#2\endcsname{#3}
   \expandafter\gdef\csname#1page#2\endcsname{#4}
   \ifpr@\expandafter\xdef\csname#1check#2\endcsname{}\fi}
\def\re@der{\ifeof\t@gsin\let\next=\relax\else
   \read\t@gsin to\t@gline\ifx\t@gline\v@idline\else
   \expandafter\m@kecs \t@gline\fi\let \next=\re@der\fi\next}
\def\pretags#1{\pr@true\pret@gs#1,,}
\def\pret@gs#1,{\def\@{#1}\ifx\@\empty\let\n@xtfile=\relax
   \else\let\n@xtfile=\pret@gs \openin\t@gsin=#1.tgs \message{#1} \re@der 
   \closein\t@gsin\fi \n@xtfile}

\newcount\sectno\sectno=0\newcount\subsectno\subsectno=0
\newif\ifultr@local \def\ultralocal{\ultr@localtrue}
\def\firstpart{\number\sectno}
\def\lastpart#1{\ifcase#1 \or a\or b\or c\or d\or e\or f\or g\or h\or 
   i\or k\or l\or m\or n\or o\or p\or q\or r\or s\or t\or u\or v\or w\or 
   x\or y\or z \fi}

\def\resetall{\global\advance\sectno by 1\subsectno=0
   \gdef\firstpart{\number\sectno}\r@s@t}
\def\resetsub{\global\advance\subsectno by 1
   \gdef\firstpart{\number\sectno.\number\subsectno}\r@s@t}
\def\newsection#1\par{\resetall\vskip0pt plus.3\vsize\penalty-250
   \vskip0pt plus-.3\vsize\bigskip\bigskip
   \message{#1}\leftline{\bf#1}\nobreak\bigskip}
\def\subsection#1\par{\ifultr@local\resetsub\fi
   \vskip0pt plus.2\vsize\penalty-250\vskip0pt plus-.2\vsize
   \bigskip\smallskip\message{#1}\leftline{\bf#1}\nobreak\medskip}


\newdimen\marginshift

\newdimen\margindelta
\newdimen\marginmax
\newdimen\marginmin

\def\margininit{       
\marginmax=3 true cm                  
				      
\margindelta=0.1 true cm              
\marginmin=0.1true cm                 
\marginshift=\marginmin
}    

\def\t@gsjj#1,{\def\@{#1}\ifx\@\empty\let\next=\relax\else\let\next=\t@gsjj
   \def\@@{p}\ifx\@\@@\else
   \expandafter\gdef\csname#1cite\endcsname##1{\citejj{##1}}
   \expandafter\gdef\csname#1page\endcsname##1{?}
   \expandafter\gdef\csname#1tag\endcsname##1{\tagjj{##1}}\fi\fi\next}
\newif\ifshowstuffinmargin
\showstuffinmarginfalse
\def\jjtags{\ifx\shlhetal\relax 
  \else
\ifx\shlhetal\undefinedcontrolseq
\else
\showstuffinmargintrue
\ifx\all\relax\else\expandafter\t@gsjj\all,\fi\fi \fi
}

\def\tagjj#1{\realstag{#1}\oldmginpar{\zeigen{#1}}}
\def\citejj#1{\rechnen{#1}\mginpar{\zeigen{#1}}}     

\def\rechnen#1{\expandafter\ifx\csname stag#1\endcsname\relax ??\else
                           \csname stag#1\endcsname\fi}

\newdimen\theight

\def\marginfont{\sevenrm}

\def\trymarginbox#1{\setbox0=\hbox{\marginfont\hskip\marginshift #1}%
		\global\marginshift\wd0 
		\global\advance\marginshift\margindelta}

\def \oldmginpar#1{%
\ifvmode\setbox0\hbox to \hsize{\hfill\rlap{\marginfont\quad#1}}%
\ht0 0cm
\dp0 0cm
\box0\vskip-\baselineskip
\else 
             \vadjust{\trymarginbox{#1}%
		\ifdim\marginshift>\marginmax \global\marginshift\marginmin
			\trymarginbox{#1}%
                \fi
             \theight=\ht0
             \advance\theight by \dp0    \advance\theight by \lineskip
             \kern -\theight \vbox to \theight{\rightline{\rlap{\box0}}%
\vss}}\fi}

\newdimen\upordown
\global\upordown=8pt
\font\tinyfont=cmtt8 
\def\mginpar#1{\smash{\hbox to 0cm{\kern-10pt\raise7pt\hbox{\tinyfont #1}\hss}}}
\def\mginpar#1{{\hbox to 0cm{\kern-10pt\raise\upordown\hbox{\tinyfont #1}\hss}}\global\upordown-\upordown}


\def\t@gsoff#1,{\def\@{#1}\ifx\@\empty\let\next=\relax\else\let\next=\t@gsoff
   \def\@@{p}\ifx\@\@@\else
   \expandafter\gdef\csname#1cite\endcsname##1{\zeigen{##1}}
   \expandafter\gdef\csname#1page\endcsname##1{?}
   \expandafter\gdef\csname#1tag\endcsname##1{\zeigen{##1}}\fi\fi\next}
\def\verbatimtags{\showstuffinmarginfalse
\ifx\all\relax\else\expandafter\t@gsoff\all,\fi}
\def\zeigen#1{\hbox{$\scriptstyle\langle$}#1\hbox{$\scriptstyle\rangle$}}


\def\margintag#1{\ifshowstuffinmargin\oldmginpar{\zeigen{#1}}\fi}

\def\(#1){\edef\dot@g{\ifmmode\ifinner(\hbox{\noexpand\etag{#1}})
   \else\noexpand\eqno(\hbox{\noexpand\etag{#1}})\fi
   \else(\noexpand\ecite{#1})\fi}\dot@g}

\newif\ifbr@ck
\def\eat#1{}
\def\[#1]{\br@cktrue[\br@cket#1'X]}
\def\br@cket#1'#2X{\def\temp{#2}\ifx\temp\empty\let\next\eat
   \else\let\next\br@cket\fi
   \ifbr@ck\br@ckfalse\br@ck@t#1,X\else\br@cktrue#1\fi\next#2X}
\def\br@ck@t#1,#2X{\def\temp{#2}\ifx\temp\empty\let\neext\eat
   \else\let\neext\br@ck@t\def\temp{,}\fi
   \def\teemp{#1}\ifx\teemp\empty\else\rcite{#1}\fi\temp\neext#2X}
\def\resetbr@cket{\gdef\[##1]{[\rtag{##1}]}}
\def\references{\resetbr@cket\newsection References\par}

\newtoks\symb@ls\newtoks\s@mb@ls\newtoks\p@gelist\n@wcount\ftn@mber
    \ftn@mber=1\newif\ifftn@mbers\ftn@mbersfalse\newif\ifbyp@ge\byp@gefalse
\def\defm@rk{\ifftn@mbers\n@mberm@rk\else\symb@lm@rk\fi}
\def\n@mberm@rk{\xdef\m@rk{{\the\ftn@mber}}%
    \global\advance\ftn@mber by 1 }
\def\rot@te#1{\let\temp=#1\global#1=\expandafter\r@t@te\the\temp,X}
\def\r@t@te#1,#2X{{#2#1}\xdef\m@rk{{#1}}}
\def\b@@st#1{{$^{#1}$}}\def\str@p#1{#1}
\def\symb@lm@rk{\ifbyp@ge\rot@te\p@gelist\ifnum\expandafter\str@p\m@rk=1 
    \s@mb@ls=\symb@ls\fi\write\f@nsout{\number\count0}\fi \rot@te\s@mb@ls}
\def\byp@ge{\byp@getrue\n@wwrite\f@nsin\openin\f@nsin=\jobname.fns 
    \n@wcount\currentp@ge\currentp@ge=0\p@gelist={0}
    \re@dfns\closein\f@nsin\rot@te\p@gelist
    \n@wread\f@nsout\openout\f@nsout=\jobname.fns }
\def\m@kelist#1X#2{{#1,#2}}
\def\re@dfns{\ifeof\f@nsin\let\next=\relax\else\read\f@nsin to \f@nline
    \ifx\f@nline\v@idline\else\let\t@mplist=\p@gelist
    \ifnum\currentp@ge=\f@nline
    \global\p@gelist=\expandafter\m@kelist\the\t@mplistX0
    \else\currentp@ge=\f@nline
    \global\p@gelist=\expandafter\m@kelist\the\t@mplistX1\fi\fi
    \let\next=\re@dfns\fi\next}
\def\symbols#1{\symb@ls={#1}\s@mb@ls=\symb@ls} 
\def\bigsymbol{\textstyle}
\symbols{\bigsymbol\ast,\dagger,\ddagger,\sharp,\flat,\natural,\star}
\def\ftnumbers{\ftn@mberstrue} \def\ftsymbols{\ftn@mbersfalse}
\def\paginal{\byp@ge} \def\resetftnumbers{\ftn@mber=1}
\def\ftnote#1{\defm@rk\expandafter\expandafter\expandafter\footnote
    \expandafter\b@@st\m@rk{#1}}

\long\def\jump#1\endjump{}
\def\ssum{\mathop{\lower .1em\hbox{$\textstyle\Sigma$}}\nolimits}

\def\qed{\nobreak\kern 1em \vrule height .5em width .5em depth 0em}
\def\newneq{\hbox{\rlap{\hbox to 1\wd9{\hss$=$\hss}}\raise .1em 
   \hbox to 1\wd9{\hss$\scriptscriptstyle/$\hss}}}
\def\subsetne{\setbox9 = \hbox{$\subset$}\mathrel{\hbox{\rlap
   {\lower .4em \newneq}\raise .13em \hbox{$\subset$}}}}
\def\supsetne{\setbox9 = \hbox{$\subset$}\mathrel{\hbox{\rlap
   {\lower .4em \newneq}\raise .13em \hbox{$\supset$}}}}

\def\vbar{\mathchoice{\vrule height6.3ptdepth-.5ptwidth.8pt\kern-.8pt}
   {\vrule height6.3ptdepth-.5ptwidth.8pt\kern-.8pt}
   {\vrule height4.1ptdepth-.35ptwidth.6pt\kern-.6pt}
   {\vrule height3.1ptdepth-.25ptwidth.5pt\kern-.5pt}}
\def\f@dge{\mathchoice{}{}{\mkern.5mu}{\mkern.8mu}}
\def\b@c#1#2{{\rm \mkern#2mu\vbar\mkern-#2mu#1}}
\def\b@b#1{{\rm I\mkern-3.5mu #1}}
\def\b@a#1#2{{\rm #1\mkern-#2mu\f@dge #1}}
\def\bb#1{{\count4=`#1 \advance\count4by-64 \ifcase\count4\or\b@a A{11.5}\or
   \b@b B\or\b@c C{5}\or\b@b D\or\b@b E\or\b@b F \or\b@c G{5}\or\b@b H\or
   \b@b I\or\b@c J{3}\or\b@b K\or\b@b L \or\b@b M\or\b@b N\or\b@c O{5} \or
   \b@b P\or\b@c Q{5}\or\b@b R\or\b@a S{8}\or\b@a T{10.5}\or\b@c U{5}\or
   \b@a V{12}\or\b@a W{16.5}\or\b@a X{11}\or\b@a Y{11.7}\or\b@a Z{7.5}\fi}}

\catcode`\X=11 \catcode`\@=12




\let\thischap\jobname

\def\partof#1{\csname returnthe#1part\endcsname}
\def\chapof#1{\csname returnthe#1chap\endcsname}

\def\setchapter#1,#2,#3;{%
  \expandafter\def\csname returnthe#1part\endcsname{#2}%
  \expandafter\def\csname returnthe#1chap\endcsname{#3}%
}

\setchapter 300a,A,II.A;
\setchapter 300b,A,II.B;
\setchapter 300c,A,II.C;
\setchapter 300d,A,II.D;
\setchapter 300e,A,II.E;
\setchapter 300f,A,II.F;
\setchapter 300g,A,II.G;
\setchapter  E53,B,N;
\setchapter  88r,B,I;
\setchapter  600,B,III;
\setchapter  705,B,IV;
\setchapter  734,B,V;

\def\cprefix#1{
\edef\theotherpart{\partof{#1}}\edef\theotherchap{\chapof{#1}}%
\ifx\theotherpart\thispart
   \ifx\theotherchap\thischap 
    \else 
     \theotherchap%
    \fi
   \else 
     \theotherchap\fi}

\def\sectioncite[#1]#2{%
     \cprefix{#2}#1}

\edef\thispart{\partof{\thischap}}
\edef\thischap{\chapof{\thischap}}

\def\lastpage of '#1' is #2.{\expandafter\def\csname lastpage#1\endcsname{#2}}


\def\spuriousreset{}


\expandafter\ifx\csname citeadd.tex\endcsname\relax
\expandafter\gdef\csname citeadd.tex\endcsname{}
\else \message{Hey!  Apparently you were trying to
\string\input{citeadd.tex} twice.   This does not make sense.} 
\errmessage{Please edit your file (probably \jobname.tex) and remove
any duplicate ``\string\input'' lines}\endinput\fi

\sectno=-1   
\localtags
\jjtags
\NoBlackBoxes
\define\mr{\medskip\roster}
\define\sn{\smallskip\noindent}
\define\mn{\medskip\noindent}
\define\rest{\restriction}
\define\bn{\bigskip\noindent}
\define\ub{\underbar}
\define\wilog{\text{without loss of generality}}
\define\ermn{\endroster\medskip\noindent}

\define\dbcu{\dsize\bigcup}
\define \nl{\newline}
\newbox\noforkbox \newdimen\forklinewidth
\forklinewidth=0.3pt   
\setbox0\hbox{$\textstyle\bigcup$}
\setbox1\hbox to \wd0{\hfil\vrule width \forklinewidth depth \dp0
             height \ht0 \hfil}
\wd1=0 cm
\setbox\noforkbox\hbox{\box1\box0\relax}
\def\unionstick{\mathop{\copy\noforkbox}\limits}
\def\nonfork#1#2_#3{#1\unionstick_{\textstyle #3}#2}
\def\nonforkin#1#2_#3^#4{#1\unionstick_{\textstyle #3}^{\textstyle #4}#2}     
%
\setbox0\hbox{$\textstyle\bigcup$}
\setbox1\hbox to \wd0{\hfil{\sl /\/}\hfil}
\setbox2\hbox to \wd0{\hfil\vrule height \ht0 depth \dp0 width
                                \forklinewidth\hfil}
\wd1=0cm
\wd2=0cm
\newbox\doesforkbox
\setbox\doesforkbox\hbox{\box1\box0\relax}
\def\nunionstick{\mathop{\copy\doesforkbox}\limits}

\def\fork#1#2_#3{#1\nunionstick_{\textstyle #3}#2}
\def\forkin#1#2_#3^#4{#1\nunionstick_{\textstyle #3}^{\textstyle #4}#2}     

\magnification=\magstep 1
\documentstyle{amsppt}

{    
\catcode`@11

\ifx\alicetwothousandloaded@\relax
  \endinput\else\global\let\alicetwothousandloaded@\relax\fi

\gdef\subjclass{\let\savedef@\subjclass
 \def\subjclass##1\endsubjclass{\let\subjclass\savedef@
   \toks@{\def\usualspace{{\rm\enspace}}\eightpoint}%
   \toks@@{##1\unskip.}%
   \edef\thesubjclass@{\the\toks@
     \frills@{{\noexpand\rm2000 {\noexpand\it Mathematics Subject
       Classification}.\noexpand\enspace}}%
     \the\toks@@}}%
  \nofrillscheck\subjclass}
} 


\expandafter\ifx\csname alice2jlem.tex\endcsname\relax
  \expandafter\xdef\csname alice2jlem.tex\endcsname{\the\catcode`@}
\else \message{Hey!  Apparently you were trying to
\string\input{alice2jlem.tex}  twice.   This does not make sense.}
\errmessage{Please edit your file (probably \jobname.tex) and remove
any duplicate ``\string\input'' lines}\endinput\fi

\expandafter\ifx\csname bib4plain.tex\endcsname\relax
  \expandafter\gdef\csname bib4plain.tex\endcsname{}
\else \message{Hey!  Apparently you were trying to \string\input
  bib4plain.tex twice.   This does not make sense.}
\errmessage{Please edit your file (probably \jobname.tex) and remove
any duplicate ``\string\input'' lines}\endinput\fi

\def\renewcommand{\newcommand}	       
\edef\cite{\the\catcode`@}%
\catcode`@ = 11
\let\@oldatcatcode = \cite
\chardef\@letter = 11
\chardef\@other = 12
%
%
%
%
\def\@innerdef#1#2{\edef#1{\expandafter\noexpand\csname #2\endcsname}}%
%
%
\@innerdef\@innernewcount{newcount}%
\@innerdef\@innernewdimen{newdimen}%
\@innerdef\@innernewif{newif}%
\@innerdef\@innernewwrite{newwrite}%
%
%
%
\def\@gobble#1{}%
%
%
%
\ifx\inputlineno\@undefined
   \let\@linenumber = \empty 
\else
   \def\@linenumber{\the\inputlineno:\space}%
\fi
%
%
%
\def\@futurenonspacelet#1{\def\cs{#1}%
   \afterassignment\@stepone\let\@nexttoken=
}%
\begingroup 
\def\\{\global\let\@stoken= }%
\\ 
\endgroup
\def\@stepone{\expandafter\futurelet\cs\@steptwo}%
\def\@steptwo{\expandafter\ifx\cs\@stoken\let\@@next=\@stepthree
   \else\let\@@next=\@nexttoken\fi \@@next}%
\def\@stepthree{\afterassignment\@stepone\let\@@next= }%
%
%
%
\def\@getoptionalarg#1{%
   \let\@optionaltemp = #1%
   \let\@optionalnext = \relax
   \@futurenonspacelet\@optionalnext\@bracketcheck
}%
%
%
\def\@bracketcheck{%
   \ifx [\@optionalnext
      \expandafter\@@getoptionalarg
   \else
      \let\@optionalarg = \empty
      \expandafter\@optionaltemp
   \fi
}%
\def\@@getoptionalarg[#1]{%
   \def\@optionalarg{#1}%
   \@optionaltemp
}%
%
%
%
\def\@nnil{\@nil}%
\def\@fornoop#1\@@#2#3{}%
\def\@for#1:=#2\do#3{%
   \edef\@fortmp{#2}%
   \ifx\@fortmp\empty \else
      \expandafter\@forloop#2,\@nil,\@nil\@@#1{#3}%
   \fi
}%
\def\@forloop#1,#2,#3\@@#4#5{\def#4{#1}\ifx #4\@nnil \else
       #5\def#4{#2}\ifx #4\@nnil \else#5\@iforloop #3\@@#4{#5}\fi\fi
}%
\def\@iforloop#1,#2\@@#3#4{\def#3{#1}\ifx #3\@nnil
       \let\@nextwhile=\@fornoop \else
      #4\relax\let\@nextwhile=\@iforloop\fi\@nextwhile#2\@@#3{#4}%
}%
%
%
%
\@innernewif\if@fileexists
\def\@testfileexistence{\@getoptionalarg\@finishtestfileexistence}%
\def\@finishtestfileexistence#1{%
   \begingroup
      \def\extension{#1}%
      \immediate\openin0 =
         \ifx\@optionalarg\empty\jobname\else\@optionalarg\fi
         \ifx\extension\empty \else .#1\fi
         \space
      \ifeof 0
         \global\@fileexistsfalse
      \else
         \global\@fileexiststrue
      \fi
      \immediate\closein0
   \endgroup
}%
%
%
%
%
\def\bibliographystyle#1{%
   \@readauxfile
   \@writeaux{\string\bibstyle{#1}}%
}%
\let\bibstyle = \@gobble
%
%
\let\bblfilebasename = \jobname
\def\bibliography#1{%
   \@readauxfile
   \@writeaux{\string\bibdata{#1}}%
   \@testfileexistence[\bblfilebasename]{bbl}%
   \if@fileexists
      \nobreak
      \@readbblfile
   \fi
}%
\let\bibdata = \@gobble
%
%
\def\nocite#1{%
   \@readauxfile
   \@writeaux{\string\citation{#1}}%
}%
\@innernewif\if@notfirstcitation
%
%
\def\cite{\@getoptionalarg\@cite}%
%
%
\def\@cite#1{%
   \let\@citenotetext = \@optionalarg
   \printcitestart
   \nocite{#1}%
   \@notfirstcitationfalse
   \@for \@citation :=#1\do
   {%
      \expandafter\@onecitation\@citation\@@
   }%
   \ifx\empty\@citenotetext\else
      \printcitenote{\@citenotetext}%
   \fi
   \printcitefinish
}%
\newif\ifweareinprivate
\weareinprivatetrue
\ifx\shlhetal\undefinedcontrolseq\weareinprivatefalse\fi
\ifx\shlhetal\relax\weareinprivatefalse\fi
\def\@onecitation#1\@@{%
   \if@notfirstcitation
      \printbetweencitations
   \fi
   \expandafter \ifx \csname\@citelabel{#1}\endcsname \relax
      \if@citewarning
         \message{\@linenumber Undefined citation `#1'.}%
      \fi
     \ifweareinprivate
      \expandafter\gdef\csname\@citelabel{#1}\endcsname{%
\strut 
\vadjust{\vskip-\dp\strutbox
\vbox to 0pt{\vss\parindent0cm \leftskip=\hsize 
\advance\leftskip3mm
\advance\hsize 4cm\strut\openup-4pt 
\rightskip 0cm plus 1cm minus 0.5cm ?  #1 ?\strut}}
         {\tt
            \escapechar = -1
            \nobreak\hskip0pt\pfeilsw
            \expandafter\string\csname#1\endcsname
             \pfeilso
            \nobreak\hskip0pt
         }%
      }%
     \else  
      \expandafter\gdef\csname\@citelabel{#1}\endcsname{%
            {\tt\expandafter\string\csname#1\endcsname}
      }%
     \fi  
   \fi
   \csname\@citelabel{#1}\endcsname
   \@notfirstcitationtrue
}%
%
%
\def\@citelabel#1{b@#1}%
%
%
\def\@citedef#1#2{\expandafter\gdef\csname\@citelabel{#1}\endcsname{#2}}%
%
%
%
\def\@readbblfile{%
   \ifx\@itemnum\@undefined
      \@innernewcount\@itemnum
   \fi
   \begingroup
      \def\begin##1##2{%
         \setbox0 = \hbox{\biblabelcontents{##2}}%
         \biblabelwidth = \wd0
      }%
      \def\end##1{}
      %
      %
      \@itemnum = 0
      \def\bibitem{\@getoptionalarg\@bibitem}%
      \def\@bibitem{%
         \ifx\@optionalarg\empty
            \expandafter\@numberedbibitem
         \else
            \expandafter\@alphabibitem
         \fi
      }%
      \def\@alphabibitem##1{%
         \expandafter \xdef\csname\@citelabel{##1}\endcsname {\@optionalarg}%
         \ifx\biblabelprecontents\@undefined
            \let\biblabelprecontents = \relax
         \fi
         \ifx\biblabelpostcontents\@undefined
            \let\biblabelpostcontents = \hss
         \fi
         \@finishbibitem{##1}%
      }%
      \def\@numberedbibitem##1{%
         \advance\@itemnum by 1
         \expandafter \xdef\csname\@citelabel{##1}\endcsname{\number\@itemnum}%
         \ifx\biblabelprecontents\@undefined
            \let\biblabelprecontents = \hss
         \fi
         \ifx\biblabelpostcontents\@undefined
            \let\biblabelpostcontents = \relax
         \fi
         \@finishbibitem{##1}%
      }%
      \def\@finishbibitem##1{%
         \biblabelprint{\csname\@citelabel{##1}\endcsname}%
         \@writeaux{\string\@citedef{##1}{\csname\@citelabel{##1}\endcsname}}%
         \ignorespaces
      }%
      %
      %
      \let\em = \bblem
      \let\newblock = \bblnewblock
      \let\sc = \bblsc
      \frenchspacing
      \clubpenalty = 4000 \widowpenalty = 4000
      \tolerance = 10000 \hfuzz = .5pt
      \everypar = {\hangindent = \biblabelwidth
                      \advance\hangindent by \biblabelextraspace}%
      \bblrm
      \parskip = 1.5ex plus .5ex minus .5ex
      \biblabelextraspace = .5em
      \bblhook
      \input \bblfilebasename.bbl
   \endgroup
}%
%
%
\@innernewdimen\biblabelwidth
\@innernewdimen\biblabelextraspace
%
%
%
\def\biblabelprint#1{%
   \noindent
   \hbox to \biblabelwidth{%
      \biblabelprecontents
      \biblabelcontents{#1}%
      \biblabelpostcontents
   }%
   \kern\biblabelextraspace
}%
%
%
%
\def\biblabelcontents#1{{\bblrm [#1]}}%
%
%
\def\bblrm{\rm}%
%
%
\def\bblem{\it}%
%
%
\def\bblsc{\ifx\@scfont\@undefined
              \font\@scfont = cmcsc10
           \fi
           \@scfont
}%
%
%
\def\bblnewblock{\hskip .11em plus .33em minus .07em }%
%
%
\let\bblhook = \empty
%
%
%
\def\printcitestart{[}
\def\printcitefinish{]}
\def\printbetweencitations{, }
\def\printcitenote#1{, #1}
%
%
%
\let\citation = \@gobble
%
%
%
\@innernewcount\@numparams
%
%
\def\newcommand#1{%
   \def\@commandname{#1}%
   \@getoptionalarg\@continuenewcommand
}%
%
%
\def\@continuenewcommand{%
   \@numparams = \ifx\@optionalarg\empty 0\else\@optionalarg \fi \relax
   \@newcommand
}%
%
%
\def\@newcommand#1{%
   \def\@startdef{\expandafter\edef\@commandname}%
   \ifnum\@numparams=0
      \let\@paramdef = \empty
   \else
      \ifnum\@numparams>9
         \errmessage{\the\@numparams\space is too many parameters}%
      \else
         \ifnum\@numparams<0
            \errmessage{\the\@numparams\space is too few parameters}%
         \else
            \edef\@paramdef{%
               \ifcase\@numparams
                  \empty  No arguments.
               \or ####1%
               \or ####1####2%
               \or ####1####2####3%
               \or ####1####2####3####4%
               \or ####1####2####3####4####5%
               \or ####1####2####3####4####5####6%
               \or ####1####2####3####4####5####6####7%
               \or ####1####2####3####4####5####6####7####8%
               \or ####1####2####3####4####5####6####7####8####9%
               \fi
            }%
         \fi
      \fi
   \fi
   \expandafter\@startdef\@paramdef{#1}%
}%
%
%
%
%
\def\@readauxfile{%
   \if@auxfiledone \else 
      \global\@auxfiledonetrue
      \@testfileexistence{aux}%
      \if@fileexists
         \begingroup
            \endlinechar = -1
            \catcode`@ = 11
            \input \jobname.aux
         \endgroup
      \else
         \message{\@undefinedmessage}%
         \global\@citewarningfalse
      \fi
      \immediate\openout\@auxfile = \jobname.aux
   \fi
}%
%
%
\newif\if@auxfiledone
\ifx\noauxfile\@undefined \else \@auxfiledonetrue\fi
%
%
%
%
\@innernewwrite\@auxfile
\def\@writeaux#1{\ifx\noauxfile\@undefined \write\@auxfile{#1}\fi}%
%
%
%
\ifx\@undefinedmessage\@undefined
   \def\@undefinedmessage{No .aux file; I won't give you warnings about
                          undefined citations.}%
\fi
%
%
\@innernewif\if@citewarning
\ifx\noauxfile\@undefined \@citewarningtrue\fi
%
%
%
\catcode`@ = \@oldatcatcode

\def\pfeilso{\leavevmode
            \vrule width 1pt height9pt depth 0pt\relax
           \vrule width 1pt height8.7pt depth 0pt\relax
           \vrule width 1pt height8.3pt depth 0pt\relax
           \vrule width 1pt height8.0pt depth 0pt\relax
           \vrule width 1pt height7.7pt depth 0pt\relax
            \vrule width 1pt height7.3pt depth 0pt\relax
            \vrule width 1pt height7.0pt depth 0pt\relax
            \vrule width 1pt height6.7pt depth 0pt\relax
            \vrule width 1pt height6.3pt depth 0pt\relax
            \vrule width 1pt height6.0pt depth 0pt\relax
            \vrule width 1pt height5.7pt depth 0pt\relax
            \vrule width 1pt height5.3pt depth 0pt\relax
            \vrule width 1pt height5.0pt depth 0pt\relax
            \vrule width 1pt height4.7pt depth 0pt\relax
            \vrule width 1pt height4.3pt depth 0pt\relax
            \vrule width 1pt height4.0pt depth 0pt\relax
            \vrule width 1pt height3.7pt depth 0pt\relax
            \vrule width 1pt height3.3pt depth 0pt\relax
            \vrule width 1pt height3.0pt depth 0pt\relax
            \vrule width 1pt height2.7pt depth 0pt\relax
            \vrule width 1pt height2.3pt depth 0pt\relax
            \vrule width 1pt height2.0pt depth 0pt\relax
            \vrule width 1pt height1.7pt depth 0pt\relax
            \vrule width 1pt height1.3pt depth 0pt\relax
            \vrule width 1pt height1.0pt depth 0pt\relax
            \vrule width 1pt height0.7pt depth 0pt\relax
            \vrule width 1pt height0.3pt depth 0pt\relax}

\def\pfeilsw{ \leavevmode 
            \vrule width 1pt height0.3pt depth 0pt\relax
            \vrule width 1pt height0.7pt depth 0pt\relax
            \vrule width 1pt height1.0pt depth 0pt\relax
            \vrule width 1pt height1.3pt depth 0pt\relax
            \vrule width 1pt height1.7pt depth 0pt\relax
            \vrule width 1pt height2.0pt depth 0pt\relax
            \vrule width 1pt height2.3pt depth 0pt\relax
            \vrule width 1pt height2.7pt depth 0pt\relax
            \vrule width 1pt height3.0pt depth 0pt\relax
            \vrule width 1pt height3.3pt depth 0pt\relax
            \vrule width 1pt height3.7pt depth 0pt\relax
            \vrule width 1pt height4.0pt depth 0pt\relax
            \vrule width 1pt height4.3pt depth 0pt\relax
            \vrule width 1pt height4.7pt depth 0pt\relax
            \vrule width 1pt height5.0pt depth 0pt\relax
            \vrule width 1pt height5.3pt depth 0pt\relax
            \vrule width 1pt height5.7pt depth 0pt\relax
            \vrule width 1pt height6.0pt depth 0pt\relax
            \vrule width 1pt height6.3pt depth 0pt\relax
            \vrule width 1pt height6.7pt depth 0pt\relax
            \vrule width 1pt height7.0pt depth 0pt\relax
            \vrule width 1pt height7.3pt depth 0pt\relax
            \vrule width 1pt height7.7pt depth 0pt\relax
            \vrule width 1pt height8.0pt depth 0pt\relax
            \vrule width 1pt height8.3pt depth 0pt\relax
            \vrule width 1pt height8.7pt depth 0pt\relax
            \vrule width 1pt height9pt depth 0pt\relax
      }


\def\widestnumber#1#2{}

\def\citewarning#1{\ifx\shlhetal\relax 
    \else
    \par{#1}\par
    \fi
}

\def\rm{\fam0 \tenrm}

\def\fakesubhead#1\endsubhead{\bigskip\noindent{\bf#1}\par}



%
%
%

%

\font\textrsfs=rsfs10
\font\scriptrsfs=rsfs7
\font\scriptscriptrsfs=rsfs5

\newfam\rsfsfam
\textfont\rsfsfam=\textrsfs
\scriptfont\rsfsfam=\scriptrsfs
\scriptscriptfont\rsfsfam=\scriptscriptrsfs

\edef\oldcatcodeofat{\the\catcode`\@}
\catcode`\@11

\def\Cal@@#1{\noaccents@ \fam \rsfsfam #1}

\catcode`\@\oldcatcodeofat


\expandafter\ifx \csname margininit\endcsname \relax\else\margininit\fi

\long\def\red#1\endred{}
\long\def\green#1\endgreen{}
\long\def\blue#1\endblue{}
\long\def\private#1\endprivate{}

\def\endred{ \unmatched endred! }
\def\endgreen{ \unmatched endgreen! }
\def\endblue{ \unmatched endblue! }
\def\endprivate{ \unmatched endprivate! }

\ifx\latexcolors\undefinedcs\def\latexcolors{}\fi

\def\emptycs{}
\def\evaluatelatexcolors{%
        \ifx\latexcolors\emptycs\else
        \expandafter\xxevaluate\latexcolors\xxfertig\evaluatelatexcolors\fi}
\def\xxevaluate#1,#2\xxfertig{\setupthiscolor{#1}%
        \def\latexcolors{#2}}


\font\smallfont=cmsl7
\def\rutgerscolor{\ifmmode\else\endgraf\fi\smallfont
\advance\leftskip0.5cm\relax}
\def\setupthiscolor#1{\edef\tmptmpcs{\noexpand\bgroup\noexpand\rutgerscolor
\noexpand\def\noexpand\currentcolor{#1}%
\noexpand}%
\expandafter\let\csname#1\endcsname\tmptmpcs
\def\tmptmpcs{\checkColorUnmatched{#1}\popthecolor}
\expandafter\let\csname end#1\endcsname\tmptmpcs}

\def\checkColorUnmatched#1{\def\expectcolor{#1}%
    \ifx\expectcolor\currentcolor   
    \else \edef\failhere{\noexpand\tryingToClose '\currentcolor' with end\expectcolor}\failhere\fi}

\def\currentcolor{???}

\def\popthecolor{\ifmmode\else\endgraf\fi\egroup}

\expandafter\def\csname#1\endcsname{}

\evaluatelatexcolors

 \let\outerhead\head
 \def\head{\innerhead}
 \let\innerhead\outerhead

 \let\outersubhead\subhead
 \def\subhead{\innersubhead}
 \let\innersubhead\outersubhead

 \let\outersubsubhead\subsubhead
 \def\subsubhead{\innersubsubhead}
 \let\innersubsubhead\outersubsubhead

 \let\outerproclaim\proclaim
 \def\proclaim{\innerproclaim}
 \let\innerproclaim\outerproclaim

 %
 %
 %
 %

\def\demo#1{\medskip\noindent{\it #1.\/}}
\def\enddemo{\smallskip}

\def\remark#1{\medskip\noindent{\it #1.\/}}
\def\endremark{\smallskip}

\pageheight{8.5truein}
\topmatter
\title{Model theory without choice?  Categoricity} \endtitle

\author {Saharon Shelah \thanks {The author would like to thank the Israel
Science Foundation for partial support of this research (Grant No.242/03)
and Alice Leonhardt for the beautiful typing. Publication 840.}\endthanks}
\endauthor  

\affil{The Hebrew University of Jerusalem \\
Einstein Institute of Mathematics \\
Edmond J. Safra Campus, Givat Ram \\
Jerusalem 91904, Israel
 \medskip
 Department of Mathematics \\
 Hill Center-Busch Campus \\
  Rutgers, The State University of New Jersey \\
 110 Frelinghuysen Road \\
 Piscataway, NJ 08854-8019 USA} 
\endaffil

\abstract  We prove \L os conjecture = Morley theorem in ZF, with
the same characterization (of first order countable theories
categorical in $\aleph_\alpha$ for some (equivalently for every ordinal)
$\alpha > 0$.  Another central result here is, in this context: the
number of models of a countable first order $T$ of cardinality
$\aleph_\alpha$ is either $\ge |\alpha|$ for every 
$\alpha$ or it has a small upper bound (independent of $\alpha$ 
close to $\beth_2$). 
\endabstract
\endtopmatter
\document

\newpage

\head {Annotated Content} \endhead
 \spuriousreset
\bn
\S0 Introduction, pg.3-8
\bn
\ub{Part I}:
\mn
\S1 Morley's proof revisited, pg.9-10
\mr
\item "{${{}}$}"  [We clarify when does Morley proof works - when
there is an $\omega_1$-sequence of reals]
\endroster
\bn
\S2 Stability and categoricity, pg.11-18
\mr
\item "{${{}}$}"  [We prove the choiceless {\L}o\'s conjecture.  This
requires a different proof as possibly there is no well 
ordered uncountable set of reals.  Note that it is harder to construct
non-isomorphic models, as e.g. we do not know whether successor
cardinals are regular and even so whether, e.g. they have a
stationary/co-stationary subset.  Also we have to use more of
stability theory.]
\endroster
\bn
\S3 A dichotomy on $\dot I(\aleph_\alpha,T)$: either bounded or $\ge
|\alpha|$, pg.19-37
\mr
\item "{${{}}$}"  [This shows that ``the lower part of the family of
functions $\{\dot I(\lambda,T):T$ first order complete countable$\}$" is nice.]
\endroster
\bn
\S4 On $T$ categorical in $|T|$, pg.38-43
\mr
\item "{${{}}$}"  [Here we get only partial results.]
\endroster
\bn
\ub{Part II}:
\mn
\S5 Consistency results, pg.44-49
\mr
\item "{${{}}$}"  [We look for cases of ``classes have few models"
which do not occur in the ZFC context.]
\endroster
\bn
\S6 Comments on model theory in ZF, pg.50-53
\bn
\S7 On powers which are not cardinals: categoricity, pg.55-60
\mr
\item "{${{}}$}"  [We deal with models of (first order) theories in so called
reasonable powers (which are not cardinals), it is equivalent to 
the completeness theorem holds.  We throw some light on ``can a 
countable first order $T$ be categorical in some reasonable power".]
\endroster
\bn
References, pg.59
\newpage

\head {\S0 Introduction} \endhead  \resetall \sectno=0
 \spuriousreset
\bigskip

I have known for long that there is no interesting model theory
without (the axiom of) choice, not an exciting question anyhow as we
all know that AC is true.   
This work is dedicated to a try to refute this opinion, i.e., this 
work throws some light on this in the contrary direction: Theorem \scite{0.1}
seriously, Theorem \scite{ss.10}, (the parallel of the ZFC theorem
\scite{0.2R}) in a stronger way. 

Lately, I have continued my work on pcf without full choice (see
\cite{Sh:835}, earlier \cite{Sh:497}, \cite{Sh:E38}, later
\cite{Sh:F728}) and saw that with suitable ``reasonable" weak version
of (the axiom of) choice essentially we can redo 
all \cite{Sh:c} (for first order classes with well ordered vocabulary; see
\scite{4.10}).

Then it seems reasonable to see if older established version suffices,
say ZF + DC$_{\aleph_0}$.  We first consider \L o\'s conjecture which
can be phrased (why only $\aleph_\alpha$'s and not other powers? 
see below)
\bn
\margintag{0.0.7}\ub{\stag{0.0.7} The choiceless \L o\'s Conjecture}:

For a countable (first order theory) $T$: 
\mr
\item "{$(*)_1$}"  $T$ is categorical in $\aleph_\alpha$ for (at
least) one ordinal $\alpha >0$ \nl
iff
\sn
\item "{$(*)_2$}"  $T$ is categorical in $\aleph_\beta$ for every ordinal
$\beta >0$.
\endroster
\bn
In \S1 we shall show that the Morley's proof works exactly when there is an 
uncountable well ordered set of reals.
In \S2 we give a new proof which works always (under ZF); it used
Hrushovski \cite{Hr89d}, so:
\proclaim{\stag{0.1} Theorem}  {\rm (ZF)} For any countable $T$ we
have: $(*)_1$ of \scite{0.0.7} iff $(*)_2$ of \scite{0.0.7}
iff $T$ is $\aleph_0$-stable with no two cardinal models
\footnote{that is, for no model $M$ of $T$, formula $\varphi(x,\bar y)
\in \Bbb L_{\tau(T)}$ and sequence $\bar a \in {}^{\ell g(\bar y)}M$, do we
have $\aleph_0 \le |\varphi(M,\bar a)| < \|M\|$.}.
\endproclaim
\bn
Note that though we have been ready enough to use ZF +DC$_{\aleph_0}$
in fact we solve the problem in ZF.  
\bn
A theorem from \cite{Sh:c} is
\proclaim{\stag{0.2R} Theorem}  {\rm [ZFC]}  For a 
countable complete (first order theory) $T$, one of the following occurs:
\mr
\item "{$(a)$}"  $\dot I(\aleph_\alpha,T) \ge |\alpha|$ for every
ordinal $\alpha$
\sn
\item "{$(b)$}"  $\dot I(\aleph_\alpha,T) \le \beth_2$ for every 
ordinal $\alpha$
(and can analyze this case: either $\dot I(\aleph_\alpha,T) = 1$ for every
$\alpha$ or $\dot I(\aleph_\alpha,T) = \text{\rm Min}
\{\beth_2,2^{\aleph_\alpha}\}$ for every $\alpha >0$).
\endroster
\endproclaim
\bn
We shall prove a similar theorem in ZF in \scite{ss.10}.

Thirdly, we consider an old conjecture from Morley \cite{Mo65}: if a complete 
(first order) $T$ is
categorical in the cardinal $\lambda,\lambda = |T| > \aleph_0$ then
$T$ is a definitional extension of some $T' \subseteq T$ of smaller
cardinality.   The conjecture actually says that $T$ is not really 
of cardinality $\lambda$.
This was proved in ZFC.   Keisler \cite{Ke71a} 
proved it when $|T| < 2^{\aleph_0}$.
By \cite{Sh:4} it holds if $|T|^{\aleph_0} = |T|$.  It is fully proven in
\cite[IX,1.19,pg.491]{Sh:c}).  
The old proof which goes by division to three cases is 
helpful but not sufficient.  Without choice (but note that $\lambda$
is an $\aleph$) the case $T$ superstable
(or just $\kappa_r(T) < \lambda$) has really a similar proof.  The other two
cases, $T$ is unstable and $T$ stable with large $\kappa_r(T)$, are
not.  Here in \S4 it is partially confirmed, e.g., when 
$\lambda$ is regular, the proofs are different though related.

In \S7 we deal with power of non-well orderable sets, in \S5 we deal
with consistency results and in \S6 we look what occurs to classical
theorems of model theory.  

We may consider isomorphism after appropriate forcing.  
Baldwin-Laskowski-Shelah \cite{BLSh:464}, Laskowski-Shelah \cite{LwSh:518}
deal with the question  
``does $T$ or even {\rm PC}$(T_1,T)$ have non-isomorphic models
which become isomorphic after some c.c.c. forcing?"  But this turns
out to be very different and does not seem related to the work here.

However, the following definition \scite{0.3} 
suggests a problem which is closely related but it may be 
easier to find examples of such objects, so called below
``cardinal cases" with ``not so nice behaviour"
\ub{than} to find forcing extension of $\bold V$ which satisfies ZF +
a failure of some hopeful theorem.  
\bigskip

\definition{\stag{0.3} Definition}  1) A cardinal case is a pair 
$(\lambda,\bold P)$ where $\lambda$ is a cardinal and $\bold P$ is a 
family of forcing notions.
\nl
2) A cardinal$^+$ case is a triple $(\lambda,\bold P,<)$ such that
$\lambda$ is a cardinal, $\bold P$ a family of forcing notions and
$<$ a partial order on $\bold P$ such that $\Bbb P_1 \le \Bbb P_2 
\Rightarrow \Bbb P_1 \lessdot \Bbb P_2$ so if we omit $<$ we mean $\lessdot$.
\nl
4) We say that a theory $T$ or more generally a (definition, absolute
enough, of a) class ${\frak K}$ of
models is categorical in the cardinal case $(\lambda,\bold P)$ \ub{when}:
for every $M_1,M_2 \in {\frak K}_\lambda$ (i.e., $\in {\frak K}$ of 
cardinality $\lambda$), for some $\Bbb P
\in \bold P$ we have $\Vdash_{\Bbb P} ``M_1 \cong M_2$".
\nl
5) We say that a theory $T$ or more generally a (definition, absolute
enough, of a) class ${\frak K}$ of models is categorical in the 
cardinal$^+$ case
$(\lambda,\bold P,<)$ \ub{when} for any $\Bbb P \in \bold P$, in $\bold
V^{\Bbb P}$ we have: if
$M_1,M_2\in {\frak K}^{\bold V[\Bbb P]}_\lambda$, then for some $\Bbb P' \in 
\bold P$ satisfying $\Bbb P \le \Bbb P'$ we have $\Vdash_{\Bbb P'/\Bbb
P} ``M_1 \cong M_2"$.
\nl
4) Similarly uncategorical, has/does not have $\mu$ pairwise 
non-isomorphic models, etc.
\nl
5) We may replace cardinal by power.
\enddefinition
\bn
\margintag{0.4N}\ub{\stag{0.4N} Question}:  Characterize countable (complete first
order) $T$ which may be categorical in
some uncountable power (say in some forcing extension of $\bold
L[T]$).  See on this \S7.  

This work may be continued in \cite{Sh:F701}.

We thank Udi Hrushovski for various comments and pointing out that 
\L o\'s conjecture proof is over after \scite{uni.3} as Kueker conjecture
is known in the relevant case (in earlier versions the proof (of the
choiceless \L o\'s conjecture) was more interesting and longer).  We thank Moti
Gitik for a discussion of the consistency results and for pointing out
\scite{4.17.1}). 

Lately, I have learned that Truss
and his students were pursuing the connection between universes with restricted
choice and model theory by a different guiding line: using model
theory to throw light on the arithmetic of Dedekind finite powers, 
works in this direction are Agatha Walczak-Typke \cite{WT05},
\cite{WT07}.  Very interesting, does not interact with the present
investigation, but may be relevant to Question \scite{0.4N}.

\bn
\centerline{$* \qquad * \qquad*$}
\bn
Recall
\definition{\stag{0.A} Definition}  A cardinal is the power of some
well ordered set (so an $\aleph$ or a natural number).
\enddefinition
\bn
In \cite{Sh:F701} we may deal with theories in a vocabulary which 
is not well ordered.
\sn 
\margintag{0.B}\ub{\stag{0.B} Convention}:  If not said otherwise
\mr
\item "{$(a)$}"  $T$ is a first order theory in a vocabulary $\tau
\subseteq \bold L$
\sn
\item "{$(b)$}"  $T$ is complete
\sn
\item "{$(c)$}"  $T$ is infinite
\sn
\item "{$(d)$}"  if $T$ is countable for simplicity $\tau,T \subseteq
{\Cal H}(\aleph_0)$ (for notational simplicity).
\ermn
This is justified by
\demo{\stag{0.C} Observation}  Assume $\tau$ is a countable
vocabulary and $T$ is a first order theory in $\tau$, i.e., $T \subseteq
\Bbb L_\tau$.
\nl
1) There is a vocabulary $\tau' \subseteq 
{\Cal H}(\aleph_0)\,(\subseteq \bold L_\omega)$ and first order theory
$T'$ in $\tau'$ (so $T' \subseteq \Bbb L_\tau \subseteq {\Cal H}(\aleph_0))$
such that for every cardinal $\lambda,T$ is categorical in $\lambda$
iff $T'$ is categorical in $\lambda$ (and even $\dot I(\lambda,T) =
\dot I(\lambda,T')$, similarly for power and the parallel of
\scite{0.C.1} below).
\nl
2) If $T$ is categorical in some cardinal $\lambda$ then $T \cup
\{(\exists^{\ge n} x)(x = x):n < \omega\}$ is complete.
\enddemo
\bigskip

\demo{\stag{0.C.1} Observation}  Assume $\tau$ is a vocabulary which
 can be well ordered (i.e., $|\tau| \in \text{ Card}$).
\nl
There is a vocabulary $\tau' \in \bold L$ (or even $\tau' \in \bold
L_{|\tau|^+}$) and a function $f$ from $\Bbb L(\tau')$ onto 
$\Bbb L(\tau)$ (note that $\Bbb L(\tau') \subseteq \bold L_{|\tau|^+})$
mapping predicates/functions symbols to predicate/function symbols
respectively with the same arity such that:
\mr
\item "{$\boxtimes_1$}"  $f$ maps the set of (complete) first order
theories in $\Bbb L_\tau$ onto the set of (complete) first order
theories in $\Bbb L_{\tau'}$ (really this is a derived map, $\hat f$)
\sn
\item "{$\boxtimes_2$}"  for some definable class $\bold F$ which is 
a function, $\bold F$ maps the class of $\tau$-models
onto the class of $\tau'$-models such that
{\roster
\itemitem{ $(a)$ }  $\bold F$ is one to one onto and Th$(\bold F(M))
= \hat f(\text{Th}(M))$
\sn
\itemitem{ $(b)$ }  $\bold F$ preserves isomorphisms and non-isomorphisms
\sn
\itemitem{ $(c)$ }  $\bold F$ preserves $M \subseteq N,M \prec N$
\sn
\itemitem{ $(d)$ }  for some function $f$, if $\bold F(M) = M'$ 
then for every sentence $\psi
\in \Bbb L(\tau'),M' \models \psi \Leftrightarrow M \models f(\psi)$
where $f(\psi)$ is defined naturally
\sn
\itemitem{ $(e)$ }   $\bold F$ preserves power, so equality and inequality of
powers (hence for any theory $T \subseteq \Bbb L_\tau$, letting $T' =
\hat f(T)$, for any set $X,(\{M/\cong:M \in \text{\rm Mod}_T$ has
power $|X|\}| = |\{M'/\cong:M' \in \text{\rm Mod}_{T'}$ has power $|X|\}|$.
\endroster}
\endroster
\enddemo
\bn
We shall use absoluteness freely recalling the main variant.
\definition{\stag{0.D} Definition}  1) We say $\varphi(\bar x)$ is
upward ZFC-absolute \ub{when}: 
if $\bold V_1 \subseteq \bold V_2$ (are transitive classes
containing the class Ord of ordinals, both models of ZFC) and 
$\bar a \in \bold V_1$ \ub{then}
$\bold V_1 \models \varphi(\bar a) \Rightarrow \bold V_2 \models
\varphi(\bar a)$.
\nl
2) Replacing upward by downward mean we use $\Leftarrow$; omitting
upward mean we use $\Leftrightarrow$. 
Similarly for version ZFC$'$ of ZFC (e.g. ZF + DC); 
but absolute means ZFC-absolute.
\enddefinition
\bigskip

\demo{\stag{0.E} Convention}  1) If not said otherwise, for a theory $T$
belonging to $\bold L[Y_0],Y_0 \subseteq \text{ Ord}$, saying  
``$T$ satisfies Pr", (``Pr" stands for ``Property") 
we mean ``for some $Y_1 \subseteq
\text{ Ord}$ for every $Y_2 \subseteq \text{ Ord},T$ satisfies Pr in $\bold
L[Y_0,Y_1,Y_2]$". 
\nl
2) But ``$T$ categorical in $\lambda$" always means in $\bold V$.
\enddemo
\bn
Recall
\definition{\stag{0.EA} Definition}  1) $\theta(A) = \text{ Min}
\{\alpha$: there is no function from $A$ onto $\alpha\}$.
\nl
2) $\Upsilon(A) = \text{ Min}\{\alpha$: there is no one-to-one function from
$\alpha$ into $A\}$.
\enddefinition
\bigskip

\definition{\stag{0.F} Definition}  1) If $T \subseteq \Bbb
L(\tau),\Gamma$ is a set of types in $\Bbb L(\tau)$, i.e., each is an
$m$-type for some $m$, \ub{then} EC$(T,\Gamma)$
is the class of $\tau$-models $M$ of $T$ which omits every $p(\bar x)
\in \Gamma$.
\nl
2) If $T \subseteq \Bbb L(\tau)$ is complete, $T \subseteq T_1
\subseteq \Bbb L(\tau_1)$ and $\tau \subseteq \tau_1$ \ub{then} 
PC$(T_1,T)$ is the
class of $\tau$-reducts of models $M_1$ of $T_1$.  Similarly for a set
$\Gamma$ of types in $\Bbb L(\tau_1)$ let PC$(T,T_1,\Gamma)$ be the
class of $\tau$-reducts of models $M \in \text{ EC}(T_1,\Gamma)$.
\enddefinition
\bn
We shall use Ehrenfuecht-Mostowski models.
\definition{\stag{0.F13} Definition}  1) $\Phi$ is proper for linear
orders \ub{when}:
\mr
\item "{$(a)$}"  for some vocabulary $\tau = \tau_\Phi =
\tau(\Phi),\Phi$ is an $\omega$-sequence, the $n$-th element a complete
quantifier free $n$-type in the vocabulary $\tau$
\sn
\item "{$(b)$}"  for every linear order $I$ there is a $\tau$-model
$M$ denoted by EM$(I,\Phi)$, generated by $\{a_t:t \in I\}$ such that
$s \ne t \Rightarrow a_s \ne a_t$ for $s,t \in I$ and
 $\langle a_{t_0},\dotsc,a_{t_{n-1}}\rangle$ realizes the quantifier free
$n$-type from clause (a) whenever $n < \omega$ and
$t_0 <_I \ldots <_I t_{n-1}$; so
really $M$ is determined only up to isomorphism but we may ignore this and
use $I_1 \subseteq J_1 \Rightarrow \text{ EM}(I_1,\Phi) \subseteq
\text{ EM}(I_2,\Phi)$.  We call $\langle a_t:t \in I\rangle$ ``the"
skeleton of $M$; of course ``the" is an abuse of notation as it is not
necessarily unique.
\ermn
2) If $\tau \subseteq \tau(\Phi)$ then we let EM$_\tau(I,\Phi)$ be the
$\tau$-reduct of EM$(I,\Phi)$.
\nl
3) For first order $T$, let 
$\Upsilon^{\text{or}}_\kappa[T]$ be the class of $\Phi$ proper
for linear orders such that
\mr
\item "{$(a)$}"  $\tau_T \subseteq \tau_\Phi$ and
$\tau_\Phi$ has cardinality $\le \kappa$
\sn
\item "{$(b)$}"  for any linear order $I$ the model 
{\rm EM}$(I,\Phi)$ has cardinality $|\tau(\Phi)|+ |I|$ and 
we have {\rm EM}$_{\tau(T)}(I,\Phi) \in K$
\sn
\item "{$(c)$}"  for any linear orders $I \subseteq J$ we have 
{\rm EM}$_{\tau(T)}(I,\Phi) \prec$ {\rm EM}$_{\tau(T)}(J,\Phi)$.
\ermn
4) We may use Skeleton $\langle \bar a_t:t \in I\rangle$ with $\alpha
= \ell g(\bar a_t)$ constant but in the definition of ``$\Phi \in
\Upsilon^{\text{or}}_\kappa[T]$; we add $\alpha < \kappa^+$.  
Alternatively $\bar a_t = \langle
F_i^{\text{EM}(I,\Phi)}(a_t):i < \alpha\rangle$, where $F_i \in
\tau_\Phi$ are unary function symbols.  We use $\Phi,\Psi$ only for such
objects.  Let $\Upsilon^{\text{or}}_T =
\Upsilon^{\text{or}}_{|T|+\aleph_0}[T]$. 
\enddefinition
\newpage

\head {\S1 Morley's proof revisited} \endhead  \resetall \sectno=1
 \spuriousreset
\bigskip

The main theorem of this section is \scite{1.1}.  The proof is just
adapting Morley's proof in ZFC.  We shall use \scite{0.C}(2) and
convention \scite{0.B} freely.
\proclaim{\stag{1.1} Theorem}  [{\rm ZF} + there is an uncountable well
ordered set of reals].
\nl
The following conditions on a countable (first order) $T$ are equivalent:
\mr
\item "{$(A)$}"  $T$ is categorical in some cardinal $\aleph_\alpha >
\aleph_0$, in $\bold V$, of course
\sn
\item "{$(B)$}"  $T$ is categorical in every cardinal $\aleph_\beta >
\aleph_0$, in $\bold V$, of course
\sn
\item "{$(C)$}"  $T$ is (in $\bold L[T]$), totally transcendental
(i.e. $\aleph_0$-stable) with no two cardinal models (i.e., for
no model $M$ of $T$ and formula 
$\varphi(x,\bar y) \in \Bbb L(\tau_T)$ and
$\bar a \in {}^{\ell g(\bar y)}M$ do we have $\aleph_0 \le
|\varphi(M,\bar a)| < \|M\|$ and $\|M\|$ is a cardinal, i.e. the set
of elements of $M$ is well-orderable hence its power is a cardinal)
\sn
\item "{$(D)$}"  if $\bold V' \subseteq \bold V$ is 
a transitive class extending $\bold L,T \in \bold V'$ and $\bold V'$
satisfies {\rm ZFC} then the conditions in (C) hold
\sn
\item "{$(E)$}"  for some $\bold V'$ clause (D) holds.
\endroster
\endproclaim
\bigskip

\demo{Proof}  By \scite{0.B} or better \scite{0.C}(2) 
\wilog \, $T$ is complete, $T \subseteq
{\Cal H}(\aleph_0)$.   Trivially $(B) \Rightarrow (A)$.  Next $(A) \Rightarrow
(C)$ by claims \scite{1.2}, \scite{1.3} below.  Lastly, $(C)
\Rightarrow (A)$ by \scite{1.4} below and $(C) \Leftrightarrow (D)
\Leftrightarrow (E)$ holds by absoluteness.  \hfill$\square_{\scite{1.1}}$
\enddemo
\bigskip

\proclaim{\stag{1.2} Claim}  [{\rm ZF} + $\exists$ a set of $\aleph_1$ reals]

If $T$ is (countable) and in $\bold L[T]$ the theory $T$ is not 
$\aleph_0$-stable and $\lambda >
\aleph_0$ \ub{then} $T$ is not categorical in $\lambda$.
\endproclaim
\bigskip

\demo{Proof}  In $\bold L[T]$ we can find E.M. models, i.e. $\Phi \in
\Upsilon^{\text{or}}_T$ such that
$\tau(\Phi)$ is countable, extends $\tau = \tau_T$ and EM$_\tau(I,\Phi)$ is a
model of $T$ (of cardinality $\lambda$) for every 
linear order $I$ (of cardinality $\lambda$) and let $M_1 = \text{
EM}_\tau(\lambda,\Phi)$ and \wilog \, the universe of $M_1$ is
$\lambda$.

In $\bold V$ let 
$\bar \eta = \langle \eta_\alpha:\alpha < \omega_1 \rangle$ be a
sequence of pairwise distinct reals.  In $\bold L[T]$ there is a
countable model $M_0$ of $T$ with ${\bold S}(M_0)$ uncountable so
containing a perfect set.  Hence also in $\bold L[T,\bar \eta],M_0$ is
a countable model of $\tau$ with ${\bold S}(M_0)$ containing a perfect
set, hence there is (in $\bold L[T,\bar \eta]$) a model $M_2$ of $T$
of cardinality $\lambda$ ($\lambda$ is still an uncountable cardinal
in $\bold L[T,\bar \eta]$) such that $M_0 \prec M_2$ and there is a
sequence $\langle a_i:i < \omega_1 \rangle,a_i \in M_2$ realizes $p_i
\in {\bold S}(M_0)$ with $\langle p_i:i < \omega^{\bold V}_1 \rangle$ pairwise 
distinct.  Without loss of generality the universe of $M_2$ is $\lambda$.

Clearly even in $\bold V$, the model $M_1$ satisfies ``if $A \subseteq
M_1$ is countable then the set $\{\text{tp}(a,A,M_1):a \in M_1\}$ is
countable" whereas $M_2$ fails this;
hence the models $M_1,M_2$ have universe $\lambda$,
are models of $T$ and are not isomorphic, so we are done.
\hfill$\square_{\scite{1.2}}$ 
\enddemo
\bigskip

\proclaim{\stag{1.3} Claim}  Assume $T$ is countable $\aleph_0$-stable and has a
two cardinal model (in $\bold L[T]$, but both are absolute).

\ub{Then}  $T$ is not categorical in $\lambda$, in fact,
$\dot I(\aleph_\alpha,T) \ge |\alpha|$ for every ordinal $\alpha$.
\endproclaim
\bigskip

\demo{Proof}  So in $\bold L(T)$ it has a model $M_1$ and a finite sequence
$\bar a \in {}^{\ell g(\bar y)}(M_1)$ and a formula 
$\varphi(x,\bar y) \in \Bbb L(\tau_T)$
such that $\aleph_0 \le |\varphi(M_1,\bar a)| < \|M_1\|$.  
If $\aleph_\beta < \lambda$, working in
$\bold L[T]$ \wilog \,  $|\varphi(M_1,\bar a)| = \aleph_\beta,\|M_1\| =
\lambda$ (by \cite{Sh:3}) and the universe of $M_1$ is 
$\lambda$.  But $\aleph_0$-stability $T$ has (in
$\bold L[T]$) a saturated model $M_2$ of cardinality $\lambda$, so
\wilog \, the universe of $M_2$ is $\lambda$.  So $\bar a' \in
{}^{\ell g(\bar y)}(M_2) \Rightarrow |\varphi(M_2,\bar a')| \notin
[\aleph_0,\lambda)$.  Clearly even in $\bold V,M_1,M_2$ are models of
$T$ of cardinality $\lambda$ and are not isomorphic.  In fact for
every $\aleph_\beta \le \lambda,T$ has an $\aleph_\beta$-saturated not
$\aleph_{\beta +1}$-saturated model $M_\beta$ of cardinality $\lambda$
such that $|\varphi(M_\beta,\bar a_\beta)| = \aleph_\beta$ for some
$\bar a \in {}^{\ell g(\bar y)}(M_\beta)$.  So
$\{|\varphi((M_\beta,\bar b)| + \aleph_0:\bar b \in {}^{\ell g(\bar
y)}(M_\beta)|\} = \{\aleph_\beta,\lambda\}$, hence $\bold V \models
M_\beta \ncong M_\gamma$ when $\aleph_\beta < \aleph_\gamma \le
\lambda$ so also the second phrase in the conclusion of the claim
holds and even $\dot I(\aleph_\alpha,T) \ge |\alpha +1|$.
\hfill$\square_{\scite{1.3}}$
\enddemo
\bigskip

\proclaim{\stag{1.4} Claim}  Assume $T$ is countable $\aleph_0$-stable with no
two cardinal models even just in $\bold L[T]$ and 
$\lambda > \aleph_0$.  \ub{Then} $T$ is
categorical in $\lambda$.
\endproclaim
\bigskip

\demo{Proof}  Let $M_1,M_2$ be models of $T$ of cardinality $\lambda$,
\wilog \, both have universe $\lambda$, clearly $\bold L[T,M_1,M_2]$ is
a model of ZFC and by absoluteness $T$ still satisfies the assumption
of \scite{1.4} in it, and $M_1,M_2$ are (also in it) uncountable
models of $T$ of the same uncountable cardinality in this universe.
But by \scite{1.1} being a theorem of ZFC clearly $M_1,M_2$ are isomorphic
in $\bold L[T,M_1,M_2]$, hence in $\bold V$.
\hfill$\square_{\scite{1.4}}$
\enddemo
\newpage

\head {\S2 Stability and categoricity} \endhead  \resetall \sectno=2
 \spuriousreset
\bigskip

Our aim in this section 
is the categoricity spectrum for countable $T$ (i.e. Th. \scite{st.0}),
but in the
claims leading to the proof we do not assume countability.  
Note that the absoluteness of various
properties is easier for countable $T$.
\bigskip

\proclaim{\stag{st.0} Theorem}  {\rm [ZF]} For countable $T$, clauses
(A),(B),(C),(D),(E) of Theorem \scite{1.1} are equivalent.
\endproclaim
\bigskip

\demo{\stag{st.1} Observation}  1) If $T$ is unstable so has the
order property, say as witnessed
by $\varphi(\bar x,\bar y)$ and, of course, $\tau = \tau_T \subseteq \bold L$
\ub{then} for some $\Phi \in \bold L[T]$
\mr
\item "{$\circledast$}"  $(a) \quad \Phi$ is proper for linear orders
\sn
\item "{${{}}$}"  $(b) \quad \tau \subseteq \tau_\Phi$ and
for every linear order $I$, EM$_\tau(I,\Phi)$ is a model of $T$  
\nl

\hskip25pt  with
skeleton $\langle \bar a_t:t \in I \rangle,\ell g(\bar a_\tau) = \ell
g(\bar x) = \ell g(\bar y)$
\sn
\item "{${{}}$}"  $(c) \quad$ EM$_\tau(I,\Phi) \models \varphi[\bar
a_s,\bar a_t]^{\text{if}(s < t)}$
\sn
\item "{${{}}$}"  $(d) \quad \tau(\Phi) \subseteq \bold L$ and $|\tau(\Phi)|
= |T|$ (if $\tau(T) \in \bold L$, \wilog 
\nl

\hskip25pt  $\tau(\Phi) \in \bold L$) and \wilog \, $\bold L[T]
\models |\tau(\Phi)| = |T|$.
\ermn
2) It follows that if $I$ is well orderable \ub{then} the universe of 
EM$(I,\Phi)$ is well orderable so it is of cardinality $|I| + |T| +
\aleph_0$ hence we can assume it has this cardinal as its universe.
\enddemo
\bigskip

\demo{Proof}  1) By \cite{Sh:c}.
\nl
2) Follows.  \hfill$\square_{\scite{st.1}}$
\enddemo
\bn
Our first aim is to derive stability from categoricity, for diversion
we give some versions.
\proclaim{\stag{st.2} Claim}  Let $\Phi$ be as in \scite{st.1}.  Then

$M_1 \ncong M_2$ \ub{when}
$\kappa_1,\kappa_2$ are regular uncountable cardinals $(> |T|)$ and for some $A
\subseteq { \text{\rm Ord\/}}$, in $\bold L[A]$
\mr
\item "{$\circledast$}"  $(a) \quad M_\ell = { \text{\rm
EM\/}}_\tau(I_\ell,\Phi)$ in $\bold L[A]$, (so $T,\Phi \in \bold
L[A],I_\ell \in \bold L[A]$) for $\ell=1,2$
\sn
\item "{${{}}$}"  $(b) \quad \bar s^1 = \langle s^1_\alpha:\alpha <
\kappa_1 \rangle$ is increasing in $I_1,\bar t^1 = \langle
t^1_\alpha:\alpha < \kappa_2 \rangle$ is \nl

\hskip25pt decreasing in $I_1$ (in $\bold L[A]$)
\sn
\item "{${{}}$}"  $(c) \quad \alpha < \kappa_1 \wedge \beta < \kappa_2
\Rightarrow s^1_\alpha <_{I_1} t^1_\beta$ but $\neg(\exists s \in I_1)
[(\forall \alpha < \kappa_1)(s^1_\alpha < r) \wedge$ \nl

\hskip25pt $(\forall \beta < \kappa_2)(s < t^1_\beta)]$
\sn
\item "{${{}}$}"  $(d) \quad$ in $I_2$ there is no pair of sequences
like $\bar s^1,\bar t^1$
\sn
\item "{${{}}$}"  $(e) \quad$ also in the inverse of $I_2$, there is
no such pair
\sn
\item "{${{}}$}"  $(f) \quad$ [only for simplicity, implies (d)+(e)] $I_2$ is 
$\cong I_2 \times \Bbb Q$ ordered
\nl

\hskip25pt  lexicographically.
\endroster
\endproclaim
\bigskip

\demo{Proof}  Without loss of generality \, the universes of $M_1,M_2$
are ordinals, and toward contradiction assume $f$ is an isomorphism
from $M_1$ onto $M_2$.  We can work in $\bold L[A,M_1,M_2,f]$ which is
a model of ZFC, so easy to contradict (as in \cite{Sh:12}, see
detailed proof showing more in \scite{ss.1.1}). \nl
${{}}$  \hfill$\square_{\scite{st.2}}$ 
\enddemo
\bigskip

\demo{\stag{st.3} Conclusion}  [ZF + $|T|^+$ is regular]  If $T$ is
categorical in some cardinal $\lambda > |T|$, \ub{then} $T$ is
stable (in $\bold L[T]$).
\enddemo
\bn
A fuller version is
\proclaim{\stag{st.5} Claim}  $M_1 \ncong M_2$ when for some $\lambda
>|T|$ we have:
\mr
\item "{$\circledast$}"  $(a) \quad M_\ell = { \text{\rm
EM\/}}_\tau(I_\ell,\Phi)$ where $T,\Phi$ are as in \scite{st.1} 
so $\Phi \in \bold L[T]$
\sn
\item "{${{}}$}"  $(b) \quad I_1 = \lambda \times \Bbb Q$ ordered
lexicographically 
\sn
\item "{${{}}$}"  $(c) \quad I_2 = \dsize \sum_{\alpha \le \lambda}
I^2_\alpha \in \bold L[T]$ where $I^2_\alpha$ is isomorphic to 
$\alpha + \alpha^*$ ($\alpha^*$ the inverse of $\alpha$)
\nl

\hskip25pt or just
\sn
\item "{${{}}$}"  $(c)^- \quad I_2$ is a linear order of cardinality
$\lambda$ such that for every limit ordinal 
\nl

\hskip25pt $\delta \le |T|^+,I_2$ has an
interval isomorphic to $\delta + \delta^*$
\sn
\item "{${{}}$}"  $(d) \quad I_1,I_2$ has cardinality $\lambda$.
\endroster
\endproclaim
\bigskip

\demo{Proof}  Let $\theta = |T|$ in $\bold L[T]$ and $\theta_1 =
(\theta^+)^{\bold V}$.
Without loss of generality $M_\ell$ has universe
$\lambda$, assume toward contradiction that $M_1 \cong M_2$ let $f$ be
an isomorphism from $M_1$ onto $M_2$ and consider the universe $\bold
L[T,M_1,M_2,f]$.  In this universe $\theta^{\bold V}_1$ may be singular
but is still a cardinal so $\delta =: (\theta^+)^{\bold L[T,M_1,M_2,f]}$
is necessarily $\le (|\theta|^+)^{\bold V}$ hence $I_2$ has an interval
isomorphic to $\delta + \delta^*$.  Now we continue as in
\scite{st.2} (see details in \scite{ss.1.1}).  \hfill$\square_{\scite{st.5}}$
\enddemo
\bigskip

\demo{\stag{st.6} Conclusion}  If $T$ is categorical in the cardinal
$\lambda > |T|$, \ub{then} $T$ is stable.
\enddemo
\bn
\margintag{st.7}\ub{\stag{st.7} Discussion}:  1) We may like to have many models.  So
for $T$ unstable if there are $\alpha$ regular cardinals $\le \lambda$ we can
get a set of pairwise non-isomorphic models of $T$ of cardinality
$\lambda$ indexed by $|{\Cal P}(\alpha)|$.

It is not clear what, e.g., we can get in $\aleph_1$.  As \scite{st.5}
indicate it is hard to have few models, i.e., to have such universe
(see more in \S3); but for our present purpose all this is peripheral,
as we have gotten two.  
\bn
On uni-dimensional see \cite[V,Definition 2.2,pg.241]{Sh:c} 
and \cite[V.Theorem 2.10,p.246]{Sh:c}.
\definition{\stag{uni.0} Definition}  
A stable theory $T$ is uni-dimensional \ub{if} there are no $M \models
T$ and two infinite indiscernible sets in $M$ which are orthogonal.
\enddefinition
\bigskip

\proclaim{\stag{uni.1} Claim}  Assume $T$ is stable 
(in $\bold L[T]$, anyhow this is $Z^-$-absolute).  
\ub{Then} for every $\lambda > |T|,T$ has a 
model $M \in \bold L[T]$ of cardinality $\lambda$ such that:
\mr
\item "{$\odot$}"  in $M$ there are no two (infinite) indiscernible
non-trivial sets each of cardinality $\ge |T|^+$ which are orthogonal.
\endroster
\endproclaim
\bigskip

\demo{Proof}  We work in $\bold L[T]$ or $\bold L[T,Y],Y \subseteq
\text{ Ord}$ and let $\kappa = |T|^{\bold L[T]}$ and $\partial = 
\theta^{\bold V}({\Cal P}(\kappa))$.  Let $\mu$ be large enough
(e.g., $\beth((2^\partial)^+)$, i.e. the $(2^\partial)^+$-th beth), 
let ${\frak C}$ be a $\mu^+$-saturated model
of $T$.  Let $\bold I = \{a_i:i < \mu\} \subseteq {\frak C}$ be an infinite
indiscernible set of cardinality $\mu$ and minimal, i.e.
Av$(\bold I,\cup \bold I)$ is a minimal type.  
Let $M_1 \prec {\frak C}$ be $\kappa^+$-prime over $\bold I$. 

More specifically
\mr
\item "{$\circledast$}"  $\langle A_\varepsilon:\varepsilon \le
\kappa^+\rangle$ is an increasing sequence of subsets of
$M_1,A_{\kappa^+} = M_1,A_0 = \{a_i:i < \mu\}$ and $\bar B =
\langle B_a:a \in M_1\rangle$, satisfies $[a \in A_0
\Rightarrow B_a = \{a\}]$ and if $a \in A_{\varepsilon +1} \backslash
A_\varepsilon$ then $B_a \subseteq A_\varepsilon$ and
tp$(a,B_a,{\frak C}) \vdash \text{ tp}(a,A_{\varepsilon +i} \backslash
\{a\},{\frak C})$ and $|B_a|
\le \kappa$ and \wilog \, $B_a = \{b_{a,j}:j < \kappa\}$ and $\zeta <
\varepsilon \Rightarrow B_a \cap A_{\zeta +1} \nsubseteq A_\zeta$ and
$a' \in B_a \Rightarrow B_{a'} \subseteq B_a$.
\ermn
Expand $M_1$ to
$M_2$ by adding $P^{M_2} = \bold I,<^{M_2} = \{(a_i,a_j):i < j <
\mu\},E^{M_2} = \{(b_1,b_2)$: for some $\varepsilon < \kappa^+$ we have
$b_1 \in A_\varepsilon \wedge b_2 \in (A_{\varepsilon +1} \backslash
A_\varepsilon)\},F^{M_2}_j(a) = b_{a,j}$ for $j < \kappa^+$, (hence $a \in
A_{\varepsilon +1} \backslash A_\varepsilon \wedge \zeta <
\varepsilon \Rightarrow c \ell_{M_2}\{a\} \cap A_{\zeta +1}
\backslash A_\zeta \ne \emptyset$) and add
Skolem functions, still $\tau(M_2) \in \bold L[T]$ 
has cardinality $\kappa = (|T| + \aleph_0)^{\bold L[T]}$.

Now (as in the proof of the omitting type theorem, see e.g.,
\cite[VII,\S5]{Sh:c}) we can find $\langle \bold I_n:n < \omega
\rangle,\bold I_n \subseteq  \bold I$ is an 
$n$-indiscernible sequence in $M_2$ of cardinality $> 2^\partial$ and
\mr
\item "{$(*)_1$}"   if for $n < \omega,M_2 \models a^n_0 < \ldots
< a^n_{n-1}$ and $\ell < n \Rightarrow a^n_\ell \in \bold I_n$ then 
$p_n = \text{ tp}(\langle a^n_0,\dotsc,a^n_{n-1} \rangle,\emptyset,M_2)
= \text{ tp}(\langle a^{n+1}_0,\dotsc,a^{n+1}_{n-1}
\rangle,\emptyset,M_2)$.
\ermn
Let $\bold I_n= \{a_\alpha:\alpha \in {\Cal U}_n\}$ and note that
\mr
\item "{$\boxtimes_1$}"  $\langle 
\bar\sigma(a_{i(\alpha,0)},\dotsc,a_{i(\alpha,m-1)}):\alpha \in Z\rangle$
is an indiscernible set in $M_1$ (equivalently in ${\frak C}$) when:
{\roster
\itemitem{ $(a)$ }  $2m \le n < \omega$
\sn
\itemitem{ $(b)$ }  $Z \subseteq {\Cal U}_m$ is infinite
\sn
\itemitem{ $(c)$ }  $i(\alpha,\ell) \in {\Cal U}_n$ is increasing with
$\ell < m$ for $\alpha \in Z$
\sn
\itemitem{ $(d)$ }  for each $\ell,k <m$ and $\alpha_1 <
\beta_1,\alpha_2 < \beta_2$ from $Z$ we have $i(\alpha_1,\ell) <
i(\beta_1,k) \Leftrightarrow i(\alpha_2,\ell) < i(\beta_2,k)$
\sn
\itemitem{ $(e)$ }  $\bar\sigma(x_0,\dotsc,x_{m-1})$ is a finite
sequence of $\tau(M_2)$-term
\sn
\itemitem{ $(f)$ }  all
$\langle a_{i(\alpha,0)},\dotsc,a_{i(\alpha,m-1)}\rangle$ for $\alpha \in Z$
realize the same type (equivalently of quantifier-free type) in $M_2$.
\endroster}
\ermn
[Why?  If $k\le n$ and $j < \ell g(\bar \sigma)$ then the truth values of
$\sigma_j(a_{i_0},\dotsc,a_{i_{k-1}}) \in A_\varepsilon$ for $i_0 <
\ldots < i_{k-1} < \mu$ such that $a_{i_0},\dotsc,a_{i_{k-1}} \in
\bold I_n$ depend on $\sigma$ only, we can prove this by induction on
max$\{\varepsilon_j:j < \ell g(\bar \sigma)\}$, using  
the properties of the $B_a$'s.  By the properties
of $\bold F^t_{\kappa^+}$-constructions \footnote{instead we
can use the conclusion derived in $\boxtimes_2$}
(\cite[IV]{Sh:c}) we are easily done\footnote{this is not the end of
the proof, we still need to show another indiscernible set does not exist}.]

Moreover
\mr
\item "{$\boxtimes_2$}"  in $\boxtimes_1$, the $\Bbb L(\tau_T)$-type
of $\langle \bar\sigma(a_{i(\alpha,0)},\dotsc,a_{i(\alpha,m-1)}):\alpha
\in Z\rangle$ depends
just on $Z,\bar\sigma$ and the truth values in (d) from $\boxtimes_1$ and
the types in (f) of $\boxtimes_1$ over acl$_{M_1}(\emptyset)$.
\ermn
[Why:  Note that acl$_{M_2}(\emptyset) \subseteq$ (the Skolem hull of
$\emptyset$ in $M_2$) and $\bold I \cap \text{ acl}_{M_2}(\emptyset)$
is infinite.]

So we can find a $\tau(M_2)$-model $M_3$ generated by the
indiscernible sequence $\langle b_\alpha:\alpha < \lambda \rangle$
such that for every $n <\omega$ and $\alpha_0 < \ldots < \alpha_{n-1} 
< \lambda$, recalling $(*)_1$ we have
$p_n = \text{ tp}(\langle b_{\alpha_0},\dotsc,b_{\alpha_{n-1}}\rangle,
\emptyset,M_3)$.  Without loss of generality the Skolem hull of
$\emptyset$ in $M_3$ is the same as in $M_2$.  Let $M_4 = M_3
\restriction \tau_T$.  Clearly $M_4$ is a model of $T$ of 
cardinality $\lambda$.
 
Now suppose that
\mr
\item "{$(*)_2$}"   in $\bold V$ we have $\bold J \subseteq M_4$ (or
even $\bold J \subseteq {}^{\omega \ge}(M_4))$, is an
indiscernible set of cardinality $\ge \kappa^+$ orthogonal to 
$P^{M_3}$ which is an infinite indiscernible set in $M_4$
(this is absolute enough).
\ermn
Let $\bold J \supseteq \{c_\alpha:\alpha < (\kappa^+)^{\bold V}\}$ with 
the $c_\alpha$'s pairwise distinct.  Let $c_i =
\sigma_i(b_{\alpha(i,0)},\dotsc,b_{\alpha(i),n(\alpha)-1})$ where $\alpha(i,0)
< \ldots < \alpha(i,n(i))$ (may be clearer in $\bold L[T,Y,\bold J]$).

So in $\bold L[T,Y,\bold J]$ for some $Z \subseteq (\kappa^+)^{\bold
V}$ of cardinality $\ge (\kappa^+)^{\bold L[T,Y,\bold J]}$ 
(so maybe $\bold V \models |Z| < \kappa^+$) 
we have $i \in Z \Rightarrow \sigma_i =
\sigma_* \wedge n(i) = n(*)$, and the truth
value of $\alpha(i_1,\ell_1) < \alpha(i_2,\ell_2)$ for $i_1 < i_2$
depend just on $(\ell_1,\ell_2)$.
In $\bold L[T,Y]$ for each $n
\ge 2n(*)$, we can find $a^n_{i,\ell} \in \bold I_n$ for $i < \partial,\ell
< n(*)$ such that $M_2 \models a^n_{i_1,\ell_1} < a^n_{i_2,\ell_2}
\Leftrightarrow \alpha(0,\ell_1) < \alpha(1,\ell_2)$ for every $i_1 <
i_2 < \partial$ and $M_2 \models a^n_{i,0} < a^n_{i,1} < \ldots <
a^n_{i,n(*)-1}$ for $i < \partial$.  By $\boxtimes_1$ we know that $\langle
\sigma_*(a^n_{i,0},\dotsc,a^n_{i,n(*)-1}):i < \partial\rangle$ is an
indiscernible set in $M_1$ hence an indiscernible set over
acl$_{M_2}(\emptyset)$.  By $\boxtimes_2$, its type over
acl$_{M_2}(\emptyset)$ does not depend on $n$.  As $|\bold I_n| >
\partial > (2^{|T|})^{\bold L[T,Y,\bold J]}$, we easily get, see
\cite[Ch.V,2.5,pg.244]{Sh:c} that this indiscernible set is not orthogonal
to the indiscernible set $\{a_i:i < \mu\}$.  Also easily letting $f:Z
\rightarrow \partial$ be one to one order preserving,  the type which
$\langle c_i:i \in Z\rangle$ realizes over acl$_{M_2}(\emptyset)$ in
$M_3$ is the same as the type of $\langle
\sigma_*(a^n_{f(i),0},\ldots):i \in Z\rangle$ 
realized in $M_2$ over acl$_{M_2}(\emptyset)$ for $n \ge
2n(*)+1$, as for formulas with $\le m$ variable we consider $n > m$.
As $\bold J$ was chosen to be indiscernible not orthogonal to
tp$(a^n_{0,0},\text{acl}_{M_2}(\emptyset),M_1)$, i.e., to the indiscernible
set $P^{M_4}$, we get a contradiction.  So there is no $\bold J$ as
in $(*)_2$.  As $\bold I$ is minimal, it follows that in $\bold V$,
if for $\ell=1,2$ the set $\bold J_\ell \subseteq {}^{n(\ell)}(M_4)$
is indiscernible of cardinality $\ge \theta^+$ \ub{then} $\bold J_1,
\bold J_2$ are not orthogonal to $\bold I$ hence $\bold J_1,\bold
J_2$ not orthogonal (e.g. works in $\bold L[T,Y,\bold J_1,\bold J_2]$).
\nl  
But this says that $M_4$ is a model as required in the conclusion of
\scite{uni.1}.  \hfill$\square_{\scite{uni.1}}$
\enddemo
\bigskip

\remark{\stag{uni.1A} Remark}  1) By $\bold
F^f_{\aleph_0}$-constructions (see \cite[IV]{Sh:c}) we can get
models with peculiar properties. 
\nl
2) On absoluteness see \scite{ss.1}.  
\nl
3) In fact by \cite[\S1]{Sh:300f}, we can assume that $\langle
\sigma(b_0,\dotsc,b_{m-1}\rangle:m<n,b_0 <^{M_2} \ldots
<^{M_2} b_{m-1}$ are from $\bold I_n\rangle$ where $\sigma = 
\sigma(x_0,\dotsc,x_{m-1})$ and $\sigma$ is a 
$\tau(M_2)$-term, is (fully)
indiscernible in the model $M_2 \restriction \tau_T$, i.e., in $M_1$,
see definition there.
But the argument above is simpler.
\endremark
\bigskip

\demo{\stag{uni.2} Conclusion}  If $T$ is stable and categorical in
$\lambda > |T|$ \ub{then}
(in $\bold L[T,Y]$ where $Y \subseteq$ Ord):
\mr
\item "{$\boxtimes$}"  $(a) \quad T$ is uni-dimensional 
\sn
\item "{${{}}$}"  $(b) \quad T$ is superstable
\sn
\item "{${{}}$}"  $(c) \quad T$ has no two cardinal models
\sn
\item "{${{}}$}"  $(d) \quad D(T)$ has cardinality $\le |T|$; moreover
$D(T) \in \bold L[T]$ and 
\nl

\hskip25pt $\bold L[T] \models |D(T)| = \lambda$.
\endroster
\enddemo
\bigskip

\demo{Proof}  Assume clause (a) fails and we shall 
produce two models of cardinality
(and universe) $\lambda$.  The first $N_1$ is from \scite{uni.1}.  The
second is a model $N_2$ such that there are indiscernible $\bold
I,\bold J \subseteq N_1$ (or ${}^{\omega >}(N_1)$) of 
cardinality $\lambda$ which are
orthogonal; this contradicts the categoricity hence clause (a).  

The superstability, i.e.,  clause (b) follows from clause (a) by
Hrushovski \cite{Hr89d}.

Clause (c), no two cardinal models follows from clause (a) 
by \cite[V,\S6]{Sh:c}.

Now $|D(T)| \le |T|$ (clause (d)) is trivial as otherwise we have two
models $M_1,M_2$ of $T$ of cardinality $\lambda$ such that some $p
\in D(T)$ is realized in one but not the other (i.e., first choose
$M_1 \in \bold L[T]$ realizing 
$\le |T|$ types.  Clearly $\{p \in D(T):p$ is realized in
$M\}$ is a well ordered set so by the assumption we can choose $p \in D(T)$
not realized in $M_1$ and lastly choose $M_2$ realizing $p$).
\hfill$\square_{\scite{uni.2}}$
\enddemo
\bigskip

\proclaim{\stag{uni.3} Claim}  If clauses (b),(c),(d) of 
\scite{uni.2} hold and $T$ is categorical in $\lambda > |T|$ \ub{then}:
\mr
\item "{$(e)$}"  any model $M$ of $T$ of cardinality $\mu$, for any
$\mu > |T|$ is $\aleph_0$-saturated.
\endroster 
\endproclaim
\bigskip

\demo{Proof}  Assume clause $(e)$ fails as exemplified by $M$ and we
shall get contradiction to clause (c) of \scite{uni.2}, so \wilog \,
the universe of $M$ is $\mu$.

For any $Y \subseteq \text{ Ord}$ working in 
$\bold L[T,Y,M]$ we can find $\bar a \subseteq M$ and formula
$\varphi(x,\bar y) \in \Bbb L_{\tau(T)}$ such that $\varphi(x,\bar a)$
is a weakly minimal formula in $M$, existence as in \cite{Sh:31}.  
Let $N_0 \prec M$ be of cardinality $|T|$ such that
$\bar a \subseteq N_0$ and $N_0 \in \bold L[T,M]$.
\enddemo
\bn
\ub{Case 1}:  $\{p \in {\bold S}(N_0):p$ is realized in $M\}$ has power
$> |T|$ in $\bold V$.  

But as $|M|$ is an ordinal this set is well
ordered so the proof of \scite{1.2} applies contradicting categoricity
in $\lambda$ and we get more than needed.
\bn
\ub{Case 2}:  Not Case 1 but there is a finite $A \subseteq M$ such that
$\bar a \subseteq A$ and $p \in {\bold S}(A),\varphi(x,\bar a) \in p$
and the type $p$ is omitted by $M$.

As in \cite{Sh:31} (using ``not Case 1"
here instead ``$T$ stable in $|T|$" there) we can find (in $\bold
L[T,M]$) a model $M'$ such that $M \prec M',M'$ omits the type $p$ and
$\|M'\| \ge \lambda$, so by DLST (= the downward Lowenheim-Skolem-Tarksi)
some $N_1 \prec M'$ has cardinality $\lambda$ and is not
$\aleph_0$-saturated.  Hence for some complete type $p(\bar x,\bar y)
\in D(T)^{\bold L[T,M]}$, for some $\bar b \in {}^{\ell g(\bar
y)}(N_1)$, the model $N_1$ omits the type $p(\bar x,\bar b)$ which is
a type, i.e. finitely satisfiable in $N_1$.  

By clause (d) of $\boxtimes$ of \scite{uni.2} we have 
$|D(T)| \le |T|$ in $\bold L[T,Y]$ and $D(T)$ is included in 
$\bold L[T,Y]$.  So in
$\bold L[T,Y]$, for every finite $A \subseteq N \models T$, 
${\bold S}(A,N)^{\bold V}$ is the same as ${\bold S}(A,N)$ computed in $\bold
L[T,Y]$ and is there of cardinality $\le |D(T)|$ hence absolute.  So 
in $\bold L[T,Y]$ we can find a model $N_2$ of $T$
of cardinality $\lambda$ which is $\aleph_0$-saturated.  
\nl
[Alternatively to this, we can choose a model $N_2$ of
cardinality $\lambda$ such that: if $\bar b' \in {}^{\ell g(\bar y)}
N_2$ realizes tp$(\bar b',\emptyset,M')$ then for some $\bar a \in
{}^{\ell g(\bar x)} N_2$ the sequence $\bar a' \char 94 \bar b'$
realizes $p(\bar x,\bar y)$.]
\nl
By the previous paragraphs this is a contradiction to categoricity.
\bn
\ub{Case 3}:  Neither Case 1 nor Case 2.
\bn
\ub{Subcase A}:  $T$ countable.

Let $N_1$ be such that
\mr
\item "{$\circledast$}"  $(a) \quad N_1 \prec M$ is countable
\sn
\item "{${{}}$}"   $(b) \quad \bar a \subseteq N_1$
\sn
\item "{${{}}$}"  $(c) \quad$ if $\bar a \subseteq A \subseteq N_1,A$ 
finite and $p$ is a non-algebraic type satisfying
\nl

\hskip25pt  $\varphi(x,\bar a) \in p \in {\bold S}(A,M)$ 
then $p$ is realized in $N_1$
\ermn
(possible as by clause (d) of $\boxtimes$ of \scite{uni.2} the set
$D(T)$ is countable and ``neither case 1 nor case 2").

Let $N_2$ be a countable saturated model of $T$ such that $N_1 \prec
N_2$.  We can build an
elementary embedding $f$ (still working in $\bold L[T,M]$) from $N_1$
into $N_2$ such that $f(\varphi(N_1,\bar a)) = (\varphi(N_2,\bar a))$.  This
contradicts clause (c) of $\boxtimes$ of \scite{uni.2}.
\bn
The last subcase is not needed for this section's main theorem
\scite{st.0}, (but is needed for \scite{uni.2}).
\bn
\ub{Subcase B}:  $T$ uncountable.  

So possibly increasing $Y \subseteq \text{\rm Ord}$, in 
$\bold L[T,Y]$ we have two models $M_1,M_2$ of $T,M_1$
is $\aleph_0$-saturated, $M_2$ is not but $\varphi(x,\bar a),M_2$
fails cases 1 and 2; we work in $\bold L[T,Y]$.  
Let $\ell g(\bar a)=n$ and $T^+ \in \bold L[T,Y]$ be the first
order theory in the vocabulary $\tau^+ = \tau_T \cup \{c_\ell:\ell < n\}
\cup \{P\}$ where $c_\ell$ an individual constant, $P$ a unary predicate such
that $M^+ = (M,c^{M^+}_0,\dotsc,c^{M^+}_{n-1},P^{M^+})$ is a model of
$T^+$ iff $M = M^+ \restriction \tau$ is a model of
$T,\varphi(M,c^{M^+}_0,\dotsc,c^{M^+}_{n-1})$ is infinite and
$\subseteq P^{M^+}$ and $M \restriction P^{M^+} \prec M$.  As $T$ is
uni-dimensional (more specifically clause (c) of $\boxtimes$ of
\scite{uni.2} + \cite[V,\S6]{Sh:c}) $T^+$ is inconsistent, hence for
some finite $\tau' \subseteq \tau,T^+ \cap \Bbb L(\tau' \cup
\{c_0,\dotsc,c_{n-1},P\})$ is inconsistent.  Now choose $\bar a_1
\in {}^n(M_1)$ realizing tp$(\bar a,\emptyset,M_2)$, let $\bar a_2 =
\bar a$ and let $\chi$ be large enough, ${\frak B} \prec ({\Cal
H}(\chi)^{\bold L[T,Y]},\in)$ be countable such that $\{M_1,M_2,\tau',\bar a_1,\bar a_2\}
\in {\frak B}$; recall that we are working in $\bold L[T,Y]$.  
Now replacing $M_1,M_2$ by $(M_1 \restriction \tau')
\cap {\frak B},(M_2 \restriction \tau') \cap {\frak B}$ we get a
contradiction as in Subcase A.  \hfill$\square_{\scite{uni.3}}$
\bigskip

\demo{Proof of Theorem \scite{st.0}} By \scite{0.C}(2) and \scite{0.C.1}
\wilog \, $T$ is complete,
$T \subseteq {\Cal H}(\aleph_0)$. Trivially $(B) \Rightarrow (A)$, by
\scite{1.4} we have $(C) \Rightarrow (B)$ and by absoluteness $(C)
\Leftrightarrow (D) \Leftrightarrow (E)$, so it suffices to prove 
$(C)$ assuming $(A)$.  By \scite{st.6} the theory $T$ is stable hence the
assumption of \scite{uni.2} holds hence its conclusion, i.e.,
$\boxtimes$ of \scite{uni.2} holds whenever $Y \subseteq \text{ Ord}$,
in particular $D(T) \in \bold L[T]$.  So by \scite{uni.3} we can
conclude: every model of $T$ of cardinality $\lambda > \aleph_0$ is
$\aleph_0$-saturated (in $\bold V$ or, equivalently, 
in $\bold L[T,M]$ when $M$ has
universe $\lambda$).    If $T$ is $\aleph_0$-stable use \scite{1.3}.
So we can assume $T$ is not $\aleph_0$-stable but is superstable (recall
clause (b) of $\boxtimes$ of \scite{uni.2}) hence $T$
is not categorical in $\aleph_0$ (even has $\ge \aleph_0$
non-isomorphic models, by a theorem of Lachlan, see, e.g.,
\cite{Sh:c}), in any $\bold L[T,Y]$.  So by Kueker 
conjecture (proved by Buechler \cite{Be84} for $T$
superstable and by Hrushovski \cite{Hr89} for 
stable $T$), we get contradiction.  \hfill$\square_{\scite{st.0}}$ 
\enddemo
\bigskip

\remark{\stag{uni.4.1} Remark}  See more in \cite{Sh:F701} about $T$
which is categorical in the
cardinal $\lambda > |T|,T$ not categorical in some $\mu > |T|$.
\endremark
\newpage

\head {\S3 A dichotomy for $\dot I(\aleph_\alpha,T)$: bounded or 
$\ge |\alpha|$}  \endhead  \resetall 
 \spuriousreset
\bigskip

Our aim is to understand the lower part of the family of functions
$\dot I(\lambda,T),T$ countable: either $(\forall \alpha)
\dot I(\aleph_\alpha,T) \ge |\alpha|$ or $\dot I(\aleph_\alpha,T)$ is 
constant and not too large (for $\alpha$ not too small), see
\scite{ss.10}.  For completeness we give a full proof of \scite{ss.1.1}. 

We need here absoluteness between
models of the form $\bold L[Y]$ and this may fail for ``$\kappa(T) >
\kappa$", ``$T$ stable uni-dimensional".  But usually more is true.
\demo{\stag{ss.1} Observation}  1)  ``$T$ is first order", 
$``\tau_T \subseteq \bold L_\omega"$, $``\tau_T \subseteq \bold L"$,
``$T$ is complete" are $Z^-$-absolute.
\nl
2) For $T$ (not necessarily $\in \bold L$ but $\tau_T$ well orderable,
our standard assumption) which is complete:
\mr
\item "{$(a)$}"   ``$T$ is stable" is $Z^-$-absolute
\sn
\item "{$(b)$}"   ``$T$ is superstable" is $(Z^- +$ DC)-absolute (and
downward $Z^-$-absolute; $Z^-$-absolute if $T \subseteq \bold L$)
\sn
\item "{$(c)$}"  ``$T$ totally transcendental" is
($Z^- +$ DC)-absolute (and downward
$Z^-$-absolute; $Z^-$-absolute if $T \subseteq \bold L$);
``$T$ is $\aleph_0$-stable, $\tau(T) \subseteq \bold L$" is $Z^-$-absolute
\sn
\item "{$(d)$}"  the appropriate ranks are $(Z^- +$ DC)-absolute
($Z^-$-absolute if $T \subseteq \bold L$) as the rank 
of $\{\varphi(x,\bar a)\}$ in $M$
depend just on $T,\varphi(x,\bar y)$ and tp$(\bar a,\emptyset,M)$
\sn
\item "{$(e)$}"   ``$M \text{ a model of } T \text{ and }
\bold I,\bold J \subseteq M$ (or ${}^{\omega >} M$) are
 infinite indiscernible sets, $\bold I,\bold J$ are orthogonal and
 where $T$ is stable", is  $Z^-$-absolute
\sn
\item "{$(f)$}"   ``$T$ is stabe not uni-dimensional" is upward $Z^-$-absolute
\sn
\item "{$(g)$}"  for countable $T$, ``$T$ is stable not
uni-dimensional" is $Z^-$-absolute when $T \subseteq \bold L$, 
\sn
\item "{$(h)$}"  `$T$ is countable stable with the OTOP 
(omitting type order property, see \scite{6.22} below)" 
is ($Z^- +$ DC)-absolute
\sn
\item "{$(i)$}"  ``$M$ is primary over $A,M$ a model 
of the (complete) stable theory $T$" is upward $Z^-$-absolute.
\endroster
\enddemo
\bigskip

\demo{Proof}  E.g. 
\nl
2) \ub{Clause (b)}:

This just asks if the tree ${\Cal T}$ has an $\omega$-branch where the
$n$-th level of ${\Cal T}$ is the set of sequence $\langle
\varphi_\ell(x,\bar y_\ell):\ell < n \rangle$ such that for every
$m,\{\varphi(x_\eta,\bar y_\nu)^{\text{if}(\nu \triangleleft
\eta)}:\ell < n,\nu \in {}^\ell m,\eta \in {}^n m\}$ is consistent
with $T$.
\mn
\ub{Clause (e)}: Recall that this is equivalent to
\mr
 \item "{$(*)_1$}"  Av$(\bold I,\bold I \cup \bold J)$ and Av$(\bold
J,\bold I \cup \bold J)$ are weakly orthogonal types
\nl
which is equivalent to
\sn
\item "{$(*)_2$}"  for every $\varphi = \varphi(\bar x,\bar y,\bar z) \in \Bbb
 L(\tau_T)$ and $\bar b \in {}^{\ell g(\bar z)}(\bold I \cup \bold J)$ for some
 $\psi_\ell(x,\bar z_\ell) \in \Bbb L(\tau_T),\bar c_\ell \in {}^{\ell
 g(\bar z_\ell)} (\bold I \cup \bold J)$ such that 
$\psi_\ell(\bar x,\bar c_\ell)$ is satisfied by infinitely many 
$\bar a \in \bold I$ if $\ell=1,\bar a \in \bold J$ if $\ell=2$ and 
truth value $\bold t$ we have $M \models (\forall \bar x,\bar y)
[\psi_1(\bar x,\bar c_1) \wedge \psi_2(\bar y,\bar c_2)
 \rightarrow \varphi(\bar x,\bar y,\bar c)^{\bold t}]$.
\ermn
\ub{Clause (h)}: 
  
We just ask for the existence of the $\Phi \in \Upsilon^{\text{or}}_T$
 so with $\tau_\Phi$ countable $\supseteq \tau_T$  and
type $p(\bar x,\bar y,\bar z)$ from $D(T)$ such that $(\ell g(\bar y)
 = \ell g(\bar z)$ and) for any linear
 order $I$, which is well orderable EM$_\tau(I,\Phi)$ is a
model of $T$ of cardinality $|T| + |I|$ and 
$p(\bar x,\bar a_s,\bar a_t)$ is realized in it iff $s <_I
t$ (so O.K. for stable). \hfill$\square_{\scite{ss.1}}$
\enddemo
\bigskip

\proclaim{\stag{ss.1.1} Claim}  If $T$ is unstable and $|T| = \aleph_{\beta_*}
< \aleph_\alpha = \lambda$ \ub{then} $\dot I(\lambda,T) \ge |\alpha-\beta_*|$.
\endproclaim
\bigskip

\demo{Proof}  In $\bold L[T]$, let $\Phi$ be as in \scite{st.1} such that
for every a linear order $I$ we have $s,t \in I \Rightarrow 
\text{ EM}(I,\Phi) \models \varphi[\bar a_s,\bar
a_t]^{\text{if}(s<t)}$, where, of course, $\varphi(\bar x,\bar z)
\in \Bbb L(\tau(T))$.  

First, we define for $\gamma \le \aleph_\alpha$

$$
J_\gamma =: \gamma + (\gamma)^*.
$$
\mn
We can specify: the set of members of $J_\gamma$ is
$\{(\gamma,\ell,\zeta):\ell \in \{0,1\},\zeta < \gamma\}$ and 
$(\gamma,\ell_1,\zeta_1) < (\gamma,\ell_2,\zeta_2)$ iff $\ell_1 = 0
\wedge \ell_2 = 1$ or $\ell_1 = \ell_2 = 0 \wedge \zeta_1 < \zeta_2$
or $\ell_1 = \ell_2 = 1 \wedge \zeta_1 > \zeta_2$.

Second, for $\beta \in [\beta_*,\alpha]$ 
let $J^\beta = \dsize \sum_{\gamma \le \aleph_\beta}
J_\gamma + J_\infty$ where $J_\infty = (\aleph_\alpha +1) \times \Bbb Q$
ordered lexicographically.

Third, let $M^\beta_1 = \text{ EM}(J^\beta,\Phi)$.

Lastly, $M^\beta := M^\beta_1 \restriction \tau_T$ - clearly a 
model of $T$ of cardinality $\aleph_\alpha$.  We like to ``recover", ``define"
$\aleph_\beta$ from $M^\beta/\cong$ at least when $\beta \ge \beta_*$.  
This is sufficient as the sequence
$\langle M_\beta:\beta \in [\beta_*,\alpha]\rangle$ exists (in fact in $\bold
L[T]$).  We shall continue after stating \scite{ss.1A}.
\enddemo
\bn
\ub{Discussion}:  1) In ZFC we could recover from the isomorphism types,
stationary subsets modulo the club filter so 
as we get $2^{\aleph_\alpha}$, if, e.g.,
$\aleph_\alpha$ is regular and there are $2^{\aleph_\alpha}$ 
subsets of $\aleph_\alpha$ any two with a
stationary difference so we get $\dot I(\aleph_\alpha,T) = 2^{\aleph_\alpha}$.
But here (ZF) the stationary subsets of a regular uncountable
$\lambda$ may form an ultrafilter or all uncountable cardinals are singulars.
\nl
2) More than \scite{ss.1A} is true in $\bold L[T,Y]$; EM$(J,\Phi)$
satisfies $\otimes_\theta$ iff $J$ has a $(\theta,\theta)$-cut
(provided $J$ has no
$(1,\theta),(\theta,1),(0,\theta),(\theta,0)$ cuts), see below.
\nl
3) See more in \scite{ss.2} on OTOP.
\nl
4) Of course, we can prove theorems saying e.g.: if $\aleph_\alpha >
|T|$ is regular, $T$ unstable then $\dot I(\aleph_\alpha,T) \ge
|{\Cal P}(\aleph_\alpha)$/(the club filter on $\aleph_\alpha$)$|$.
\bigskip

\proclaim{\stag{ss.1A} Subclaim}  If $J = J^\beta,M = M^\beta$ are as
above and $Y \subseteq { \text{\rm Ord\/}}$ satisfies $M \in \bold
L[T,Y]$ \ub{then} in $\bold L[T,Y]$ for any 
regular cardinal $\theta$ (of $\bold L[T,Y]$)
\mr
\item "{$(*)$}"  $\otimes_\theta \Leftrightarrow \theta >
\aleph_\beta$ where
{\roster
\itemitem{ $\otimes_\theta$ }  if $p$ is a set of 
$\Delta$-formulas with parameters from $M$ of cardinality $\theta$ where
$$
\Delta =: \{\varphi(\bar x,\bar z_1) \wedge \neg \varphi(\bar x,\bar z_1)\}
$$
and any subset $q$ of $p$ of cardinality $< \theta$ is realized in 
$M$ \ub{then} some $q \subseteq p,|q| = \theta$ is realized in $M$.
\endroster}
\endroster
\endproclaim
\bigskip

\demo{Proof of Claim \scite{ss.1.1} from the subclaim \scite{ss.1A}}
Why does this subclaim help us to prove the Theorem?  Assume $\beta_* \le
\beta_1 < \beta_2 \le \alpha$ and we consider
$M^{\beta_1},M^{\beta_2}$ as above and toward a contradiction we
assume that there is an isomorphism $f$ from $M^{\beta_1}$ onto
$M^{\beta_2}$.

Let $Y \subseteq \text{ Ord}$ code $T,M^{\beta_1},M^{\beta_2}$ and $f$.
So $\bold L[T,M^{\beta_1},Y] = \bold L[Y] = \bold L[T,M^{\beta_2},Y]$.  In
this universe let $\theta$ the first cardinal greater than the ordinal 
$> \aleph^{\bold V}_{\beta_1}$ so $\aleph^{\bold V}_{\beta_1} < 
\theta \le \aleph^{\bold V}_{\beta_1 +1} < \aleph_{\beta_2}$.
\bn
\ub{Question}:  Why we cannot prove that 
$\theta = \aleph^{\bold V}_{\beta_1 +1}$?  As possibly
$\bold L[Y] \models \aleph^{\bold V}_{\beta_1 +1}$ is singular or just a
limit cardinal.
\bn
\ub{Note}:  Maybe every $\bold L[Y]$-cardinal from 
$(\aleph^{\bold V}_{\beta_1},\aleph^{\bold V}_{\beta_1 +1})$
have cofinality $\aleph_0$ in $\bold V$!

But in $\bold L[Y],\aleph^{\bold V}_\beta,\aleph^{\bold V}_{\beta +1}$
are still cardinals so the successor of $\aleph^{\bold V}_\beta$ 
in $\bold L[Y]$ is $\le \aleph^{\bold V}_{\beta +1}$ but in $\bold L[Y]$ this
successor, $\theta$ is regular.
(In $\bold V,\theta$ may not be a cardinal at all).
In $\bold L[T]$ there are many possibilities for $\theta$ (it was
defined from $Y$!) and we have built $M^\beta_1$
\ub{before} knowing who they will be in $\bold L[Y]$ so

$$
\theta > \aleph_{\beta_1} \Leftrightarrow M^{\beta_1} \models
\otimes_\theta \Leftrightarrow M^{\beta_2} \models \otimes_\theta
\Leftrightarrow \theta > \aleph_{\beta_2}
$$

$$
\text{(the first } \Leftrightarrow \text{ by } (*) \text{ of the
subclaim and the second } 
\Leftrightarrow \text{ as } f \text{ is an isomorphism)}
$$
\mn
but $\aleph_{\beta_1} < \theta \le \aleph_{\beta_2}$; \ub{contradiction}.
\enddemo
\bigskip

\demo{Proof of the subclaim \scite{ss.1A}}  I.e., in $\bold L[T,Y]$ we
have to prove:
\mr
\item "{$(*)$}"  $[\otimes_\theta \Leftrightarrow \theta > \aleph_\beta]$ 
\ermn
\ub{First} we will prove:
\mr
\item "{$(*)_1$}"   $\theta \le \aleph_\beta \Rightarrow \neg
\otimes_\theta$.  
\ermn
By the choice of $J = J^\beta$ clearly $J_\theta$ is an interval of $J$ so let

$$
p =: \{\varphi(\bar x,\bar a_{(\theta,1,i)} \wedge \neg \varphi(\bar x,\bar
a_{(\theta,0,i)}):i<\theta\}.
$$
\mn
Let $q \subseteq p,|q|< \theta$ now as $\theta$ is regular (in $\bold
L[T,Y]$) for some $j<\theta$ we have

$$
q \subseteq p_j = \{\varphi(\bar x,\bar a_{(\theta,1,i)}) \wedge 
\neg \varphi(\bar x,\bar a_{(\theta,0,i)}):i < j\}.
$$
\mn
We have a natural candidate for a sequence realizing $q$: the sequence 
$\bar a_{(\theta,1,j)}$.  Now

$$
i<j \Rightarrow (\theta,1,j) <_{J_\theta} (\theta,1,i) \Rightarrow M \models
\varphi[\bar a_{(\theta,1,j)},\bar a_{(\theta,1,i)}]
$$

$$
i<j \Rightarrow (\theta,0,i) <_{J_\theta} (\theta,1,j) 
\Rightarrow M \models \neg \varphi[\bar a_{(\theta,1,j)},
\bar a_{(\theta,0,i)}].
$$
\mn
So we have proved that every $q \subseteq p,|q| < \theta$ is realized
in the model.  Secondly, we need to show:
\mr
\item "{$\otimes$}"  no $\bar a \in M$ satisfies $\theta$ of formulas
from $p$.
\ermn
Assume toward contradiction that $\bar a$ is a counterexample.
\nl
So we can find $n<\omega$, a finite sequence of terms 
$\bar \sigma(\bar x_0,\dotsc,\bar x_{n-1})$ 
from $\tau(\Phi)$ and $t_0 <_J t_1 <_J \ldots <_J t_{n-1}$ such that 
$\bar a = \bar \sigma(\bar a_{t_0},\dotsc,\bar a_{t_{n-1}})$.  
Now for each $\ell$ for some $i_\ell <
\theta,t_\ell$ is not in the interval $((\theta,0,i_\ell),
(\theta,1,i_\ell))_J$.

Let:

$$
j^* = \max[\{i_\ell+1:\ell < n\} \cup \{1\}].
$$
\mn
Now consider $\varphi(\bar x,\bar a_{(\theta,1,j)} \wedge \neg \varphi(\bar
x,\bar a_{(\theta,0,j)}))$ for $j \in [j^*,\theta)$.
\nl
So $t_\ell <_J (\theta,1,j)
\equiv t_\ell <_J (\theta,0,j)$ for $\ell=0,\dotsc,n-1$ hence $M \models
\varphi[\bar \sigma(\bar a_{t_0},\dotsc,\bar a_{t_{n-1}}),
\bar a_{(\theta,1,j)}] \Leftrightarrow M \models \varphi
[\bar \sigma(\bar a_{t_0},\dotsc,\bar a_{t_{n-1}}),\bar a_{(\theta,0,j)}]$.  So
$\bar a = \sigma(\bar a_{t_0},\dotsc,\bar a_{t_{n-1}})$ 
fail the $j$-th formula from
$p$ for $j \in [j^*,\theta)$.  So $p$ really exemplifies the $\neg
\otimes_\theta$.  So we have proved $(*)_1$ which is one 
implication of the Subclaim.
\enddemo
\bn
Now we will prove:
\mr
\item "{$(*)_2$}"  if $\bold L[T,Y] \models ``\theta$ is regular 
$> \aleph_\beta"$ then $\otimes_\theta$.
\ermn
So let $p = \{\varphi(\bar x,\bar a_i) \wedge \neg \varphi
(\bar x,\bar b_i):i < \theta\} \in \bold L[T,Y]$ be given.  
For $j <\theta$ let

$$
p_j =: \{\varphi(\bar x,\bar a_i) \wedge \neg \varphi(\bar x,\bar
b_i):i<j\}.
$$
\mn
So some $c_j \in M$ realizes it and let $(\bar a_i,\bar b_i,\bar c_i) = \langle
\bar \sigma^k_i(a_{t^i_0},\dotsc,a_{t^i_{n_i-1}}):k = 0,1,2 \rangle$
where $\bar \sigma^k_i$ is a finite sequence of 
terms from $\tau(\Phi)$ and
$J \models t^i_0 < t^i_1 < \ldots < t^i_{n_i-1}$; note that we can
make $\langle t^i_\ell:\ell < n_i\rangle$ not to depend on $k$ because
we can add dummy variables.

As $\tau(\Phi)$ is of cardinality $< \theta = \text{\rm cf}(\theta)$ (in
$\bold L[T,Y]$), for some $\sigma^k_*,n_*$ the set $S =
\{i:\sigma^k_i = \sigma^k_*$ for $k = 0,1,2$ and $n_i=n_*\}$
is unbounded in $\theta$.

Recall

$$
J^\beta = \dsize \sum_{\gamma \le \aleph_\beta} J_\gamma +
(\aleph_\alpha +1) \times \Bbb Q.
$$
\mn
So for some $m_i \le n_*$

$$
t^i_\ell \in \dsize \sum_{\gamma \le \aleph_\beta} J_\gamma
\Leftrightarrow \ell < m_i
$$
\mn
shrinking $S$ \wilog \, $i \in S \Rightarrow m_i = m_*$.
\sn
\ub{Now} $\bold L[T,Y] \models ``|\dsize \sum_{\gamma \le \aleph_\beta}
J_\gamma| \le \dsize \sum_{\gamma \le \aleph_\beta} |J_\gamma| =
\dsize \sum_{\gamma \le \aleph_\beta} (|\gamma| + \aleph_0) \le
\aleph_\beta < \theta= \text{ \rm cf}(\theta)$".

So \wilog
\mr
\item "{$\circledast_1$}"  $\ell < m_* \Rightarrow t^i_\ell =
t^*_\ell$ for $i \in S$ and for $\ell \in [m_*,n_*)$ let $t^i_\ell =
(\varepsilon^i_\ell,q^i_\ell)$ where $q^i_\ell \in \Bbb Q$.
\ermn
Clearly for $q^i_\ell$ there are $\aleph_0$ possibilities so \wilog,
\, for each $\ell \in [m_*,n_*)$
\mr
\item "{$\circledast^\ell_2$}"  $q^i_\ell = q^*_\ell$ for $i \in S$,
\sn
\item "{$\circledast^\ell_3$}"  $\langle \varepsilon^i_\ell:i \in S
\rangle$ is \ub{constant} say $\varepsilon^*_\ell$ \ub{or is strictly
increasing} with limit $\varepsilon^*_\ell$ and is strictly increasing
iff $\ell \in u$
\ermn
so \wilog 
\mr
\item "{$\circledast_4$}"  $(i) \quad$ if $\ell_1 \ne \ell_2$ are in
the interval $[m_*,n_*)$ and $\varepsilon^*_{\ell_1} < 
\varepsilon^*_{\ell_2}$ then
\nl

\hskip25pt $i,j \in S \Rightarrow \varepsilon^i_{\ell_1} \le
\varepsilon^*_{\ell_1} < \varepsilon^j_{\ell_2}$
\sn
\item "{${{}}$}"  $(ii) \quad$ if $\ell_1 \ne \ell_2 \in [m_*,n_*)$
and $\varepsilon^*_{\ell_1} = \varepsilon^*_{\ell_2} \wedge \ell_1 \in
u \wedge \ell_2 \notin u$ and $i<j$ are in 
\nl

\hskip25pt $S$ then $\varepsilon^i_{\ell_1} < \varepsilon^j_{\ell_2}$.
\ermn
We choose $t_0 <_J t_1 <_J \ldots <_J t_{n-1}$ which satisfies
\mr
\item "{$\circledast_5$}"   $(a) \quad$ if 
$\ell < m_*$ then $t_\ell = t^*_\ell$
\sn
\item "{${{}}$}"   $(b) \quad$ if 
$\ell \in [m_*,n_*)$ and $\langle t^i_\ell:i \in S
\rangle$ is constant then $t_\ell = t^*_\ell$
\sn
\item "{${{}}$}"  $(c) \quad$ if $\ell \in [m_*,n_*),\langle t^i_\ell:i \in S
\rangle$ is not constant (i.e. $\ell \in u$) then:
\nl

\hskip25pt (recall that $\langle q^i_\ell:i \in S \rangle$ is 
constantly $q^*_\ell,\langle \varepsilon^i_\ell:i \in S \rangle$
\nl

\hskip25pt  is strictly increasing
with limit $\varepsilon^*_\ell$) we choose
$t_\ell = (\varepsilon_\ell,q_\ell)$ \nl

\hskip25pt such that $\varepsilon_\ell =
\varepsilon^*_\ell,q_\ell = \text{ min}(\{0\} \cup 
\{q^*_k:k \in [m^*,n^*)\}) - n^*+\ell$
\nl

\hskip25pt (the computation is in $\Bbb Q!$)
\ermn
Hence
\mr
\item "{$\circledast_6$}"  $(\alpha) \quad q_\ell$ is $< q^*_k$ for
every $k \in [m^*,n_*)$ when $\ell \in u$
\sn
\item "{${{}}$}"  $(\beta) \quad$ if $\varepsilon^*_\ell = 
\varepsilon^*_k$ and $\ell,k \in u$ then $q_\ell < q_k \equiv \ell < k$.
\ermn
\ub{Now} note that:
\mr
\item "{$\circledast_7$}"  for $\varepsilon < \zeta < \theta$ from $S$, in
$J$ the quantifier free types of $\langle
t^\varepsilon_\ell:\ell < n_* \rangle {}^\frown \langle t_\ell:\ell < n_*
\rangle$ and $\langle t^\varepsilon_\ell:\ell < n_* \rangle {}^\frown
\langle t^\zeta_\ell:\ell < n_* \rangle$ are equal [all the shrinking
was done for this].
\ermn
Now for $\varepsilon < \zeta$ from $S$, by 
the original choice above $M^\beta \models \varphi[\bar c_\zeta,
\bar a_\varepsilon] \wedge \neg \varphi[\bar c_\zeta,\bar b_\varepsilon]$
that is: $M^\beta_1 \models \varphi[\bar
\sigma^0_*(a_{t^\zeta_0},\ldots),
\bar \sigma^1_\varepsilon (a_{t^\varepsilon_0},\ldots)] 
\wedge \neg \varphi[\bar \sigma^0_* (a_{t^\zeta_0},\ldots),
\bar \sigma^2_\varepsilon(a_{t^\varepsilon_0},\ldots)]$.

By the last sentence and 
$\circledast_7 +$ indiscernibility of $\langle \bar a_t:t \in
J\rangle$ in $M^\beta_1$ we have 
\nl
$M \models \varphi[\bar \sigma^0_*(a_{t_0},\ldots),\bar \sigma^1_\varepsilon
(a_{t^\varepsilon_0},\ldots)]
\wedge \neg \varphi[\bar \sigma^0_*(a_{t_0},\ldots),
\bar \sigma^2_\varepsilon(a_{t^\varepsilon_0},\ldots)]$.

Let $\bar c = \bar \sigma^0_*(a_{t_0},\ldots)$ in $M^\beta_1$-sense, so
$\varepsilon \in S \Rightarrow M \models \varphi[\bar c,\bar
a_\varepsilon] \wedge \neg \varphi[\bar c,\bar b_\varepsilon]$.  Hence
$\{\varphi(\bar x,\bar a_\varepsilon) \wedge \neg \varphi(\bar x,\bar
b_\varepsilon):\varepsilon \in S\}$ is realized in $M^\beta$ and 
$\bold L[T,Y]
\models ``|S| = \theta"$ as promised.
\nl
${{}}$  \hfill$\square_{\scite{ss.1.1}}$
\bigskip

\proclaim{\stag{ss.1.2} Claim}  If $T$ is stable not uni-dimensional,
$|T|= \aleph_\beta < \aleph_\alpha = \lambda$ \ub{then} $\dot I(\lambda,T) \ge
|\alpha-\beta|$.
\endproclaim
\bigskip

\demo{Proof}  As in \scite{uni.1}; if $\gamma \in [\beta,\alpha]$
\ub{then} there is a model $M$ of $T$ of cardinality $\lambda$ such
that $M$ satisfies $(*)_\gamma$ but not $\beta \le \gamma_1 < \gamma
\Rightarrow \neg(*)_{\gamma_1}$ where
\mr
\item "{$(*)_\gamma$}"  if $\bold I,\bold J \subseteq M$ are infinite
orthogonal indiscernible sets and $|\bold I| = \lambda$ then $|\bold
J| \le \aleph_\gamma$.  \hfill$\square_{\scite{ss.1.2}}$
\endroster
\enddemo
\bigskip

\demo{\stag{ss.1.3} Conclusion}  If $\lambda = \aleph_\alpha >
\aleph_\beta = |T|$ and $\dot I(\lambda,T) < |\alpha-\beta|$ \ub{then}
(in $\bold L[T,Y]$ when $Y \subseteq$ Ord)
\mr
\item "{$\boxtimes_T$}"  $(a) \quad T$ is stable and 
uni-dimensional 
\sn
\item "{${{}}$}"  $(b) \quad T$ is superstable
\sn
\item "{${{}}$}"  $(c) \quad T$ has no two cardinal models
\sn
\item "{${{}}$}"  $(d) \quad D(T)$ has cardinality $\le |T|$ or
cardinality $< |\alpha-\beta|$.
\endroster
\enddemo
\bigskip

\demo{Proof}  $T$ is stable by \scite{ss.1.1} and uni-dimensional by
\scite{ss.1.2} so clause (a) holds.  This implies clause (c), see
\cite[V,\S6]{Sh:c}.  Clause (d) is trivial by now and clause (b)
follows from clause (a) by Hrushovski \cite{Hr89d}.
\hfill$\square_{\scite{ss.1.3}}$ 
\enddemo
\bigskip

\proclaim{\stag{ss.2} Claim}  In \scite{ss.1.3} we can add to
$\boxtimes_T$ also clause (e) and if $T$ is countable also clause (f) where
\mr
\item "{$\boxtimes_T$}"  $(e) \quad T$ fails the {\rm OTOP} (see
\cite[XII,Def.4.1,pg.608]{Sh:c} or \scite{6.22}(1) below) 
\sn
\item "{${{}}$}"  $(f) \quad T$ has the prime existence property \nl

\hskip25pt (see \cite[XII,Def.4.2,pg.608]{Sh:c} or \scite{6.22}(2)
below) hence
\nl

\hskip25pt for ${\frak C}_T$ a model of $T$ with universe $|{\frak C}_T|
\subseteq \bold L$:
{\roster
\itemitem{ ${{}}$ }  for any non-forking tree $\langle N_\eta:\eta
\in {\Cal T} \rangle$ of models 
$N_\eta \prec {\frak C}_T$, there is a prime (even
primary, i.e. $\bold F^t_{\aleph_0}$-primary) model $N \prec {\frak C}_T$ 
over $\cup\{N_\eta:\eta \in {\Cal T}\}$, it is
unique up to isomorphism over $\cup\{N_\eta:\eta \in {\Cal T}\}$.
\endroster} 
\endroster
\endproclaim
\bigskip

\demo{Proof}  Clause (e) holds exactly as for stability, i.e., as in
\scite{ss.1.1} only the formulas $\varphi(\bar x,\bar y)$ are not first
order but an infinite conjunction of such formulas.  
Clause (f) follows by \cite[XII]{Sh:c}, i.e., it holds
in any $\bold L[T,Y]$ which suffices.  \hfill$\square_{\scite{ss.2}}$ 
\enddemo
\bigskip

\definition{\stag{6.22} Definition}  1) $T$ has OTOP if for some type
$p = p(\bar x,\bar y,\bar z)$ in $\Bbb L(\tau_T)$
the theory $T$ has it for $p$, which means that for every
$\lambda$ for some model $M$ of $T$ with well ordered universe and 
$\bar b_\alpha \in {}^{\ell g(\bar y)}M, 
\bar c_\alpha \in {}^{\ell g(\bar z)}M$,
for $\alpha,\beta < \lambda$ we have: for any $\alpha,\beta < \lambda$
the model $M$ realizes
the type $p(\bar x,\bar b_\alpha,\bar c_\beta)$ iff $\alpha < \beta$.
\nl
2) $T$ has the prime existence property \ub{when} for every
 triple $(M_0,M_1,M_2)$ in stable amalgamation in a model 
${\frak C}_T$ of $T$ such that $|{\frak C}_T|$ is well orderable 
(so $M_\ell \prec {\frak C}_T$), the set of isolated types is 
dense in ${\bold S}^m(M_1 \cup M_2)$ for every $m$.  
\enddefinition
\bigskip

\proclaim{\stag{ss.3} Claim}   [$T$ countable]  We can add clause
(g) below to $\boxtimes_T$ from \scite{ss.1.3} + \scite{ss.2}:
\mr
\item "{$\boxtimes_T$}"  $(g) \quad$ if clause (A) \ub{then} for some $M'$ 
clause (B) below holds 
\nl

\hskip25pt (both in $\bold L[T,Y]$) where
{\roster
\itemitem{ $(A)$ }  $(\alpha) \quad M_\emptyset \prec M_{\{i\}} \prec
M^*$ are countable models of $T$ for $i < \omega \times 2$
\sn
\itemitem{ ${{}}$ }  $(\beta) \quad (M_{\{i\}},c)_{c \in M_\emptyset}
\cong (M_{\{0\}},c)_{c \in M_\emptyset}$ for $i < \omega$ that is
$M_{\{i\}}$ is 
\nl

\hskip60pt isomorphic to $M_{\{0\}}$ over $M_\emptyset$ for $i < \omega$
\sn
\itemitem{ ${{}}$ }  $(\gamma) \quad (M_{\{\omega +i\}},c)_{c \in
M_\emptyset} \cong (M_{\{\omega\}},c)_{c \in M_\emptyset}$ that is
$M_{\{\omega +i\}}$ is isomorphic 
\nl

\hskip45pt to $M_{\{\omega\}}$ over $M_\emptyset$ for $i < \omega$
\sn
\itemitem{ ${{}}$ }  $(\delta) \quad \{M_{\{i\}}:i < \omega \times
2\}$ is independent over $M_\emptyset$ inside $M^*$
\sn
\itemitem{ ${{}}$ }  $(\varepsilon) \quad M^*$ is prime over $\cup
\{M_{\{i\}}:i < \omega \times 2\}$
\sn
\itemitem{ $(B)$ }  $(\alpha) \quad M_\emptyset \prec M' \prec M^*$
\sn
\itemitem{ ${{}}$ }  $(\beta) \quad (M',c)_{c \in M_\emptyset} \cong
(M_{\{0\}},c)_{c \in M_\emptyset}$
\sn
\itemitem{ ${{}}$ }  $(\gamma) \quad \langle M_{\{i\}}:i < \omega \rangle
{}^\frown \langle M'\rangle$ is independent over $M_\emptyset$.
\endroster}
\endroster
\endproclaim
\bigskip

\remark{\stag{ss.3.1} Remark}  1) We can formulate (B) closer to
$\circledast_6$ inside the proof of \scite{ss.4}.
\nl
2) We can omit ``$T$ countable" but then have to change $Y$ with the
same proof.
\nl
3) We know more on $T$'s satisfying $\boxtimes_T$ of \scite{ss.1.3}
 by Laskowski \cite{Las88} and Hart-Hrushovski-Laskowski \cite{HHL00}.
\endremark
\bigskip

\demo{Proof}  Note that $|T| = \aleph_0$ and choose the ordinals
 $\beta_*,\alpha_*$ such that
$\beta_* = 0,\lambda = \aleph_{\alpha_*}$; most of the proof we do not
use $\beta_* = 0$ but we use $\boxtimes_T(a)-(f)$.

We do more than is strictly necessary for the proof; we use $\odot$ to
denote definitions, working in $\bold L[T,Y]$ if not said otherwise and
${\frak C}_Y$ is a monster for $T$ in $\bold L[T,Y]$:
\mr
\item "{$\odot_1$}"  $(a) \quad$ for a model $M \prec {\frak C}_Y$ let
$\bold S^{c,\theta}_Y(M) = \{\text{tp}(\bar a,M,N):M \prec N 
\prec {\frak C}$, 
\nl

\hskip25pt $\|N\| \le \theta$ and $\bar a$ enumerates $N\}$, 
omitting $\theta$ means some $\theta$
\sn
\item "{${{}}$}"  $(b) \quad$ in this case we say $N$ realizes $p$
\sn
\item "{${{}}$}"  $(c) \quad$ if $p = \text{ tp}(\bar a,M,N)$ is as
above, \ub{then} we denote 
\nl

\hskip5pt  $|p| = \|N\|$, 
\sn
\item "{$\odot_2$}"  for $\bar\alpha = \langle
\alpha_\varepsilon:\varepsilon < \zeta\rangle$ and $\bar p = \langle
p_\varepsilon:\varepsilon < \zeta\rangle,p_i \in \bold S^c_Y(M)$, we
say $N$ is $(\bar p,\bar\alpha)$-constructed over $M$ when there is
$\bar M$ such that
\sn
\item "{${{}}$}"  $(a) \quad \bar M = \langle M_{\{i\}}:i <
\alpha^\zeta\rangle$, where 
$\alpha^\varepsilon = \dsize \sum_{\xi < \varepsilon}
\alpha_\xi$ for $\varepsilon \le \zeta$
\sn
\item "{${{}}$}"  $(b) \quad M_{\{i\}}$ realizes $p_\varepsilon$ if $i
\in [\alpha^\varepsilon,\alpha^\varepsilon + \alpha_\varepsilon)$
\sn
\item "{${{}}$}"  $(c) \quad \langle M_{\{i\}}:i < \alpha^\zeta\rangle$ is
independent over $M$
\sn
\item "{${{}}$}"  $(d) \quad N$ is primary over $\dbcu_{i < \alpha^\zeta}
M_{\{i\}}$
\sn
\item "{$\odot_3$}"  we say $N$ is $\bar p$-constructed over $M$ if
this holds for some $\bar\alpha$
\sn
\item "{$\odot_4$}"  if $M \prec N \prec {\frak C}_Y,p \in \bold
S^c_Y(M)$ then we say $q$ lifts $p$ or $(p,M)$ to $N$ \ub{when} $q \in \bold
S^c_Y(N)$ and for some $M_1,N_1$ realizing $p,q$ respectively,
tp$(M_1,N)$ does not fork over $M$ and $N_1$ is primary over $N \cup
M_1$
\sn
\item "{$\odot_5$}"  for $M \prec {\frak C}$ and $p_1,p_2 \in \bold
S^c_Y(M)$ we say $p_2$ pushes $p_1$ (in $\bold L[T,Y]$) \ub{when} for some
ordinals $\alpha_1,\alpha_2$ there are $M'_{\{i\}}$ for $i < \alpha_2 +
\alpha_1$ and $\bar M,M^*,M'$ satisfying
\sn
\item "{${{}}$}"  $(a) \quad M^*$ is $(\langle p_1,p_2\rangle,\langle
\alpha_1,\alpha_2\rangle)$-constructed over $M$ as witnessed by
\nl

\hskip25pt $\bar M' = \langle M'_{\{i\}}:i < \alpha_1 + \alpha_2\rangle$
\sn
\item "{${{}}$}"  $(b) \quad M \prec M' \prec M^*$
\sn
\item "{${{}}$}"  $(c) \quad M'$ realizes $p_1$
\sn
\item "{${{}}$}"  $(d) \quad \langle M_{\{i\}}:i < \alpha_1\rangle \char
94 \langle M'\rangle$ is independent over $M$
\sn
\item "{$\odot_6$}"  $(\alpha) \quad$ assume
$p_\varepsilon,q_\varepsilon \in \bold S^c_Y(M)$ for $\varepsilon <
\varepsilon(*)$; we say $(\bar p,\bar \alpha)$ is equivalent to
\nl

\hskip25pt  $(\bar q,\bar \beta)$ when $\bar \alpha = \langle
\alpha_\varepsilon:\varepsilon < \varepsilon(*)\rangle,\bar \beta =
\langle \beta_\varepsilon:\varepsilon < \varepsilon(*)\rangle$ and
\nl

\hskip25pt  there is $M'$ which is both $(\bar p,\bar \alpha)$-constructed over
$M$ and 
\nl

\hskip25pt $(\bar q,\bar \beta)$-constructed over $M$
\sn
\item "{${{}}$}"  $(\beta) \quad$ we may write $p$ instead of $\langle
p \rangle,q$ instead of $\langle q\rangle$, and omitting
\nl

\hskip25pt  $\bar\alpha,\bar \beta$ means ``for some $\bar\alpha,\bar\beta$".
\sn
\item "{$\circledast_1$}"   if $p_1,p_2 \in \bold S^c_Y(M)$ 
and $q$ pushes $p$ \ub{then} in $\odot_5$ 
\wilog \, $\alpha,\beta \le \|M\| + |T| + |p_1| + |p_2|$.
\ermn
[Why?  By the DLST argument.]
\mr
\item "{$\odot_7$}"  $(a) \quad$ let AP$^\theta_Y = \{(M,p_1,q_1)$: in
$\bold L[T,Y],M \prec {\frak C}_Y$ and $p_1,q_1 \in \bold S^c_Y(M)$ have
\nl

\hskip25pt cardinality $\le \theta\}$
\sn
\item "{${{}}$}"  $(b) \quad$ AP$^\theta_Y \models ``(M_1,p_1,q_1) \le
(M_2,p_2,q_2)"$ means that
{\roster
\itemitem{ ${{}}$ }  $(\alpha) \quad$ both triples are from AP$^\theta_Y$
\sn
\itemitem{ ${{}}$ }  $(\beta) \quad M_2$ is $(p_1,q_1)$-constructed over $M_1$
\sn
\itemitem{ ${{}}$ }  $(\gamma) \quad p_2,q_2$ lift $p_1,q_1$ over
$M_2$ respectively
\endroster}
\sn
\item "{$\circledast_2$}"  if AP$^\theta_Y \models ``(M_1,p_1,q_1) \le
(M_2,p_2,q_2)"$ and $q_2$ pushes $p_2$ \ub{then} $q_1$ pushes $p_1$.
\ermn
[Why?  Straight.]
\mr
\item "{$\circledast_3$}"  if $(M,p_1,q_1) \in 
\text{ AP}^{\aleph_{\beta(*)}}_Y$ and $p_1$ does not push $q_1$ \ub{then}
we can find $\mu_0,\mu_1,M_*,p_2,q_2$ and $r$ such that
{\roster
\itemitem{ $(a)$ }  AP$^{\aleph_0}_Y \models ``(M,p_1,q_1) \le
(M_*,p_2,q_2)"$ hence by $\circledast_2$ the type $p_2$ does not push
$q_2$, (this is the only point where we use ``$p_1$ does not push $q_1$")
\sn
\itemitem{ $(b)$ }  $\|M_*\| = \mu_0$
\sn
\itemitem{ $(c)$ }  $r \in \bold S^{c,\mu_0}_Y(M_*)$ and
$\aleph_{\beta(*)} \le \mu_0 < \mu_1 < \lambda$
\sn
\itemitem{ $(d)$ }  $(\langle p_2,q_2\rangle,\langle
\lambda,\mu_1\rangle)$ is equivalent to 
$(\langle r \rangle,\langle \lambda\rangle)$, see $\odot_6$.
\endroster}
\ermn
[Why?  For every $\mu \in [\aleph_{\beta(*)},\lambda)$ let $N^\mu$ 
be $(\langle p_1,q_1\rangle,\langle \lambda,\mu\rangle)$-constructed
over $M$ as witnessed by $\langle N_i:i < \lambda + \mu\rangle$.]

As we are assuming that $\dot I(\lambda,T) < |\alpha_* - \beta_*|$, 
there are $\mu_0,\mu_1$ such that 
$|T| = \aleph_{\beta(*)} \le \mu_0 < \mu_1 < \lambda$ and 
there is an isomorphism $f \in \bold V$ from $N^{\mu_0}$ onto 
$N^{\mu_1}$; of course $f$ is not necessarily from 
$\bold L[T,Y]$.  We now work in $\bold
L[T,Y,f]$ and in the end we use absoluteness (here we use ``$T$ countable").

Now by the DLST argument and properties of $\bold F^t_{\aleph_0}$-primary
 we can find $(u_0,u_1,M^0,M^1)$ such that
\mr
\item "{$(*)_4$}" $(a) \quad u_\ell$ is a subset of
$\lambda + \mu_\ell$ of cardinality $\mu_0$ 
satisfying $|u_\ell \cap \lambda| = \mu_0 =
|u_\ell \backslash \lambda|$ 
\nl

\hskip25pt and $[\lambda,\lambda + \mu_0) \subseteq u_\ell$ for $\ell=0,1$
\sn
\item "{${{}}$}" $(b) \quad M^\ell \prec N^{\mu_\ell}$ is primary
over $M \cup \{N_i:i \in u_\ell\}$ for $\ell = 0,1$
\sn
\item "{${{}}$}" $(c) \quad N^{\mu_\ell}$ is primary over $M^\ell
\cup\{N_i:i \in (\lambda + \mu_\ell) \backslash u_\ell\}$ for $\ell = 0,1$
\sn
\item "{${{}}$}" $(d) \quad f$ maps $M^0$ onto $M^1$.
\ermn
For $i \in (\lambda + \mu_\ell) \backslash u_\ell$ let $N_{\ell,i}
\prec N^{\mu_\ell}$ be primary over $M^\ell \cup N_{\{i\}}$ such that
$N^\ell$ is primary over $\cup\{N_{\ell,j}:j \in (\lambda + \mu_\ell)
\backslash u_\ell\}$; clearly
\mr
\item "{$(*)_5$}" for $\ell=0,1$
{\roster
\itemitem{ $(a)$ }  $M^\ell \prec N_{\ell,i}
\prec N^\ell$
\sn
\itemitem{ $(b)$ }  $(N_{\ell,i},c)_{c \in M^\ell} 
\cong (N_{\ell,j},c)_{c \in M^\ell}$ when $i,j
\in \lambda \backslash u_\ell$ or $i,j \in (\lambda + \mu_\ell)
\backslash \lambda \backslash u_\ell$
\sn
\itemitem{ $(c)$ }   $\langle N_{\ell,i}:i \in (\lambda +
\mu_\ell) \backslash u_\ell\rangle$ is independent over
$M^\ell$ and $N^\ell$ is primary over their union.
\endroster}
\ermn
Choose $\gamma_1 \in \lambda \backslash u_1,\gamma_2 \in
[\lambda,\lambda + \mu_1) \backslash u_1$, so
$(M^1,\text{tp}(N_{1,\gamma_1},M^1_\emptyset)$,
tp$(N_{1,\gamma_2},M^1_\emptyset))$ can serve as $(M_*,p_1,q_1)$ and
$r$ is $f(\text{tp}(M_{0,\gamma},M^0))$ for any $\gamma \in
\lambda \backslash u_0$.

So we have finished proving $\circledast_3$.]
\mr
\item "{$\circledast_4$}"  assume $\bar p,\bar q$ are sequences of members of
$\bold S^c_Y(M)$ and $(\bar p,\bar \alpha),(\bar q,\bar \beta)$ are
equivalent and $M \prec N$ and $p'_\varepsilon \in \bold S^c_Y(N)$
lift $p_\varepsilon$ for $\varepsilon < \ell g(\bar p)$ and
$q'_\varepsilon \in \bold S^c_Y(N)$ lift $q_\varepsilon$ for
$\varepsilon < \ell g(\bar q)$ \ub{then} $(\langle
p'_\varepsilon:\varepsilon < \ell g(\bar p)\rangle,\bar\alpha)$ and
$(\langle q'_\varepsilon:\varepsilon < \ell g(q)\rangle,\bar\beta)$
are equivalent.
\ermn
[Why?  By properties of ``primary".]
\mr
\item "{$\circledast_5$}"  if $p,q \in \bold S^{c,\theta}_Y(M)$ are
equivalent then $(p,\theta),(q,\theta)$ are equivalent.
\ermn
[Why?  By DLST.]

Note 
\mr
\item "{$\circledast_6$}"   in $\circledast_3$ we can conclude $p_2,r$
are equivalent.
\ermn
[Why?  In clause (d) of $\circledast_3$, let $N$ be the model and
let a witness for $N$ being
$(r,\lambda)$-constructed be $\langle N^r_i:i < \lambda\rangle$ and
for $N$ being $(\langle p_2,q_2\rangle,\langle \lambda,
\mu_1\rangle)$-constructed be $\langle N^*_i:i< \lambda + 
\mu_1\rangle$.  Let $u_0 \subseteq \lambda,u_1 \subseteq
\lambda + \mu_1$ be of cardinality $\mu_1,[\lambda,\lambda + \mu_1)
\subseteq u_1$ and $M'_*$ be such that:
\mr
\item "{$(*)_6$}"  $(a) \quad M' * \prec N$
\sn
\item "{${{}}$}"  $(b) \quad M'_*$ is primary over $\cup\{N^2_i:i \in u_1\}$
\sn
\item "{${{}}$}"  $(c) \quad M'_*$ is primary over $\cup\{N^*_i:i \in
u_2\}$
\sn
\item "{${{}}$}"  $(d) \quad N$ is primary over $\cup\{N^2_i:i \in
\lambda \backslash u_1\} \cup M'_*$
\sn
\item "{${{}}$}"  $(e) \quad N$ is primary over $\cup\{N^*_i:i \in
\lambda + \mu_1 \backslash u_2\}$.
\ermn
The liftings $r',p'_2$ of $r,p_2$ to $M'_*$ are equivalent, so we 
``collapse" to cardinality $\mu_0$ getting $M''_*$ so $M''_*$ is
$(r,\mu_0)$-constructed over $M_*$ and $(p_2,\mu_0)$-constructed over
$M_*$.  Then find liftings $r'',p''_2 \in \bold S^c_Y(M'_*)$ of 
$r,p$ respectively, so $r'',p''_2$ are equivalent naturally but 
$M''_*,M_*$ are isomorphic over $M$ by an isomorphism mapping
$r'',p''_2$ to $r,p_2$ so we get that $r,p_2$ are equivalent as required.]
\mr
\item "{$\circledast_7$}"   if $M'$ is $(p_2,\mu_0)$-constructed over
$M_*$ then $q_2$ is realized in $M'$.
\ermn
[Why?  Assume $M'$ is a counterexample.  Look again at the proof of
$\circledast_3$, so $M_*$ is $M^1$ there, and so $M'$ is
$(p_2,\lambda)$-constructed over $M^1 = M_*$ and $N^1$ is
$(p_2,\lambda)$-constructed over $M^1 = M_*$, so by uniqueness of
primary also in $N^1$ we cannot find $N' \prec N^1$ realizing $q_2$.
But for any $\gamma \in [\lambda,\lambda + \mu_1] \backslash u_2$ the
model $f^{-1}(N_{1,\gamma})$ contradict this.]

Now we can prove \scite{ss.3}.  Let $p,q \in 
\bold S^{c,\aleph_{\beta(*)}}_Y(M_\emptyset)$ be types which
$M_{\{0\}},M_{\{\omega\}}$ respectively realizes.  Let $(M,p_1,q_1) =
(M_\emptyset,p,q)$. 

So by $\circledast_7$, there are $\mu_0 < \lambda$ and $Y_1$ and $M_2
\in \bold L[T,Y_1],M_2$ which is $(\langle p,q\rangle,\langle
\mu_0,\mu_0\rangle)$-constructed over $M_\emptyset$ and
$p_2,q_2$ lifting of $p,q$ in $\bold S^{c,\mu_0}_{Y_1}(M_2)$ as there.
So by DLST we can find such $M'_2,p'_2,q'_2$ for the case 
$\mu_0 = \aleph_{\beta(*)}$, but this is absolute as $\beta(*)=0$.  
Also it gives the required result .   \hfill$\square_{\scite{ss.3}}$ 
\enddemo
\bigskip

\proclaim{\stag{ss.4} Theorem}  {[\rm ZF]}
If $T$ is countable and $\boxtimes_T$ below holds, \ub{then}
(recalling \scite{0.EA}(2))  
in every cardinal $\mu \ge \Upsilon({\Cal P}(\omega))$ 
we have $\dot I(\mu,T)$ is $\le |{\Cal F}^*/E|$ where ${\Cal F}^* = \{f:f$ a
function from ${\Cal P}(\omega)$ to $\omega +1\}$, for some equivalence
relation $E$ on the set of those functions where:
\mr
\item "{$\boxtimes_T$}"  $(a)-(d) \quad$ from \scite{ss.1.3}
\sn
\item "{${{}}$}" $(e)-(f) \quad$ from \scite{ss.2}
\sn
\item "{${{}}$}" $(g)$ \hskip30pt from \scite{ss.3}.
\endroster
\endproclaim
\bigskip

\remark{Remark}  1) Countability of $T$ is not used (if we write
${\Cal P}(|T|)$ instead of ${\Cal P}(\omega)$), but the gain is not
substantial.  This applies to \scite{ss.5}, \scite{ss.9}, too.
\nl
2) Fuller more accurate information is given in \scite{ss.9}.
\endremark
\bigskip

\demo{Proof}  Let $N$ be a model of $T$ of cardinality $\mu$ so
\wilog \, with universe $\mu$, we work in $\bold L[T,N]$ and we shall
analyze it.  Now we first choose a countable $M_\emptyset \prec
N_\ell$.  As $T$ is superstable, uni-dimensional we can 
find $\varphi(x,\bar y) \in \Bbb L(\tau_T)$ and $\bar a 
\in {}^{\ell g(\bar y)}(M_\emptyset)$ such that $\varphi(x,\bar a)$ is
weakly minimal.

We can find $\langle a_\alpha:\alpha < \mu \rangle$ such that:
\mr
\item "{$\circledast_1$}"  $(a) \quad a_\alpha \in
\varphi(N,\bar a_\ell) \backslash M_\emptyset$
\sn
\item "{${{}}$}" $(b) \quad \{a_\alpha:\alpha < \mu\}$ is
independent in $N$ over $M_\emptyset$ (in particular \nl

\hskip25pt with no repetitions)
\sn
\item "{${{}}$}" $(c) \quad$ modulo (a) + (b) the set
$\{a_\alpha:\alpha < \mu\}$ is maximal hence
\sn
\item "{${{}}$}" $(d) \quad \varphi(N,\bar a) \subseteq 
\text{ acl}(M_\emptyset \cup \{a_\alpha:\alpha < \mu\})$.
\ermn
Let $f \in \bold L[T,N]$ be a function from $\mu$ to $\mu$ such that
$f(\alpha) \le \alpha$ and $(\forall \beta < \mu)(\exists^\mu \alpha <
\mu)(f(\alpha) = \beta)$.  Now we try to choose
$(M_{\{\alpha\}},b_\alpha)$ by induction on $\alpha < \mu$ such that
\mr
\item "{$\circledast_2$}"  $(a) \quad b_\alpha \in
\varphi(N,\bar a)$
\sn
\item "{${{}}$}" $(b) \quad b_\alpha \notin \text{ acl}(M_\emptyset
\cup \{b_\beta:\beta < \alpha\})$
\sn
\item "{${{}}$}" $(c) \quad M_{\{\alpha\}} \prec N$ is $\bold
F^c_{\aleph_0}$-primary over $M_\emptyset \cup \{b_\alpha\}$, see
\cite[IV]{Sh:c} 
\sn
\item "{${{}}$}" $(d) \quad$ if $\alpha = 2 \beta +1$ and we can find
$(M_{\{\alpha\}},b_\alpha)$ satisfying (a) + (b) + (c) \nl

\hskip25pt and $(M^{N_\alpha}_{\{\alpha\}},c)_{c \in M_\emptyset} \cong
(M_{\{f(\beta)\}},c)_{c \in M_\emptyset}$ then
$(M_{\{\alpha\}},b_\alpha)$ satisfies this
\sn
\item "{${{}}$}"  $(e) \quad$ if $\alpha = 2 \beta$ and $\gamma_\alpha
= \text{ Min}\{\gamma:a_\gamma \notin \text{ acl}(M_\emptyset \cup
\{b_\varepsilon:\varepsilon < 2 \beta\})$ \nl

\hskip25pt then $b_\alpha = a_{\gamma_\alpha}$.
\ermn
[Why can we can carry the induction?   We can ignore clause (d) as
if its hypothesis hold, then clause (e) is irrelevant, and this
hypothesis says that we can fulfil clause (a),(b),(c),(d).  Also if
$\alpha = 2 \beta +1$ and the further assumption of (d) fail then we
can act as in clause (e).  Also in
all cases by cardinality considerations recalling $|\varphi(M,\bar a)|
= \|M\|$ by (c) of $\boxtimes_T$ of \scite{ss.1.3} there is $b_\alpha$
satisfying clauses (a) + (b) and if clause (e)'s assumption holds,
\wilog \, also its conclusion.

Let $B_\alpha = \text{ acl}_M(M_\emptyset \cup \{b_\alpha\})$.  By the choice
of $\varphi(x,\bar a)$ if $\varphi(x,\bar a) \in p \in \bold
S(B_\alpha,M)$ then \ub{either} $p$ forks over $\bar a$ hence is
algebraic hence realized in $B_\alpha$ \ub{or} $p$ does not fork over
$\bar a$ hence is finitely satisfiable in $M_\emptyset$.  Let $\langle
b_{\alpha,i}:i < i_\alpha\rangle$ be a maximal sequence of members of
$M$ such that for
each $i$ for some formula $\varphi(x,\bar c_{\alpha,i}) \in \text{
tp}(b_{\alpha,i},B_\alpha \cup\{b_{\alpha,j}:j<i\})$ hence no extension
in $\bold S(B_\alpha \cup \{b_{\alpha,j}:j<i\})$ forking over $\bar
c_{\alpha,i}$.  By \cite{Sh:31} there is $M_{\{\alpha\}} \prec N$ with
universe $B_\alpha \cup \{b_{\alpha,i}:i < i_\alpha\}$. 

So we are done.]
\mr
\item "{$\circledast_3$}"  $\langle b_\alpha:\alpha < \mu \rangle$
satisfies the requirements on $\langle a_\alpha:\alpha < \mu \rangle$.
\ermn
[Why?  Easy to check.]

So $\langle M_{\{\alpha\}}:\alpha < \mu \rangle$ is independent over
$M_\emptyset$ inside $N$ hence (by $\boxtimes_T(f))$ there is
$N' \prec N$ primary over $\cup\{M_{\{\alpha\}}:\alpha < \mu\}$
and by $\circledast_1(d)$ include $\varphi(N,\bar a)$ hence by
$\boxtimes_T(c)$ we have $N' = N$.

We can find a set $S$ and a partition $\langle I_t:t \in S \rangle$ of
$\mu$ such that: for $\ell=1,2$ and $\alpha,\beta < \mu$ we have
\mr
\item "{$\circledast_4$}" $(M_{\{\alpha\}},a_\alpha,c)_{c \in
M_\emptyset}$ is isomorphic to $(M_{\{\beta\}},a_\beta,c)_{c \in
M_\emptyset}$ iff $\dsize \bigvee_{t \in S} \{\alpha,\beta\} \subseteq I_t$.
\ermn
Now how large can $|S|$ be?  It is, in $\bold V, \le 
|{\Cal P}(\omega)|^{\bold L[T,N]} \le |{\Cal P}(\omega)|^{\bold V}$ 
(there is a function from a subset of ${\Cal P}(\omega)$ onto this
set).  So $|S| < \theta({\Cal P}(\omega))$, but $\bold L[T,N] \models
``ZFC + |S| \le 2^{\aleph_0} = |{\Cal P}(\omega)|"$ so there is a well
ordering of ${\Cal P}(\omega) \cap \bold L[T,N] = {\Cal
P}(\omega)^{\bold L[T,N]},|S| \le |{\Cal P}(\omega)^{\bold L[T,N]}| <
|S| < \Upsilon({\Cal P}(\omega))$ and
$\bold L[T,N] \models ``|S| \le 2^{\aleph_0}"$.  Now we shall prove:
\mr
\item "{$\circledast_5$}"  if $t \in S$ and $\bold L[T,N] \models
``\aleph_0 \le |I_t| < \mu"$ \ub{then} for some $\alpha < \mu$ we have
$(\forall s \in S)(|I_s \backslash \alpha| < \aleph_0)$ and $\mu < 
\Upsilon({\Cal P}(\omega))$.
\ermn
Clearly $\circledast_5$ helps because ``$\mu < \Upsilon({\Cal P}(\omega))$"
contradict an assumption on $\mu$".

Why $\circledast_5$ holds?  
Let $\alpha(*) = \text{ Min}(I_t)$, it is well defined as $I_t
\ne \emptyset$ because ``$\aleph_0 \le |I_t|$" was assumed.  
Let $J = \{2\beta +1:f(\beta) = \alpha(*)\}$.  If $2
\beta +1 \in J \Rightarrow 2 \beta +1 \in I_t$, then we get $|I_t| \ge
|\{2 \beta +1:f(\beta) = \alpha(*)\}| = \mu$ hence the assumption
``$|I_t| < \mu$"  is contradicted, so assume
that $\alpha = 2 \beta +1 \in J \backslash I_t$.  By clause (d) of
$\circledast_2$  apply
to $\alpha = 2 \beta +1$, we know that if $(N',b)$ satisfies the
demands on $(M_{\{\alpha\}},b_\alpha)$ in clauses (a),(b),(c) (i.e.,
$b_\alpha \in \varphi(N,\bar a) \backslash \text{ acl}(M_\emptyset
\cup \{b_\varepsilon:\varepsilon < \alpha\})$ and $N' \prec N$ is
$\bold F^c_{\aleph_0}$-primary over $M_\emptyset \cup \{b\})$ then
$(N',c)_{c \in M_\emptyset} \ncong (M_{\{\alpha(*)\}},c)_{c \in
M_\emptyset}$.  This implies that $I_t \subseteq \alpha$.  As it is
infinite, by $\boxtimes_T(g)$ we get $(\forall s \in S)(|I_s \backslash
\alpha| < \aleph_0)$ and recall $|\mu \backslash \alpha| = \mu$.  
So $\{\text{Min}(I_s \backslash \alpha):s \in
S$ and $I_s \nsubseteq \alpha\}$ is a subset of $\mu$ of cardinality $\mu$
(working in $\bold L[T,N]$) and there is a one-to-one mapping from it
into ${\Cal P}(\omega)$ (using the isomorphism types of
$(M_{\{\alpha\}},c)_{c \in M_\emptyset})$.  This gives $\mu <
\Upsilon({\Cal P}(\omega))$.  But this contradicts an assumption on $\mu$.

So we know
\mr
\item "{$\circledast_6$}"  if $I_t$ is infinite then it has cardinality $\mu$.
\ermn
Let $f =f_N = f_{N,M_\emptyset,\varphi(x,\bar a)}$ be the 
partial function from ${\Cal P}(\omega)$ into
$\omega +1$ defined as follows:  if $t \in S$ and 
$\eta \in {\Cal P}(\omega)$ codes
\footnote{see more details on this and similar points in the proof of
\scite{ss.9}}
 a model isomorphic to
$(M_\alpha,c)_{c \in M_\emptyset}$ for $\alpha \in I_t$ then
$f_N(\eta) = |I_t|$ if $I_t$ is finite and $f_N(\eta) = \omega$
otherwise: of course, the choice of $f_N$ is unique if we use the
canonical well
ordering of $\bold L[T,N]$ to make our choices in particular of
$M_\emptyset,\varphi(x,\bar a)$, but we could use ``any such $f$" so
increasing ${\Cal F}_\mu$ below (and fix the coding).

Now in $\bold V$ for any model $M$ of $T$ of cardinality $\mu$ we
define

$$
\align
{\Cal F}_M = \{(f_N,N \restriction \omega,\varphi(x,\bar a)):&N 
\text{ is a model with universe } \mu \text{ isomorphic to } M \\
  &\text{ such that } N \restriction \omega \prec N \text{ so can
serve as } M_\emptyset \\
  &\text{and } \bar a \in {}^{\omega >}N \text{ and } \varphi(x,\bar
a) \text{ is weakly minimal}\}
\endalign
$$

$$
\bold F^* = \cup\{{\Cal F}_M:M \text{ a model of } T \text{ of
cardinality } \mu\}.
$$
\mn
Clearly
\mr
\item "{$\circledast_7$}"  $(a) \quad {\Cal F}_M$ depends just on
$M/\cong$
\sn
\item "{${{}}$}"  $(b) \quad$ if ${\Cal F}_{M_1} \cap {\Cal F}_{M_2}
\ne \emptyset$ then $M_1 \approx M_2$ hence ${\Cal F}_{M_1} = {\Cal F}_{M_2}$
so there is an 
\nl

\hskip25pt equivalence relation $E_{T,\mu}$ on a subset of 
${\Cal F}^*$  such that the 
\nl

\hskip25pt ${\Cal F}_M$'s are its equivalence classes
\sn
\item "{${{}}$}"  $(c) \quad$ the number of models of $T$ in $\mu$ up
to isomorphism is equal to the
\nl

\hskip25pt  number of $E_{T,\mu}$-equivalence classes.
\ermn
So we are done.  \hfill$\square_{\scite{ss.4}}$
\enddemo
\bigskip

\proclaim{\stag{ss.5} Claim}  [$T$ countable] The 
demand $\boxtimes_T$ from \scite{ss.4} is absolute (property of $T$).
\endproclaim
\bigskip

\demo{Proof}  The new point is $\boxtimes_T(g)$ which should be clear.
\hfill$\square_{\scite{ss.5}}$
\enddemo
\bigskip

\remark{\stag{ss.8} Remark}  1) The proof of \scite{ss.4}, \scite{ss.5} is
really a particular case of ``the number of special dimensions" from
\cite[XIII,\S3]{Sh:c} the number being here 1; see more on this
Hrushovski Hart Laskowski \cite{HHL00}.
\nl
2) The ``primary over $\cup\{N_{\{\alpha\}}:\alpha\}$" is a special case of
decompositions. 
\endremark
\bigskip

\proclaim{\stag{ss.9} Theorem}  If $T$ is countable and $\boxtimes_T$
from \scite{ss.4} holds \ub{then}:  
\mr
\item "{$(a)$}"  $\dot I(\mu,T)$ is the same whenever $\mu \ge \mu_* =:
\theta({\Cal F}^*)$ recalling ${\Cal F}^* = \{f:f$ a function from
${}^\omega 2$ to $\omega$ with {\rm supp}$(f) = \{\eta:f(\eta) \ne 0\}$ well
orderable$\}$.
\endroster
\endproclaim
\bigskip

\demo{Proof}  We elaborate some parts done in passing in the proof of
\scite{ss.4} (and add one point).

We can interpret $\eta \in {}^\omega 2$ as a triple
$(M_0,M_1,\varphi(x,\bar a)) = (M^\eta_0,M^\eta_1,\varphi_\eta(x,\bar
a_\eta))$ such that $M_0 \prec M_1$ are models of
$T,M_1$ with universe $\omega,M_0$ with universe $\{2n:n < \omega\}$ and 
$\varphi(x,\bar a)$ a weakly minimal formula in $M_0$.
So the equivalence relation $E_1$ is
$\Sigma^1_1$ where $\eta E_1 \nu \Leftrightarrow [M^\eta_0 =
M^\nu_0,\varphi_\eta(x,\bar a_\eta) = \varphi_\nu(x,\bar a_\nu)$
and $M^\eta_1,M^\nu_1$ are isomorphic over $M^\eta_0 = M^\nu_0]$ and $E_0$
a Borel equivalent relation where 
$\eta E_0\nu \Leftrightarrow M^\eta_0 = M^0_\eta$.

Let

$$
\align
{\Cal P_1} = \{A:&A \subseteq {}^\omega 2 \text{ is not empty, any two
members are} \\
  &E_0 \text{-equivalent not } E_1 \text{ equivalent and} \\
  &A \text{ is well orderable}\}.
\endalign
$$
\mn
Let

$$
{\Cal F} = \{f:\text{ for some } A \in {\Cal P}_1,f \text{ is a
function from } A \text{ to } \omega +1 \text{ such that } \omega \in
\text{ Rang}(f)\}.
$$
\mn
Let

$$
\theta_* = \theta({\Cal F})(\le \theta({\Cal F}^*)).
$$
\mn
For $N$ a model of $T$ of cardinality $\ge \theta_*$ let ${\Cal F}_N
\subseteq {\Cal F}$ be defined as in the proof of \scite{ss.4} but we
can write $f$ and not $(f,M_\emptyset,\varphi(x,\bar a))$ as
$M_0,\varphi(x,\bar a)$ are determined by Dom$(f)$.  Let

$$
E^2_\mu = E^2_{T,\mu} =\{(f_1,f_2):\text{ there is } N \in \text{ Mod}_{T,\mu}
\text{ for which } f_1,f_2 \in{\Cal F}_N\}.
$$
\mn
(Recalling Mod$_{T,\mu} = \{M:M$ is a model of $T$ of
cardinality $\mu\}$).  

Now
\mr
\item "{$(*)_1$}"   if $\mu \ge \theta_*$ and $f \in {\Cal F}$ then
for some model $N$ of $T$ of cardinality $\mu$ we have $f \in {\Cal
F}_N$
\sn
\item "{$(*)_2$}"  if $N_1 \cong N_2$ are from $\text{Mod}_{T,\mu}$
and $\mu \ge \theta_*$ then ${\Cal F}_{N_1} = {\Cal F}_{N_2}$
\sn
\item "{$(*)_3$}"  if $N_1,N_2 \in \text{ Mod}_{T,\mu},\mu  \ge
\theta_*$ and ${\Cal F}_{N_1} \cap {\Cal F}_{N_2} \ne \emptyset$ then
$N_1 \cong N_2$
\sn
\item "{$(*)_4$}"  $E^2_\mu$ is an equivalence relation on ${\Cal F}$
\sn
\item "{$(*)_5$}"  ${\Cal F}_N$ for $N$ a model of $T$ of cardinality
$\ge \theta_*$ is an $E_2$-equivalence class
\sn
\item "{$(*)_6$}"  $E^2_\mu$ is the same for all $\mu \ge \theta_*$.
\ermn
[Why?  Assume $N^1,N^2$ are models of $T$ with universe $\mu,f_\ell
\in {\Cal F}_{N_\ell}$ and let
$N^i_\emptyset,a^\ell_\alpha,N^\ell_{\{\alpha\}}(\alpha < \mu)$ be as
in the proof of \scite{ss.4} exemplifying this.  Let $\theta_* \le \mu_1 <
\mu_2$.  If $\mu = \mu_1,f_1E^2_{T,\mu_2} f_2 \Rightarrow f_1
E^2_{T,\mu_1} f_2$ by the LS argument.  The other direction, i.e., if
$\mu = \mu_2$ is similar to the proof of \scite{ss.4}, i.e., we blow
up $\langle a_\alpha:\alpha \in I_t\rangle$ for some $t$ (or every
$t$) such that $|I_t| = \mu$ and continue as in \scite{ss.3}.]  
\hfill$\square_{\scite{ss.9}}$
\enddemo
\bigskip

\demo{\stag{ss.10} Conclusion}  For every countable complete first
order theory $T$, one of the following occurs
\mr
\item "{$(A)$}"   for every $\alpha,\dot I(\aleph_\alpha,T) \ge 
|\alpha|$,
in fact there is a sequence $\langle M_\beta:\beta < \alpha \rangle$
of pairwise non-isomorphic models of $T$ of cardinality
$\aleph_\alpha$
\sn
\item "{$(B)$}"   for all $\mu \ge \mu_* =: \theta({\Cal F}^*)$ (which
$\le \theta({}^{{\Cal P}(\omega)} \omega)),\dot I(\mu,T)$ is 
the same and has the form ${\Cal F}^*/E$ for some 
equivalence relation $E$ (see more in \scite{ss.10} and its proof).
\endroster
\enddemo
\bn
\margintag{ss.11}\ub{\stag{ss.11} Problem}  [ZF] Give complete classification of $\dot
I(\lambda,T)$ for $T$ countable by the model theoretic
properties of $T$ and the set theoretic properties of the universe.
\bn
But it maybe wiser to make less fine distinctions.
\definition{\stag{ss.11.1} Definition}  1) Let $|X| \precsim |Y|$ mean
that $X = \emptyset$ or there is a function from $Y$ onto $X$ (so $|X|
\le |Y|$ implies this).
\nl
2) Let $|X| \approx |Y|$ if $|X| \precsim |Y| \precsim |X|$ (so this
 weakens $|X| = |Y|$ and is an equivalence relation) and $|X|/\approx$
is called the essential power.
\enddefinition
\bn
\margintag{ss.12}\ub{\stag{ss.12} Thesis}:  It is most reasonable to interpret ``determining
$\dot I(\lambda,T)$" as finding $\dot I(\lambda,T)/\approx$ which is the
essential power 
$|\{M/\cong:M$ a model of $T$ with universe $\lambda\}|/\approx$.
\bigskip

\proclaim{\stag{ss.13} Claim}  Assume $\boxtimes_T$ of \scite{ss.4} and
$T$ is countable.
\nl
1) If $T$ is $\aleph_0$-stable \ub{then} 
$\dot I(\aleph_\alpha,T) = 1$ for every
$\alpha > 0$.
\nl
2) If $D(T)$ is uncountable and $\alpha > 0$ \ub{then}:
\mr
\item "{$(a)$}" $|\{A \subseteq {}^\omega 2:|A| \le \aleph_\alpha\}| \le
 \dot I(\aleph_\alpha,T)\}$
\sn
\item "{$(b)$}"  $\dot I(\aleph_\alpha,T) \text{ is } \underset \sim {}\to <
\text{-below } |\{A \subseteq {}^\omega 2:|A| \le \aleph_\alpha\}|$
\ermn
(note: $|A| \le \aleph_\alpha \Rightarrow A$ is well ordered).
\nl
3) If $D(T)$ is countable, $T$ is not $\aleph_0$-stable and there is a
set of $\aleph_1$ reals and $\alpha > 0$ \ub{then}

$$
\dot I(\aleph_\alpha,T) \approx |\{A \subseteq {}^\omega 2:|A| \le
\aleph_\alpha\}|.
$$
\endproclaim
\bigskip

\demo{Proof}  As in \cite{Sh:c}.  (E.g. in (2) the first inequality
holds as in $\bold L[T,Y]$ we can find countable complete $T_1
\supseteq T$ with Skolem functions $M_1 \models T_1,a_\eta \in {}^m
M_1$ for $\eta \in {}^\omega 2$ and $b_n \in M_1$ for $n < \omega$
such that letting $\alpha = \omega,A = ({}^\omega 2)^{\bold L[T,Y]}$
we have
\mr
\item "{$(*)^\alpha_A$}"  $(a) \quad \langle b_n:n < \alpha \rangle$
is a non-trivial indiscernible sequence in $M_1$ over 
\nl

\hskip25pt $\{\bar a_\eta:\eta \in A\}$
\sn
\item "{${{}}$}"  $(b) \quad \langle \text{ tp}(\bar
a_\eta,\emptyset,M_1 \restriction \tau_T):\eta \in {}^\omega 2 \rangle$
are pairwise distinct
\sn
\item "{${{}}$}"  $(c) \quad \langle \bar a_\eta:\eta \in {}^\omega 
2\rangle$ is indiscernible in $(M_1,b_n)_{n < \omega}$ in the weak sense of
\nl

\hskip25pt \cite[VII,\S2]{Sh:c}
\sn
\item "{${{}}$}"  $(d) \quad$ tp$(\bar a_\eta,\emptyset,M_1 \rest
\tau_T)$ is not realized in $M_1 \rest \text{ acl}(\{\bar a_\nu:\nu 
\in {}^\omega 2 \backslash \{\eta\}\} \cup \{b_n:n< \omega\})$.
\ermn
So in bigger universe this $M_1$ has a natural extension.
So we can define $M^+_1,\langle a_\eta:\eta \in ({}^\omega 2)^{\bold
V} \rangle,\langle b_\alpha:\alpha \in [\omega,\mu)\rangle$ naturally
such that $(*)^\mu_{{\Cal P}(\omega)}$ and define $M_A$ for $A
\subseteq {}^\omega 2$ as Sk$(\{\bar a_\eta:\eta \in A\} \cup 
\{b_\alpha:\alpha < \mu\},M_1)\} 
\restriction \tau_T$; if $A$ is well orderable then $M_A$ has
cardinality $\mu$.  \hfill$\square_{\scite{ss.13}}$
\enddemo
\bn
We now look at a well known example in our context.
\nl
\margintag{3c.38}\ub{\stag{3c.38} Example}:  There is countable stable, not superstable
$T$ with $D(T)$ countable such that: if
there are no sets of $\aleph_1$ reals then $\dot I(\aleph_\alpha,T)$
is ``manageable"
\mr
\item "{$(A)$}"    let $G$ be an infinite abelian group, each
element of order 2.  So ${}^\omega G$ is also such a group.  We define a
model $M$:
{\roster
\itemitem { $(a)$ }  its universe: $G \cup {}^\omega G$ (assuming $G \cap
{}^\omega G = \emptyset$)
\sn
\itemitem{ $(b)$ }  predicates $P^M = G,Q^M = {}^\omega G$
\sn
\itemitem{ $(c)$ }  the partial two-place function $H^M_1$ which is the
addition of $G$ (you may add $x \notin G \wedge y \notin G \Rightarrow x
+{}^My =x$)
\sn
\itemitem{ $(d)$ }  $H^M_2$ is the addition on ${}^\omega G$ (coordinatewise)
\sn
\itemitem{ $(e)$ }  a partial unary function $F^M_n$ such that
$\eta \in {}^\omega G \Rightarrow F^M_n(\eta) = \eta(n)$
\sn
\itemitem{ $(f)$ }  individual constants $c_1,c_2$ the zeroes of $G$
and ${}^\omega G$ respectively
\endroster}
\sn
\item "{$(B)$}"  let $T = \text{ Th}(M)$.  Let $K^* = \{N:N \models T$ and $N$
omit $\{Q(x) \wedge Q(y) \wedge F_n(x) = F_n(y) \wedge x \ne y:n <
\omega\}$
\sn
\item "{$(C)$}"  if 
$\langle M_t:t \in I \rangle$ is a sequence of models of $T$
we can naturally define their sum $\oplus_{t \in I} M_t$.  Clearly
$K^*$ is closed under sum (i.e., 
\nl

$|M| = P^M \cup Q^M$,
\nl

$P^M = \{f:f$ is a function with domain $I$ such that 
$f(t) \in P^{M_t}$ and $f(t)$
\nl

\hskip35pt  is the zero $c^{M_t}_1$ of the abelian group 
$(P^{M_t},H^{M_t}_1)$ 
\nl

\hskip35pt  for all but finitely many $t$'s$\}$,

$Q^M = \{g:g$ a function with domain $I,f(t) \in Q^{M_t}$ and for all
but finitely 
\nl

\hskip35pt many $t \in I$ we have $f(t) = c^{M_t}_1\}$
\nl

(we ignore that for $I$ finite, formally $P^M \cap Q^M \ne \emptyset\}$, etc.)
\sn
\item "{$(D)$}"   (ZF) If $M$ is a model from 
$K^*$ of cardinality $\lambda$ and $\lambda$ is a $(< \lambda)$-free 
cardinal (see Definition \scite{4.17A} below)
 \ub{then} $M = \underset{i<\lambda} {}\to \bigoplus M_i$ for some sequence
$\langle M_i:i < \lambda\rangle$ such that $i < \lambda \Rightarrow
\|M_i\| < \lambda$
\sn
\item "{$(E)$}"  in $(D)$ if the cardinal 
$\lambda$ is $(< \mu)$-free we can add 
$\|M_i\| < \mu$ (see Definition \scite{4.17A} below)
\sn
\item "{$(F)$}"  $(a) \quad$ define for 
$M \models T$ a two-place relation $E_M$ on $M$:
\nl

\hskip25pt  $a E_M b \Leftrightarrow (a=b) \vee (Q(a) \wedge Q(b) \wedge
\dsize \bigwedge_n F_n(x) = F_n(y))$.
\nl

\hskip25pt   It is an equivalence relation on $M$
\sn
\item "{${{}}$}"  $(b) \quad$ define $M/E_M$ naturally
\sn
\item "{$(G)$}"  $(a) \quad M/E_M \in K^*$ for any model $M$ of $T$
\sn
\item "{${{}}$}"  $(b) \quad M_1 \cong M_2 \Rightarrow M_1/E_{M_1} \cong
M_2/E_{M_1}$.
\sn
\item "{$(H)$}"  $(a) \quad$ if $M \models T$ 
and $a,b \in Q^M$ then $|a/E_M| = |b/E_M|$
\sn
\item "{${{}}$}"  $(b) \quad$ we can consider $a/E_M$ (where 
$M \models T,a \in Q^M$) an abelian group
\nl

\hskip25pt  $G_{a,M}$ with every element of order 2 except that the
zero is not given
\sn
\item "{${{}}$}"  $(c) \quad$ if $a,b \in Q^M$ then $G_{a,M},G_{b,M}$ as vector
spaces of $\Bbb Z/2 \Bbb Z$ without zero 
\nl

\hskip25pt has the same dimension
\sn
\item "{${{}}$}"  $(d) \quad$ call this dimension $\lambda(M)$
\sn
\item "{$(I)$}"  if $M_1,M_2$ are models of $T$ of cardinality
$\aleph_\alpha$ then $M_1 \approx M_2$ iff
$M_1/E_{M_2} \approx M_2/E_{M_2}$ and $\lambda(M_1) = \lambda(M_2)$
hence $\dot I(\aleph_\alpha,T) = \dot I(\le \aleph_\alpha,K^*) 
\times |\omega + \alpha|$ where $\dot I(\le \aleph_\alpha,K^*) =
\Sigma\{\dot I(\aleph_\beta,K^*):\beta \le \alpha\} = \Sigma\{\dot
I(\aleph_\beta,K^*):\beta \le \alpha,\aleph_\beta < \theta({\Cal
P},\omega)\}$.  \hfill$\square_{\scite{3c.38}}$
\endroster
\newpage

\head{\S4 On $T$ categorical in $|T|$} \endhead  \resetall \sectno=4
 \spuriousreset
\bigskip

The ZFC parallel of \scite{1} - \scite{3} is the known ``$|D(T)|
< |T|$ implies $T$ is the definitional extension of some $T' \subseteq
T,|T'| < |T|$", see Keisler \cite{Ke71a}, 
which in Boolean algebra terms say ``the number of
ultrafilters of an infinite Boolean algebra $B$ is $\ge |B|$".
\bigskip

\demo{\stag{0} Convention}   For \scite{1}-\scite{5}, $T$ is first
order (with $\tau_T$ not necessarily well-orderable).
\enddemo
\bigskip

\definition{\stag{1} Definition}  For a first order $T$
in the vocabulary $\tau = \tau_T$, usually for simplicity closed 
under deduction, we define the equivalence relation $E_T$ on $\tau_T$ by
\mr
\item "{$(a)$}"  \ub{for predicates $P_1,P_2 \in \tau$}
\nl
$P_1 E P_2$ \ub{iff}: $P_1,P_2 \in \tau$ are predicates with
the same arity and $(\forall \bar x)(P_1(\bar x) \equiv P_2(\bar x)) \in T$
\sn
\item "{$(b)$}"  \ub{for function symbols $F_1,F_2 \in \tau$},
e.g. individual constants
\nl
$F_1 E F_2$ \ub{iff}: $P_1,P_2 \in \tau$ are function symbols with
the same arity and $(\forall \bar x)(F_1(\bar x) = F_2(\bar x)) \in T$
\sn
\item "{$(c)$}"  no predicate 
$P \in \tau$ is $E$-equivalent to a function symbol $F \in \tau$.
\endroster
\enddefinition
\bigskip

\definition{\stag{2} Definition}  Let $\tau = \tau_T,T$ a first order
theory.
\nl
1) The theory $T$ is called reduced if $E_T$ is
the equality.
\nl
2) Let $\tau/E_T$ be the vocabulary with predicates $P/E_T,P \in
\tau$ a predicate with arity$(P/E_T) = \text{ arity}_\tau(P)$ and similarly
$F/E_T$.
\nl
3) For a $\tau$-model $M$ of $T$ we define $M^{[E_T]}$ naturally, i.e.,
\nl
$N=M^{[E_T]}$ iff they have the same universe, $N$ is a
$(\tau/E_T)$-model, $M$ is a $\tau$-model $M \models T$ and
$(R/E_T)^N = R^M$ for every predicate $R \in \tau(T)$ and $(F/E_T)^N =
F^M$ for any function symbol $F \in \tau(T)$.
\nl
4) For $N$ a $(\tau/E_T)$-model, $M = {}^{[E_T]}N$ is the $\tau$-model
such that $N = M^{[E_T]}$ if one exists.
\nl
5) Let $T/E_T$ be the set of $\psi \in \Bbb L(\tau(T)/E_T)$ such that
if we replace any predicate $R/E_T$ appearing in $\psi$ by some $R'
\in R/E_T$ and similarly for $F/E_T$, we get a sentence from 
$\{\psi \in \Bbb L(\tau_T):T \vdash \psi\}$, see \scite{1}.
\enddefinition
\bigskip

\demo{\stag{3} Observation}  {\rm [ZF]}  For every first order $T$ (as
in \scite{1}) in $\Bbb L(\tau_T)$
\mr
\item "{$(a)$}"  $E_T$ is an equivalence relation on $\tau$
\sn
\item "{$(b)$}"  if $M$ is a $\tau$-model of $T$ then $M^{[E_T]}$ is a
uniquely determined $(\tau/E_T)$-model of 
$T/E$ and ${}^{[E_t]}(M^{[E_T]}) = M$
\sn
\item "{$(c)$}"  for every ($\tau/E$)-model $M$ of $T/E_T$ the
$\tau$-model ${}^{[E_t]}M$ uniquely determined and 
is a model of $T$ and $({}^{[E_T]}M)^{[E_T]} = M$
\sn
\item "{$(d)$}"   $T/E_T$ is a reduced first order theory
\ermn
See hopefully more on such $T$'s in \cite{Sh:F701}.
\enddemo
\bigskip

\demo{\stag{cat.0} Hypothesis}  $\tau(T) \subseteq \bold L$ as usual.
\enddemo
\bigskip

\proclaim{\stag{5} Claim}   {\rm [ZF]} If $T$ is a complete 
first order theory in $\Bbb L(\tau)$ and $T$ is reduced and $Y
\subseteq \bold L$, \ub{then} $T \in \bold L[Y] \Rightarrow
|D(T)|^{\bold L[T]} \ge |T|$.
\endproclaim
\bigskip

\demo{Proof}  By the ZFC case (see Keisler \cite{Ke71a}).  
\hfill$\square_{\scite{5}}$
\enddemo
\bigskip

\proclaim{\stag{cat.1} Claim}  If $T \subseteq \bold L$ 
is categorical in $\lambda$ and $Y \in
{ \text{\rm Ord\/}}$ \ub{then} in $\bold L[T,Y]$ the following is impossible
\mr
\item "{$\circledast$}"  $(a) \quad T$ stable, $\lambda \ge |T| +
\aleph_1 + \mu$,
\sn
\item "{${{}}$}"  $(b) \quad M \prec {\frak C}$ is 
$\bold F^f_{\aleph_0}$- primary over $\emptyset$, see \cite[IV]{Sh:c}
\sn
\item "{${{}}$}"  $(c) \quad \bar a_i \in {}^n M$ for $i < \mu$
\sn
\item "{${{}}$}"  $(d) \quad$ {\rm tp}$(\bar a_\delta,\cup\{\bar a_i:i
< \delta\})$ forks over $\cup \{a_j:j < \alpha\}$ whenever \nl

\hskip20pt $\alpha < \delta < \mu,\delta$ a limit ordinal from $S$
\sn
\item "{${{}}$}"  $(e) \quad$ every type over $\cup\{\bar a_i:
i < \mu\}$ which is realized in $M$ does not fork \nl

\hskip25pt over some $\cup\{\bar a_i:i < \alpha\}$ for some $\alpha < \mu$
\sn
\item "{${{}}$}"  $(f) \quad$ in $\bold L[T,Y]$ we have:
\nl

\hskip25pt $\mu$ regular uncountable, $S \subseteq \mu$ stationary.
\endroster
\endproclaim
\bigskip

\demo{Proof}  Work in $\bold L[T,Y]$; \wilog \, $M$ has cardinality
$\lambda$, and toward contradiction assume
$\circledast$ holds.  By clause (b) there is $\bar{\bold c}$ such that
\mr
\item "{$(*)_1$}"  $\bar{\bold c} = 
\langle \bar c_i: i < i^* \rangle$ 
\sn
\item "{$(*)_2$}"  $M = \cup\{\bar c_i:i < i^*\}$ and tp$(\bar
c_i,\cup\{\bar c_j:j < i\})$ does not fork over 
some finite $B_i \subseteq \cup\{\bar c_j:j<i\}$ for each $i<i^*$
\ermn
So by the properties of non-forking (or of $\bold
F^f_{\aleph_0}$-constructions, \cite[IV]{Sh:c}) without loss of
generality we have ($i^* \ge \mu$ and) $\cup\{\bar a_i:i < \mu\} \subseteq 
\cup\{\bar c_j:j < \mu\}$.  Hence for some club $E$ of $\mu$ we have $\bar a_i
\subseteq \dbcu_{j < \delta} \bar c_j \Leftrightarrow i < \delta$ for
$i < \mu,\delta \in E$; clearly
tp$(\bar a_\delta,\cup\{\bar c_i:i < \delta\})$ does
not fork over some finite $C_\delta \subseteq \cup\{\bar c_j:j <
\delta\}$.  Hence there is stationary $S_1 \subseteq S \cap C$ 
such that $\delta \in
S \Rightarrow C_\delta = C_*$, and let $\bar c$ list $C_*$.

By clause (e) of the assumption for some $\alpha_* < \mu$,
\mr
\item "{$(*)_2$}"   tp$(\bar c,\cup\{\bar a_i:i < \mu\})$ does not
fork over $\cup\{\bar a_i:i < \alpha_*\}$
\ermn
hence by the non-forking calculus
\mr
\item "{$(*)_3$}"   for $\delta \in S_1 \backslash (\alpha +1)$ the type
tp$(\bar a_\delta,\cup\{\bar a_i:i < \delta\})$ does not fork over
$\cup\{\bar a_i:i < \alpha_*\}$.
\ermn
By this contradicts clause (d) of the
assumption.   \hfill$\square_{\scite{cat.1}}$ 
\enddemo
\bigskip

\proclaim{\stag{cat.2} Claim}  If $T$ is stable, categorical in
$\lambda$ and $\lambda = |T| > \aleph_0$ \ub{then}
\sn
\ub{Case $(\alpha)$}:  if $(\exists Y \subseteq \text{\rm Ord})
(\aleph_1 = \aleph^{\bold L[T,Y]}_1)$ then $\kappa_r(T) = \aleph_0$,
i.e. $T$ is superstable.
\sn
\ub{Case $(\beta)$}:  if $(\forall Y \subseteq \text{\rm Ord})
(\aleph_1 > \aleph^{\bold L[T,Y]}_1)$ then for every $Y 
\subseteq \text{\rm Ord}$ we have $\bold L[T,Y] \models
``\kappa(T) < \aleph^{\bold V}_1"$.
\endproclaim
\bigskip

\demo{Proof}  \ub{Case $(\alpha)$}:

Assume the conclusion fails.  Fix $Y \subseteq \text{ Ord}$ such that
$T \in \bold L[Y],\aleph_1 = \aleph_1^{\bold L[Y]}$ and ${\frak C} =
{\frak C}_T \in \bold L[Y]$ is a $\chi$-saturated (in $\bold L[Y]$) 
model of $T$ and $\bold L[Y]
\models \kappa(T) \ge \aleph_1$ where $\chi$ is large enough and
regular in $\bold L[Y]$; and we shall work inside $\bold L[Y]$.
\nl
Let $\mu = \aleph_1^{\bold L[Y]} = \aleph^{\bold V}_1$.  
We can find $\langle \bar a_n:n < \omega\rangle,
\bar a_n \in {}^{\omega >}{\frak C}_Y$ and a type $p =
\{\varphi_n(x,\bar a_n):n < \omega\}$ such that $\varphi_n(x,\bar
a_n)$ forks over $\cup\{\bar a_m:m < \omega\}$. Let $\langle \eta_i:i
< \mu\rangle$ list ${}^{\omega >}(\mu)$ such that $\eta_i
\triangleleft \eta_j \Rightarrow i <j$ and for every limit ordinal
$\delta < \mu$ we have ${}^{\omega >} \delta =
\{\eta_i:i < \delta\}$.  We choose $\bar \nu = \langle
\nu_\delta:\delta < \mu$ limit$\rangle$ such that
$\nu_\delta$ is increasing with limit $\delta$.

We choose $\langle \bar a_\eta:\eta \in {}^{\omega >}\mu\rangle$
such that:  if $\ell g(\eta)=n$ then $\hat a_{\eta
\restriction 0} \char 94 \bar a_{\eta \restriction 1} \char 94 \ldots
\char 94 \bar a_\eta$ and $\bar a_0 \char 94 \ldots \char 94 \bar a_n$
realize the same type in ${\frak C}$ and tp$(\bar a_\eta,\cup\{\bar
a_\nu:\nu \in {}^{\omega >} \mu$ and 
$\neg(\eta \trianglelefteq \nu)\})$ does not fork over $\cup\{\bar
a_{\eta \restriction k}:k < \ell g(\eta)\}$.  For limit $\delta <
\mu$ we choose $b_\delta$ which realizes
$\{\varphi_n(x,\bar a_{\nu_\delta \restriction n}):n < \omega\}$ such
that tp$(b_\delta,\cup\{\bar a_\nu:\nu \in{}^{\omega >} \mu\} \cup
\{b_{\delta'}:\delta' < \delta$ is limit$\})$ does not fork over
$\cup\{\bar a_{\nu_\delta \restriction n}:n < \omega\}$. 

Lastly, let $a'_i$ be $b_\delta$ if $i = \delta$ and be $\bar a_{\eta_j}$ if
$i=j+1$, and be $<>$ if $i=0$.

Let $M_1 \prec {\frak C},M_1 \in \bold L[Y]$ be a 
model of cardinality $\lambda$ which is
$\bold F^f_{\aleph_0}$-primary over $\emptyset$.

Let $M_2 \prec {\frak C}_T,M_2 \in \bold L[Y]$ be 
$\bold F^f_{\aleph_0}$-primary over $\cup\{\bar a'_i:
i < \mu\}$ of cardinality $\lambda$ (see \cite[IV]{Sh:c}).  Now
by \scite{cat.1} for $\mu = \aleph_1$, the models $M_1,M_2$ 
are not isomorphic even in $\bold L[Y,Y_1]$ for any $Y_2 \subseteq
\text{\rm Ord}$ (as $\aleph_1^{\bold L[Y,Y_0]} = \aleph_1^{\bold L[Y]}
= \aleph^{\bold V}_1$), contradiction.
\sn
\ub{Case $(\beta)$}:  Assume that the conclusion fails for $Y$.  Clearly
$\aleph^{\bold V}_1$ is a limit cardinal in $\bold L[T,Y']$ for 
every $Y' \subseteq \text{ Ord}$.  
So for every $\mu \in \text{ Card}^{\bold L[T,Y]} \cap
\omega^{\bold V}_1$ we can find (in ${\frak C}^{\bold L[T,Y]}_T \in
\bold L[T,Y]$ chosen as above) a sequence 
$\bar{\bold a}_\mu = \langle \bar a_{\mu,i}:i < \mu \rangle$ such that
$\bar a_{\mu,i} \in {}^{\omega>}{\frak C}$ for $i < \mu$ and
 a type $p=\{\varphi_{\mu,i}(x,\bar a_{\mu,i}):i < \mu\}$ in ${\frak
C}$ such that $\varphi_i(x,\bar a^\mu_i)$ forks over $\cup\{\bar a_{\mu,j}:j <
i\}$ for every $i$.  Choose by induction on $i < \mu$ an element
$b^\mu_i \in {\frak C}$ which realizes $\{\varphi_{\mu,j}
(x,\bar a_{\mu,j}):j <i\}$ but tp$(b^\mu_i,\cup\{\bar
a_{\mu,j}:j < \omega_1\}\cup\{b^\mu_j:j<i\})$ does not fork over 
$\cup\{\bar a_{\mu,j}:j<i\}$.  
Let $\bar a^\mu_i = \bar a_{\mu,i} {}^\frown \langle b^\mu_i\rangle$ 
so $\langle \bar a^\mu_i:i < \mu\rangle$ is as in clauses (c) + (d) of 
\scite{cat.1}.  Note that the function $(\mu,i)
\mapsto \bar a^\mu_i$ belongs to $\bold L[T,Y]$.  
Without loss of generality $\{\bar{\bold a}_\mu:\mu
\in \text{ Card}^{\bold L[T,Y]} \cap \omega^{\bold V}_1\}$ is independent over
$\emptyset$ in ${\frak C}^{\bold L[T,Y]}_Y$.  In $\bold L[T,Y]$ let 
$M_1 \prec {\frak C}^{\bold L[T,Y]}_T$ be of
cardinality $\lambda,\bold F^f_{\aleph_0}$-primary over
$\emptyset$.  Let $M_2 \prec {\frak C}^{\bold L[T,Y]}_T$ be of
cardinality $\lambda$ and $\bold F^f_{\aleph_0}$-primary over
$\cup\{\bar{\bold a}_\mu:\mu \in \text{ Card}^{\bold L[T,Y]} \cap
\omega^{\bold V}_1\}$.  But $T$ is categorical in $\lambda$ so there
is an isomorphism $\bold f \in \bold V$ from $M_1$ onto
$M_2$ and now we shall work in $\bold L[T,Y,f]$ and let $\mu_* = \aleph^{\bold
L[T,Y,\bold f]}_1$, clearly $\mu_* \in \text{ Reg}^{\bold L[T,Y]} \cap
\omega^{\bold V}_1$ so $\bar{\bold a}_{\mu_*}$ is well defined.  By the
non-forking calculus, the statement $\circledast$ of \scite{cat.1}
holds for $\mu_*$ so we are done.  \hfill$\square_{\scite{cat.2}}$
\enddemo
\bigskip

\remark{\stag{cat.2.3} Remark}  Assume $T$ is stable, (complete with
infinite models of course), 
$\lambda = |T| \ge \aleph_\alpha > \aleph_0$ and for some $Y
\subseteq \text{\rm Ord}$ we have $\bold L[Y] \models ``\kappa(T) >
\aleph_\alpha$ or $\aleph_\alpha$ is a limit cardinal and $\kappa(T)
\ge \aleph_\alpha"$.  \ub{Then} $\dot I(\lambda,T) \ge |\alpha|$.  The proof
is similar.
\endremark
\bigskip

\proclaim{\stag{cat.3} Claim}  $T$ is not categorical in $\lambda = |T|
> \aleph_0$ \ub{when} for some $Y \subseteq { \text{\rm Ord\/}}$:
\mr
\item "{$\circledast$}"  $(a) \quad T$ is stable
\sn
\item "{${{}}$}"  $(b) \quad \bold L[T,Y] \models ``|D(T)| \ge \lambda =
|T|"$ (holds if $T$ is reduced, see \scite{5})
\sn
\item "{${{}}$}"  $(c) \quad$ the conclusion of \scite{cat.2} holds,
(or just for every $Y' \subseteq \text{\rm Ord}$ we have
\nl

\hskip25pt $\lambda > \kappa(T)^{\bold L[Y',Y,T]}$).
\endroster
\endproclaim
\bigskip

\demo{Proof}  Choose $Y \subseteq \text{ Ord}$ which exemplify
the assumption of case $(\alpha)$ of \scite{cat.2} if it holds.
In $\bold L[T,Y]$ letting $\kappa = \kappa(T)^{\bold L[T,Y]}$ let:
\mr
\item "{$(*)_1$}"  $M_1$ be $\bold F^a_\kappa$-constructible over
$\emptyset$ of cardinality $\lambda$, i.e., for some sequence $\langle
a_i,B_i:i < \lambda \rangle$ we have $M_1 = \{a_i:i < \lambda\}$ and $B_i
\subseteq \{a_j:j<i\}$ has cardinality $< \kappa$ and 
stp$(a_i,B_i) \vdash \text{ stp}(a_i,\{a_j:j<i\})$ 
(not necessarily $\bold F^a_\kappa$-saturated!) 
\sn
\item "{$(*)_2$}"  $M_2$ be a model of $T$ of cardinality $\lambda$
with $\bold I \subseteq M_2$ indiscernible of cardinality $\lambda$.
\ermn
[Why $(*)_2$ is possible?  E.g., we can have $\|M_1\| = \lambda$ because 
$\bold L[T,Y] \models ``|D(T)| \ge \lambda"$.]

So assume toward
contradiction that $M_1,M_2$ are isomorphic, let $\bold f:M_1
\underset{\text{onto}} {}\to {\overset\text{iso} {}\to
\longrightarrow} M_2$ be such an
isomorphism and work in $\bold L[T,Y,\bold f]$.  Now $\kappa(T)^{\bold
L[T,Y,\bold f]}$ may be $> \kappa = \kappa(T)^{\bold L[T,Y]}$ and
$\kappa$ may be not a cardinality still
the properties of $M_1,M_2$ from $(*)_1,(*)_2$ respectively holds in
$\bold L[T,Y,\bold f]$ for $\kappa = \kappa(T)^{\bold L[T,Y,\bold f]}$. 
Now we can get a contradiction as in  \cite[IV]{Sh:c}.    
\hfill$\square_{\scite{cat.3}}$  
\enddemo
\bn
Putting together Claims \scite{cat.2}, \scite{cat.3}. 
\demo{\stag{cat.7N} Conclusion}  If $T$ is stable in $\lambda = |T| \le
|D(T)|$ \ub{then} $T$ is not categorical in $\lambda$.
\enddemo
\bn
\centerline{$* \qquad * \qquad *$}
\bigskip

\demo{\stag{7.0} Free Models}  Let $T$ be complete and stable.  ${\frak C}
= {\frak C}_{Y,T}$ a monster for $T$ in $\bold L[T,Y]$.

The proofs above (and actually \cite{Sh:c}) 
suggest that we look more into free models.
\enddemo
\bigskip

\definition{\stag{7.1} Definition}  1) A model $M$ of $T$, (a stable
theory) is free \ub{when} we
can find a sequence $\langle a_i:i < \alpha \rangle$ enumerating $M$
such that for each $i < \alpha$ the type tp$(a_i,\{a_j:j<i\},M)$ does
not fork over some finite subset say $B_i$.
\nl
2) We call $\langle (A_i,a_i,B_i):i < \alpha\rangle$ is a free
   representation of $M$ where $A_i = \{a_j:j<i\}$.
\enddefinition
\bigskip

\remark{Remark}  So free is the same as being $\bold
F^f_{\aleph_0}$-constructible over $\emptyset$.
\endremark
\bigskip

\proclaim{\stag{7.2} Claim}  If $A \subseteq {\frak C},\lambda 
= |A|$ is singular and every $A'
\subseteq A$ of cardinality $< \lambda$ is free \ub{then} $A$ is free.
\endproclaim
\bigskip

\demo{Proof}  By compactness in singular (\cite{Sh:54},
\cite{Sh:E18}).  \hfill$\square_{\scite{7.2}}$
\enddemo
\newpage

\head {\S5 Consistency results} \endhead  \resetall \sectno=5
 \spuriousreset
\bigskip

In spite of the evidence of \S1,\S4, without choice characterization
for the number of non-isomorphic models is different without choice.
We look for consistency results for ``there are few models in cases
impossible by ZFC", in particular
we ask (and give a partial answer):
\nl
\margintag{4.1}\ub{\stag{4.1} Question}:  1) Is it consistent with ZF that 
for some/many $\kappa > \aleph_0$ we have: every two strongly
$\aleph_0$-homogeneous linear orders of cardinality $\kappa$, are
isomorphic?  (Add ``$\kappa$ singular or $\kappa$ regular"; or add
cf$(\kappa) = \aleph_0$.)
\nl
2) Similarly is it consistent with ZF that

``if $M_1,M_2 \subseteq ({}^\omega \lambda,E_n)_{n <
\omega}$ are strongly $\aleph_0$-homogeneous of cardinality $\kappa$ then they
are isomorphic".
\nl
3) Instead categoricity proves the consistency of all models has
nice descriptions, (see below):

Clearly \scite{4.21} below proves that our use of elementary classes in
the proof for stable, un-superstable $T$ 
is necessary, that is we could not prove too good theorems on PC
classes parallel to the ZFC case.
\bn
Toward \scite{4.1}(2) we consider:
\definition{\stag{4.17A} Definition}  1) A cardinality $\lambda$ is
free or $\omega$-sequence-free when every subset of
${}^\omega \lambda$ of cardinality $\lambda$ is free, where
\nl
2) A subset $A \subseteq {}^\omega \lambda$ is free when there is a
one-to-one function $f:A \rightarrow {}^{\omega >}\lambda$ 
such that $\eta \in A \Rightarrow f(\eta) \triangleleft \eta$.
\nl
3) A cardinal $\lambda$ is $(< \mu)$-free when every subset of
${}^\omega \lambda$ of cardinality $\le \lambda$ is $(< \mu)$-free where
\nl
4)  We say $``A \subseteq {}^\omega \lambda$ is $(< \mu)$-free if
there is a function $f:A \rightarrow {}^{\omega >}\lambda$ which in
some $\bold L[Y]$ is 
$(< \mu)$-to-one" and $\eta \in A \Rightarrow f(\eta) \triangleleft \eta$.
\enddefinition
\bn
\margintag{4.17B}\ub{\stag{4.17B} Question}:  1) Is it consistent (with ZF) that for
arbitrarily large $\mu,\mu^+$ is $\mu^+$-free?  ($\aleph_0$ always is).
\nl
2) Is it consistent with ZF that all cardinals are free?
\bigskip

\proclaim{\stag{4.17} Claim}  [ZF + DC]  Let $\kappa = \aleph_1$.
The following is a sufficient condition
for $\lambda$ being $(< \kappa)$-free (equivalently - free)
\mr
\item "{$\boxdot_{\lambda,\mu}$}"  for every $A \subseteq \lambda$ for some $B
\subseteq \lambda$ we have:
{\roster
\itemitem{ $(*)_1$ }   if $\bold L[A] \models ``\mu$ is a
cardinal $< \lambda$ but $\ge \kappa$ such that $\mu < \mu^{\aleph_0},\mu' =
{ \text{\rm Min\/}}\{\lambda,\mu^{\aleph_0}\}"$ then $\bold L[A,B] \models
``\mu'$ is an ordinal of cardinality $\le \mu"$
\sn 
\itemitem{ $(*)_2$ }  if $\bold L[A] \models ``\mu \le \lambda$ is
regular uncountable $\ge \kappa$ and 
$S = \{\delta < \lambda:\text{\rm cf}(\delta) =
\aleph_0\}"$ then $\bold L[A,B] \models ``S$ is a non-stationary subset
of $\mu"$.
\endroster}
\endroster
\endproclaim
\bigskip

\remark{Remark}  1) A condition for $\kappa > \aleph_1$ will be more
complicated.
\nl
2) $(< \aleph_1)$-free is equivalent to free (note that ``in some
$\bold L[Y]$" in Definition \scite{4.17A}).
\endremark
\bigskip

\demo{Proof}  So assume that $A$ is a subset of ${}^\omega \lambda$ of
cardinality $\le \lambda$.

Let

$$
\align
\Xi = \{(Y,f):&Y \subseteq \text{ Ord and } f \in \bold L[Y] 
\text{ is a function from } A \text{ to }
{}^{\omega >}\lambda \\
  &\text{such that } \eta \in A \Rightarrow f(\eta) \triangleleft \eta\}
\endalign
$$
\mn
and let $\bar \mu_{Y,f} = \langle \mu^{Y,f}_{\bar \eta}:\eta \in
A\rangle$ be defined by $\mu^{Y,f}_\eta = |\{\eta' \in A:f(\eta') =
f(\eta)\}|^{\bold L[f,Y]}$.

So the role of $Y$ is in determining where we compute $\mu^{Y,f}_\eta$.

Now it suffices to prove
\mr
\item "{$\circledast$}"  if $(Y,f) \in \Xi$ then there is $(Z,g) \in
\Xi$ such that $\eta \in A \Rightarrow \bar \mu^{Y_1,f_1}_\eta <
\mu^{Y,f}_\eta \vee \mu^{Y,f}_\eta < \kappa$.
\ermn
[Why it suffices?  If so by DC we can find $\langle (Y_n,f_n):n < \omega\rangle
\in V$ such that $(Y_n,f_n) \in \Xi$ and
\mr
\item "{$(*)$}" $\eta \in A \Rightarrow (\mu^{Y_n,f_n}_\eta >
\mu^{Y_{n+1},f_{n+1}}_\eta) \vee (\mu^{Y_n,f_n}_\eta < \kappa)$.
\ermn
Let $Y_* = \{\text{cd}(\langle 1,n,\ell g(\eta)\rangle \char 94 \eta
\char 94 \langle f_n(\eta)\rangle):\eta \in {}^{\omega >}\lambda$ and
$n < \omega\} \cup \{\text{cd}(2,n,\alpha):\alpha \in Y_n$ and $n <
\omega\}$ where cd is a one-to-one definable function in $\bold L$
from ${}^{\omega >}\text{Ord}$ into Ord.

Clearly $\langle f_n:n < \omega\rangle \in 
\bold L[Y_*]$ and define $h:A
\rightarrow \omega$ by $h(\eta) = \text{ Min}\{n:\mu^{Y_n,f_n}_\eta <
\kappa\}$, it clearly exists by $\circledast$.

Lastly, let $f:A \rightarrow {}^{\omega >} \lambda$ be defined by
$f(\eta) = \eta \restriction \text{ pr}(h(\eta),\ell
g(f_{h(\eta)}(\eta))$ where pr$(n,m)$ is, e.g. $(n + m+1)^2 + n$.

Now check.]
\enddemo
\bigskip

\demo{Proof of $\circledast$}:  Let $Z$ be like $B$ is the claim's
assumption with $Y$ playing the roles of $A$; we work in $\bold
L[Y,Z]$, \wilog \, $Y \in \bold L[Z]$.  Let $\langle
\eta_\alpha:\alpha < |A|\rangle$ list $A$ with no repetitions.  Let
${\Cal U} = \{\alpha < |A|$: for no $\beta < \alpha$ do we have
$f(\eta_\beta) = f(\eta_\alpha)\}$ and let $\mu_\alpha =
|\{\beta:f(\eta_\beta) = f(\eta_\alpha)\}$ for $\alpha \in {\Cal U}$
and let (so $\langle \mu_\alpha:\alpha \in {\Cal U}\rangle \in \bold
L[Y]$.

In $\bold L[Y]$ let $\left< \langle \eta_{\alpha,\varepsilon}:\varepsilon <
\mu_\alpha\rangle:\alpha \in {\Cal U}\right>$ be such that for each
$\alpha \in {\Cal U}$ the sequence $\langle
\eta_{\alpha,\varepsilon}:\varepsilon < \mu_\alpha\rangle$ list
$A_\alpha := \{\beta:f(\beta) = f(\alpha)\}$.  Now
\mr
\item "{$\boxtimes$}"  it suffices to prove that in $\bold L[Z]$,
for every $\alpha \in {\Cal U}$ there is $f_\alpha:A_\alpha
\rightarrow {}^{\omega >} \lambda$ such that $\eta \in A_\alpha
\Rightarrow |\{\nu \in A_\alpha:f_\alpha(\nu) = f_\alpha(\eta)\}| < 
\mu_\alpha$.
\ermn
Note that in $\bold L[Z],\mu_\alpha$ is not necessary a cardinal,
in this case $f_\alpha = f \restriction A_\alpha$ can serve!
\mn
[Why?  In $\bold L[Z]$ we can choose $\langle f_\alpha:\alpha \in
{\Cal U}\rangle$ in $\circledast$ and then put together $f$ and
$\cup\{f_\alpha:\alpha \in {\Cal U}\}$ as above.]

The proof of the condition in $\oplus$ is by cases (on $\alpha$):
\bn
\ub{Case 1}:  $\alpha \in {\Cal U}$ and $\mu_\alpha$ is not a cardinal
in $\bold L[Z]$ or $\mu_\alpha < \kappa$.

Trivial.
\sn
Hence by clause (a) of the assumption
\mr
\item "{$(*)_2$}"  \wilog \, $\bold L[Z] \models ``\mu_\alpha$ is a
cardinality".
\ermn
\ub{Case 2}:  In $\bold L[Z],\mu_\alpha$ is regular $> \kappa$.

Let $B_{\alpha,\varepsilon} = \{\eta_{\alpha,\zeta}(n):n < \omega
\text{ and } \zeta < \varepsilon\}$, so in $\bold L[Z],\langle
B_{\alpha,\varepsilon}:\varepsilon < \mu_\alpha\rangle$ is
$\subseteq$-increasing continuous and let $C_{\alpha,0} = \{\delta <
\lambda:\delta$ is a limit ordinal and 
for every $\varepsilon < \mu$ we have $\varepsilon < \delta$
\ub{iff} for some $\zeta < \delta$, Rang$(\eta_{\alpha,\varepsilon})
\subseteq B_{\alpha,\zeta}\}$.

In $\bold L[Z]$ there is a club $C_\alpha = \{\beta_\xi:\xi <
\mu_\alpha\}$ of $\mu_\alpha$ such that $\delta \in C = \text{
cf}(\delta)^{\bold L[Y]} > \aleph_0$ and $C_\alpha \subseteq C$ and
$\beta_0 = 0$.

For $\varepsilon < \mu_\alpha$ let $\xi = \xi(\varepsilon)$ be maximal
such that $\varepsilon \ge \beta_\xi$ and easily
$\eta_{\alpha,\varepsilon} \notin
{}^{\omega}(B_{\alpha,\beta_\varepsilon})$, and let
$g(\eta_{\alpha,\varepsilon})$ be the shortest $\nu \trianglelefteq
\eta_{\alpha,\varepsilon}$ which $\notin B_{\alpha,\beta_{\xi(\varepsilon)}}$.

Now check. 
\bn
\ub{Case 3}:  cf$^{\bold L[Z]}(\mu_\alpha) \ge \kappa$.

Similarly.
\bn
\ub{Case 4}:  cf$^{\bold L[Z]}(\mu_\alpha)) = \aleph_0$.

Here we can find an increasing sequence $\langle B_n:n <
\omega\rangle$ of subsets of $\lambda$ of cardinality $< \mu$ such
that $A_\alpha \subseteq \dbcu_{n < \omega} {}^\omega(B_n)$.

So we can proceed as above.  \hfill$\square_{\scite{4.17}}$
\enddemo
\bn
\ub{Discussion}:  Question \scite{4.17B} seems to me to 
call for iterating Radin forcing but for
$\aleph_2$ there is a short cut.  For this we quote.
\bigskip

\proclaim{\stag{4.17.1} Theorem}  Assume ZF + DC + AD and $\kappa =
\aleph_1$.  \ub{Then}
\mr
\item "{$(*)_\kappa$}"   for every $A \subseteq \kappa$ for
some $\eta \in {}^\omega 2$ we have $A \in \bold L[\eta]$ 
and $\eta^\#$ (hence $A^\#$) exist.
\endroster
\endproclaim
\bigskip

\demo{Proof}  Well known.
\enddemo
\bigskip

\proclaim{\stag{4.17.2} Claim}  [ZF] 1) If DC + AD + $\kappa = \aleph_1$
or just $(*)_\kappa$ from \scite{4.17.1} holds, \ub{then} $\aleph_1$
is free.
\nl
2) Also $\kappa$ is Ord-free (see Definition \scite{4.23} below).
\endproclaim
\bigskip

\demo{Proof}  1) We can easily check the criterion from \scite{4.17}
as for $M$ a model with universe $\kappa$ and vocabulary $\subseteq
\bold L_\omega$, let $\eta \in {}^\omega 2$
be such that $M \in \bold L[\eta]$ and can work in $\bold
L[\eta,\eta^\#]$.
\nl
2) Easy, too.  \hfill$\square_{\scite{4.17.2}}$
\enddemo
\bigskip

\demo{\stag{7.6} Observation}  [ZFC + DC] If $(*)_{\lambda,\partial}$
then $\boxdot_{\lambda,\partial}$ where
\mr
\item "{$(*)_{\lambda,\theta}$}"   for every $A \subseteq \lambda$ there is $B
\subseteq \partial$ such that $A \in \bold L[B]$ and $B^\#$ exists (so
 $(*)_\kappa$ is $(*)_{\kappa,\aleph_0}$)
\sn
\item "{$\boxdot_{\lambda,\theta}$}"    every model $M$ of 
cardinality $\lambda$ with vocabulary of cardinality $\le \partial$ 
(so $\tau_M$ well ordered) is isomorphic to a model of the form
EM$_\tau(\lambda,\Phi)$ for some template $\Phi$ with 
$|\tau_\Phi| \le \theta$ (so $\tau_\Phi$ well ordered).  
\endroster
\enddemo
\bigskip

\remark{Remark}  This includes $(\lambda,<_\alpha)$ where $<_\alpha$
is a well order of $\lambda$ of order type 
$\alpha \in [\lambda,\lambda^+]$.  
\endremark
\bigskip

\proclaim{\stag{4.21} Claim}   Assume $T \subseteq T_1$ are countable
complete first order theories.
\nl
1)  If $T$ is stable not superstable and $\lambda > \aleph_0 + |T_1|$
 is not free (see Definition \scite{4.17A}) \ub{then} 
{\rm PC}$(T_1,T)$ is not categorical in $\lambda$.
\nl
2) If $T$ is unstable and $\lambda > \aleph_0$ \ub{then} {\rm
   PC}$(T_1,T)$ is not categorical in $\lambda$.
\endproclaim
\bigskip

\demo{Proof}  Without loss of generality $T,T_1 \subseteq \bold L_\omega$.
\nl
1) Working in $\bold L[T_1,T]$ we can find 
$\Phi$ proper for trees with $\omega +1$ levels as in \cite[VII]{Sh:c},
i.e., $\tau_\Phi \in \bold L[T_1,T]$, EM$_{\tau(T_1)}(I,\Phi)$ a model
of $T_1$ (e.g. for $I \subseteq {}^{\omega \ge} \lambda$) satisfying
EM$({}^{\omega \ge}\lambda,\Phi) \models \varphi_n(\bar a_\eta,\bar
a_\nu)^{\text{if}(\nu=\eta\restriction n)}$ when $\eta \in
{}^\omega \lambda,\nu \in {}^n \lambda$. 

Let $F:\lambda \rightarrow {}^\omega \lambda$ exemplify that $\lambda$
is not free, i.e., its range is not free.  
Working in $\bold L[T,T_1,F]$ (so \wilog \, $F$ is
one to one), let $M_1 = \text{ EM}_\tau({}^1 \lambda,\Phi),M_2 =
\text{ EM}_\tau({}^{\omega >} \lambda \cup \text{ Rang}(F),\Phi)$ and
assume toward contradiction that $f$ is an isomorphism from $M_1$ onto
$M_2$ and we shall work in $\bold L[T,T_1,F,f]$, in this universe let
${\Cal U} \subseteq \lambda$ be of minimal cardinality such that
$\{F(\alpha):\alpha \in {\Cal U}\}$ is not free (in the same sense).  
By \cite{Sh:52} (or \cite{Sh:E18}), $|{\Cal U}|$ is a
regular uncountable cardinal, so by renaming \wilog \, ${\Cal U} = \mu =
\text{ cf}(\mu) > \aleph_0$.  Let $W \subseteq \lambda,|W| =
\mu,\{f(a_\alpha):\alpha \in {\Cal U}\} \subseteq \text{ EM}({}^1 W,
\Phi)$ and let $\langle w_\alpha:\alpha < \mu \rangle$ be a  
filtration of $W$.  Clearly $M_1$ satisfies $A \subseteq
M_1 \wedge |A| < \mu \Rightarrow {\bold S}(A,M) = \{\text{tp}(a,A,M):a
\in M_1\}$ has cardinality $\le |A| + \aleph_0 < \mu$.  This holds in
$M_2$ hence $(\forall \alpha < \mu)(\exists \beta < \mu)
[\forall \gamma \in w_\alpha)[\{f(a_{F(\gamma) \restriction
n}):n \le \omega\}\subseteq \text{ EM}({}^1(w_\beta),\Phi)$.  We
continue as in \cite[VIII,\S2]{Sh:c} and get contradiction.
\nl
2) As in \scite{st.5}.  \hfill$\square_{\scite{4.21}}$
\enddemo
\bigskip

\proclaim{\stag{4.22} Claim}  Assume $\lambda > \aleph_0$ is a free cardinal.
\nl
1) For $T = { \text{\rm Th\/}}({}^\omega \omega,E_n)_{n < \omega},E_n
 = \{(\eta,\nu):\eta,\nu \in {}^\omega \omega,\eta \restriction n =
\nu \restriction n)$ for some countable complete $T_1 \supseteq T$,
 {\rm PC}$(T_1,T)$ is categorical in $\lambda$ ($T_1$ does not depend
 on $\lambda$).
\nl
2) There is a countable complete stable not superstable $T$ such that
   if $M \models T,\|M\| \le \lambda$, \ub{then} the isomorphism type
 of $M$ is determined by two dimensions. 
\endproclaim
\bigskip

\demo{Proof}  1) As in \cite{Sh:100}, $T_1$ will guarantee that for
any $M \in \text{ PC}(T_1,T)$ we have:
\mr
\item "{$(*)_1$}"   if $a \in M$ then $\{b \in M:M \models b E_n a$
for every $n < \omega\}$ has cardinality $\|M\|$
\sn
\item "{$(*)_2$}"  if $a \in M,n < \omega$ then $\{b/E^M_{n+1}:b \in
a/E_n\}$ has cardinality $\|M\|$.
\ermn
So suppose $M_1,M_2 \in \text{ PC}(T_1,T)$ has universe $\lambda$ and
we work
in $\bold L[T,M_1,M_2]$.  There is $M'_\ell \cong M_\ell$ of
cardinality $\lambda$ and $A_\ell \subseteq {}^\omega \lambda,|A_\ell|
= \lambda$ for $\ell=1,2$ such that
\mr
\item "{$(*)_3$}"  $|M'_\ell| = A_\ell \times \lambda,(\eta,\alpha)
E_n(\nu,\beta)$ iff $(\eta,\nu \in A_\alpha,\alpha,\beta <
\lambda$ and) $\eta \restriction n = \nu \restriction n$
\sn
\item "{$(*)_4$}"  $\nu \in {}^{\omega >}\lambda \Rightarrow
(\exists^\lambda \eta)(\nu \triangleleft \eta \in A_\ell)$.
\ermn
By the assummption ``$\lambda$ is free" (see Definition \scite{4.17A}) 
we can find $g_\ell:A_\ell \rightarrow \omega$ such
that $\langle \eta \restriction g_\ell(\eta):\eta \in A_\ell \rangle$
is with no repetitions and we shall work in 
$\bold L[T_2,M_1,M_2,A_1,A_2,g_1,g_2]$.  
For $\kappa < \mu$ let ${\Cal F}_\kappa$ be
the family of functions $h$ such that
\mr
\item "{$(*)^5_h$}"  $(a) \quad h$ is a partial one-to-one function
from $A_1$ into $A_2$
\sn
\item "{${{}}$}"  $(b) \quad |\text{Dom}(h)| = \kappa$
\sn
\item "{${{}}$}"  $(c) \quad$ for $\eta_1,\eta_2 \in \text{ Dom}(h)$ and
$n < \omega$ we have $\eta_1 \restriction n = \eta_2 \restriction n
\Leftrightarrow h(\eta_1) \restriction$
\nl

\hskip25pt  $n = h(\eta_2)$
\sn
\item "{${{}}$}"  $(d) \quad$ if $\ell \in \{1,2\}$ and $\nu \in
A_\ell$ and $(\forall n <\omega)(\exists \eta \in A_\ell)(\nu
\restriction n = \eta \restriction n)$ then \nl

\hskip25pt $\nu \in A_\ell$.
\ermn
Let $A_\ell = \{\eta^\ell_\alpha:\alpha < \lambda\}$.  It is easy to
choose $h_\alpha \in {\Cal F}_{\aleph_0+|\alpha|}$ by induction on
$\alpha$ increasing continuous with $\alpha$ such that $\eta^1_\alpha
\in \text{ Dom}(h_{\alpha +1}),\eta^2_\alpha \in \text{ Rang}(h_{\alpha
+1})$.
\nl
2) As in Example \scite{3c.38} using $\omega$-power. 
\enddemo
\bigskip

\definition{\stag{4.23} Definition}  1) We say $\lambda$ is Ord-$\mu$-free
\ub{when}:
\block
for every linear order $M = (\lambda,<^M),\lambda$ for some $B
\subseteq \lambda$ in $\bold L[A,B],M$ can be represented as
$\cup \{M_i:i < \mu\},M_i$ embeddable into $({}^n
\lambda,<_{\text{even}})$ where $\eta <_{\text{even}} \nu
\Leftrightarrow (\exists m < n)(m = \ell g(\eta) = \ell g(\nu) \wedge
(\eta(m) \ne \nu(m)) \wedge (\eta(m) < \nu(m) \equiv m$ even)
(see Laver \cite{Lv71}, \cite[XII,\S2]{Sh:e}).
\endblock
\sn
2) If $\mu = \aleph_0$ we may omit it.
\enddefinition
\bigskip

\proclaim{\stag{4.24} Claim}  If $\lambda$ is {\rm Ord}-free \ub{then} any
two strongly $\aleph_0$-homogeneous linear orders (see below) of cardinality
$\lambda$ of the same cofinality are isomorphic.
\endproclaim
\bigskip

\definition{\stag{4.24.1} Definition}  $I$ is a strongly
$\aleph_0$-homogeneous if $I$ is infinite dense isomorphic to any open
interval and its interval.
\enddefinition

\demo{Proof}  See above.
\enddemo
\newpage

\head {\S6 Comments on model theory in ZF} \endhead  \resetall \sectno=6
 \spuriousreset
\bigskip

Before we comment on model theory without choice we write up the
amount of absolute which holds.
\demo{\stag{6x.0} Observation}   Let $T$ be countable complete first order
theory, without loss of generality 
$\Bbb L_{\tau(T)} \subseteq {\Cal H}(\aleph_0)$
(or if you like $\subseteq \omega$), so $T \subseteq {\Cal H}(\aleph_0)$.
\nl
1) ``$T$ is stable" is a Borel relation.
\nl
2) ``$M$ is a countable model of $T,q(\bar y) \in \bold S^{<
\omega}(M),p(\bar x) \in \bold S^{< \omega}(M)$ and $M_\ell \prec
M$ for $\ell = 0,1,2$ and for stable 
$T,\nonforkin{M_1}{M_2}_{M_0}^{M},p$ does not fork over $M_0$, all
coded naturally as a subset of $\omega$" are Borel.
\nl
3) In part (2), ``$p \perp q"$ is Borel as well as 
$``p \underset{\text{wk}} {}\to \perp q"$ is Borel, also $``p \perp M_0"$
by clause (e) of part (3A).
\nl
3A) Let $T^{\text{eq}}$ be $T$ when we add predicates naming  the
equivalence classes so have a predicate $P_{\varphi(\bar x,\bar y)}$
equivalent to every $\varphi(\bar x) \in \Bbb L(\tau_T)$,
(\cite[III]{Sh:c}) and $T^{\text{eq}}_\forall$ be the universal part
(pedantically the consequences of $T^{\text{eq}}$), 
so $\tau(T),\Bbb L(\tau(T^{\text{eq}})),
T^{\text{eq}},T^{\text{eq}}_\forall$ are Borel definable from $T$.
Also the following are  Borel
\mr
\item "{$(a)$}"  $A$ is a model of $T^{\text{eq}}_\forall$ in this
observation with universe $\subseteq \omega$ and we 
use $A,A_\ell$ to denote such models
\sn
\item "{$(b)$}" $A_1 \subseteq A_2$ are models of $T^{\text{eq}}_\forall,A_2 =
acl(A_1)$ (in any $M,A_2 \subseteq M \models T^{\text{eq}}$), and
computing such $A_2$ naturally defined
\sn
\item "{$(c)$}"  $p \in \bold S^m(A),A$ a model of $T^{\text{eq}}_\forall$;
i.e. $\{\text{tp}(\bar a,A,M):A \subseteq M \models
T^{\text{eq}},\bar a \in {}^r M\}$ we may write acl$(A)$
\sn
\item "{$(d)$}"   computing $p \restriction A_2$ from $A_1
\subseteq A_2$ and $p \in \bold S^m(A_2)$
\sn
\item "{$(e)$}"  computing $R^m(p,\Delta,2),R^m(p,\Delta,\aleph_0)$,
rk$^m(p,\Delta,\aleph_0)$ for $p$ an $m$-type over $A$, (a model of
$T^{\text{eq}}_\forall$), $\Delta \subseteq \Bbb L(\tau_T)$ finite)
\sn
\item "{$(f)$}"  $A_1 \subseteq A_2,p(\bar x)$ an $m$-type over $A_2$
(in clause (c)'s sense) and $p(\bar x)$ does not fork over $A_2$
\sn
\item "{$(g)$}"  $A_1 \subseteq A_2,p(\bar x)$ an $m$-type over $A_2$,
the type $p(\bar x)$ does not fork over $A_2$ and is stationary over $A_1$
\sn
\item "{$(h)$}"  in $(g)$ computing the unique extension $q \in \bold
S^{\ell g(\bar x)}(A_2)$ of $p(\bar x)$ not forking over $A_1$ and
tp$(\bar a_0 \char 94 \bar a_1 \char 94 \bar a_1 \char 94 \dotsc
A_2,M)$ when $A_2 \subseteq M \models T^{\text{eq}},\bar a_n$
realizes $p(\bar x)$ in $M$ and $p^\omega_{A_2} = \text{ tp}((\bar
a_n,A_2 \cup \{\bar a_0,\dotsc,\bar a_{n-1}\},M)$ does not fork over
$A_1$
\sn
\item "{$(i)$}"  $p_\ell(\bar x_\ell) \in \bold S^{m(\ell)}(A)$ for
$\ell=1,2$ are weakly orthogonal
\sn
\item "{$(j)$}"  for $A_\ell \subseteq A$, the stationary types $p_\ell(\bar
x) \in \bold S^{m(\ell)}(A_\ell)$ for $\ell=1,2$ are orthogonal
\sn
\item "{$(k)$}"  from $A \subseteq A_\ell,A_2$ such that
tp$(A_\ell,A)$ is stationary for $\ell=1,2$ and computing
$A',(f_\ell,A'_\ell)$ for $\ell=1,2$ such that $A \subseteq A',f_\ell$
an isomorphism from $A_\ell$ onto $A'_\ell$ over $A,A'_\ell \subseteq
A$ for $\ell=1,2$ and $\nonforkin{A'_1}{A'_2}_{A}^{A'}$
\sn
\item "{$(l)$}"  $A_1 \subseteq A_2,A \subseteq A_2$ and $p(\bar x) \in
\bold S^m(A)$ is orthogonal to $A_1$.
\ermn
4) ``$T$ has DOP" is a $\Sigma^1_1$-relation (so NDOP is $\Pi^1_2$).
\nl
5) ``$T$ has DIDIP" is $\Sigma^1_1$ (so NDIDIP is $\Pi^1_1$).
\nl
6) ``$T$ has OTOP" in $\Sigma^1_1$.
\enddemo
\bigskip

\demo{Proof}  Sometimes we give equivalent formulations to prove.

(1),(2) are obvious; for ``dnf" see clause (f) of part (3A).
\nl
3) 
\mr
\item "{$(a)$}"    $p \underset{\text{wq}} {}\to \perp q$ just says:
for some $A \subseteq M,p,q \in \bold S^{< \omega}(A,M)$ (or even $p,q
\in \bold S^{\le \omega}(A,M))$ satisfying: if
 $\varphi(\bar x,\bar y,\bar z) \in \Bbb L(\tau_T)$ and $\bar a$ from
$A$ we have 
$p(\bar x) \cup q(\bar y) \vdash \varphi(\bar x,\bar y,\bar a)$
or $p(\bar x) \cup p(\bar y) \vdash \neg \varphi(\bar x,\bar y,\bar a)$
and remember compactness
\sn
\item "{$(b)$}"  $p \perp q$, see clause (j) of
part (3A)
\sn
\item "{$(c)$}"  $p \perp M_0$ by clause $(\ell)(\beta)$ of part (3A).
\endroster
\sn
3A) E.g.
\mn
\ub{Clause (e)}:  Because $\Delta$ is finite, the value is a natural
number and for $\theta \le \aleph_0,k < \omega$ we have $R^m(p(\bar
x),\Delta,\theta) \ge k$ \ub{iff} some Borel set of formulas, see
\cite[II,\S2]{Sh:c} is consistent.  Simiarly for $(R^m(p(\bar x),\Delta,\theta)
> k_1) \vee (R^m(p(\bar x),\Delta,\theta) = k_1) \wedge \text{
Mlt}^m(p(\bar x),\Delta,\theta) \ge k_2)$.
\mn
\ub{Clause (f)}:  This is equivalent to ``if $A_3 = acl(A_2),\Delta
\subseteq \Bbb L(\tau_{T^{\text{eq}}})$ is finite then there is $q \in
\bold S^m_\Delta(A_2)$, i.e. a definition of such type which extend $p
\restriction \Delta$ and is definable over acl$(A_1)$.
\mn
\ub{Clause (g)}:  We can use the definition: for every finite $\Delta$
there is $q \in \bold S^m_\Delta(A_3)$ definable over $acl(A_1,A_3)$
such that $R^m(p(\bar x),\Delta,2) = R^m(p(\bar x) \cup q(\bar x),\Delta,2)$.
\mn
\ub{Clause (j)}:  This is equivalent to: $p_\ell \in \bold
S^{<\omega}(A_\ell)$ for $\ell=1,2,A_\ell \subseteq M$, and for every
$n < \omega$ and finite $\Delta_1 \subseteq \Bbb L(\tau_T)$ for some
finite $\Delta_2 \subseteq \Bbb L(\tau_T)$, if $\langle
a^\ell_0,\dotsc,a^\ell_{n-1}\rangle$ is as in clause (h) with
$(A_\ell,A_1 \cup A_2,p_\ell)$ here standing for $(A_1,A_2,p)$ 
\ub{but} we have finitely many possibilities for each, \ub{then}
tp$_{\Delta_2}(\bar a^1_0 \char 94 \ldots \char 94 a^1_{n-1},A_1 \cup
A_2,M)$,tp$_{\Delta_1}(\bar a^2_0 \char 94 \ldots \char 94 \bar
a^2_{n-1},A_1 \cup A_2,M)$ determine the $A_1$-type of $\langle \bar
a^1_0 \char 94 \ldots \char 94 \bar a^1_{n-1} \char 94 \bar a^2_0
\char 94 \ldots \char 94 \bar a^2_{n-1}\rangle$ over $A_1 \cup A_2$ in $M$.
\mn
\ub{Clause $(\ell)$}:  First assume $A_1,A_2,A$ are
algebraically closed.  We know that $p \perp A_1$ iff there are $f,M$
such that $A_2 \subseteq M \models T^{\text{eq}},f \supseteq \text{
id}_A(M,M)$-elementary mapping (i.e. an automorphism of $M$) and
mapping $A_2$ to $A'_2,\nonforkin{A_1}{A'_2}_{A}^{M}$ such that $p
\perp f(p)$.  In the general case as in clause (j) work with ``for
every finite $\Delta_1$ ...".
\nl
4) Obvious by (3) and by the definition (there are countable models of
 $T,M_\ell(\ell \le 3)$ such that
$\nonforkin{M_1}{M_2}_{M_1}^{M_3},M_3$ is $\bold
F^\ell_{\aleph_0}$-constructible over $M_1 \cup M_2$ and 
$p \in \bold S^{< \omega}(M_3)$ non-algebraic such that $p
 \perp M_1,p \perp M_2$).
\nl
5) Obvious by (3) and the definition (equivalent to: there are 
countable models $M_n$ of $T,M_n \prec M_{n+1}$, and countable $N$
which is $\bold F^\ell_{\aleph_0}$-atomic over $\cup\{M_n:n <
\omega\}$ and non-algebraic $p \in \bold S^{< \omega}(N)$ such that $n <
\omega \Rightarrow p_n \perp M_n$).
\enddemo
\bigskip

\proclaim{\stag{6x.1} Claim}: 0) \ub{Convention}:
\mr
\item "{$(a)$}"  $T$'s vocabulary, $\tau = \tau_T$ is well orderable
and for simplicity $\subseteq \bold L$
\sn
\item "{$(b)$}"  $M,N$ denote models of $T$ with universe a set of
ordinals
\sn
\item "{$(c)$}"  $T$ a theory in $\Bbb L(\tau)$ so $|T|$ is a cardinal;
\wilog \, $\tau \subseteq \bold L_\lambda,\lambda = |T| + \aleph_0$
\sn
\item "{$(d)$}"  ``a model of $T$" means one with well ordered
universe so \wilog \, a set of ordinals.
\ermn
1) DLST and ULST holds (for models as in clause (b)), short for the
downward L\"owenheim-Skolem-Tarski and the upward
L\"owenheim-Skolem-Tarski theorems respectively.  If 
$T$ is categorical in $\lambda \ge |T|$ then $T \cup
\{\exists^{\ge n} x(x=x):n < \omega\}$ is complete, etc., all
that takes place in some $\bold L[Y]$ is fine. 
\nl
2) [$T$ complete] $T$ has an $\aleph_0$-saturated model iff every
 model of $T$ has an $\aleph_0$-saturated elementary extension 
iff $D(T)$ can be well ordered.
\nl
3) Define $\kappa(T) =: \sup\{\kappa(T)^{\bold L[T,Y]}:Y$ a set of
ordinals$\}$ but probably better to use
\nl
$\kappa^+(T) = \cup\{(\kappa(T)^+)^{\bold L(T,Y)}:Y$ a set of ordinals$\}$.
\nl
4) Assume $T$ is complete.  Every model $M$ of $T$ of cardinality 
$\le \lambda$ has a $\kappa$-saturated elementary extension of 
cardinality $\le \lambda$ iff $|D(T)| \le \lambda$ and $(a) \vee (b)$ where
\mr
\item "{$(a)$}"  $|{}^{\kappa >} \lambda| = \lambda$,
i.e. $[\lambda]^{< \kappa}$ is well ordered
\sn
\item "{$(b)$}"  $T = \text{\rm Th}(M)$ is stable, 
$|\lambda^{< \kappa(T)}| = \lambda$ 
 and $|{\Cal P}(\omega)|$ is a cardinal $\le \lambda$ if some 
$p \in \bold S(B),B \subseteq M \models T,M$ well orderable, $|B| < 
\kappa(T)$ has a perfect set of stationarization and $\lambda > \aleph_0$.
\ermn
\eightpoint
[Why?  As in \cite[III]{Sh:c}, particularly section 5, hopefully
see the proof of \cite[1.1]{Sh:F701}.]
\tenpoint
\nl
5) [$T$ complete]
\mr
\item "{$(a)$}"  if $\neg(|D(T)| \le |T|)$ \ub{then} there is a family
${\Cal P}$ of subsets of $D(T)$, each of cardinality $\le |T|$ and $\cup\{\bold
P:\bold P \in {\Cal P}\} = D(T)$ and for each $\bold P \in {\Cal P}$
there is a $\Phi$ proper for linear orders, with $\tau(T),\tau(\Phi)
\subseteq \bold L$
such that every model EM$_{\tau(T)}(I,\Phi)$ satisfies: the model 
realizes $p \in D(T)$ iff $p \in \bold P$ 
\ermn
6) Assume there is no set of $\aleph_1$ reals.  If $T$ 
is complete countable, $D(T)$ uncountable, $M \models T$ and
$\bold P_M = \{p \in D(T):M$ realized $p\}$ \ub{then} $\bold P_M$
 is countable.
\nl
6A) Of course, it is possible that $|D(T)|^{\bold L[T,Y]}$ is large in
$\bold L[T,Y]$! (e.g. there is a set of $|T|$ independent formulas see 12)(c)).
\nl
7) If $T$ is complete not superstable, $T_1 \supseteq T$ complete,
$\lambda = \text{\rm cf}(\lambda) > |T_1|,
\lambda \ge \theta({\Cal P}(\omega))$ 
and axiom Ax$^3_\lambda$ (see \cite{Sh:835}, i.e.,
$|[\lambda]^{\aleph_0}|$ a cardinal) \ub{then} there is $\langle M_u:u
\subseteq \lambda \rangle$ such that
\mr
\item "{$(a)$}"  $M_u \in \text{\rm PC}(T_1,T)$
\sn
\item "{$(b)$}"  $\|M_u\| = \lambda$
\sn
\item "{$(c)$}"  $u \ne v \subseteq \lambda \Rightarrow M_u \ncong M_v$.
\nl
[Why? There is a sequence $\langle C_\delta:\delta \in S
\rangle,S \subseteq S^\lambda_{\aleph_0}$ stationary $C_\delta
\subseteq \delta = \sup(C_\delta)$, otp$(C_\delta) = 0$ (hence we can
partition $S$ to $\lambda$ stationary sets)].
\ermn
8) Define $\beth'_\alpha(\lambda)$ by $\beth'_0(\lambda) =
\lambda,\beth'_{\alpha +1}(\lambda) = \theta({\Cal
P}(\beth'_\alpha(\lambda)),\beth'_\delta(\lambda) =
\cup\{\beth'_\alpha(\lambda):\alpha < \delta\}$, it is a cardinal.
If $T$ is countable, $\Gamma$ is a countable set of 
$\bold L(\tau_T)$-types and for every $\alpha < \omega_1$ there is $M
\in \text{\rm EC}_{\beth_\alpha}(T,\Gamma)$ so $|M| \subseteq
\bold L$ of power $\beth'_\alpha$ or just of power 
$\ge \beth^{\bold L[M]}_\alpha$ (but we do not say that an
$\omega_1$-sequence of such models exists!), \ub{then} there is
an $\Phi \in \Upsilon^{\text{or}}_{\aleph_0}[T]$ so 
$|\tau_\Phi| = \aleph_0$ such that {\rm EM}$_{\tau(T)}(I,\Phi) \in 
\text{\rm EC}(T,\Gamma)$ for every linear order $I$.
\nl
\eightpoint
[Why?  See proof of (9), but here the members of the tree are finite
set of formulas hence the tree is $\subseteq \bold L[T]$ and we can
define the rank in $\bold L[T]$ \ub{but}: we let $\langle \Delta_n:n
< \omega\rangle$ be an increasing sequence of finite sets of formulas, 
each $\varphi \in \Delta_n$ has a set of free variables $\subseteq
\{x_0,\dotsc,x_{n-1}\}$ and $\dbcu_n \Delta_n = \Bbb
L(\tau_T),\Delta_n$ is closed under change of free variables (modulo
the restriction above).  We define ${\Cal T}_n$ as in the proof of
part (9) by $p \in {\Cal T}_n$ is a complete $(\Delta_n,n)$-type.  The
tree is really $\subseteq {}^{\omega >}\omega$.]
\tenpoint 
\nl
9) [DC] Assume $T$ is an (infinite) theory with Skolem functions,
$\Gamma$ a set of $\Bbb L(\tau_T)$-types and for every $\alpha <
\theta({\Cal P}(|T|))$ there is $M \in \text{\rm EC}(T,\Gamma)$ of
power $\ge \beth_\alpha$ in $\bold L[T,M]$, \ub{then} there is $\Phi$ such 
that EM$_{\tau(T)}(I,\Phi)
\in \text{ EC}(T,\Gamma)$ for every linear order $I$.
\nl
\eightpoint
[Why?  A wrong way is to assume
$\theta({\Cal P}(|T|))$ is regular and in stage $n$ we have an
$n$-indisernible sequence $\bold I^n_\alpha \subseteq M$ of
cardinality $\beth'_\alpha$ for $\alpha < \theta({\Cal P}(|T|))$ with
$n$-tuple from $\bold I^n_\alpha$ realizing $p_n$, as in the ZFC
proof.  The problem is
that there may be no regular cardinal $\ge \theta({\Cal
P}(|T|))$.  But more carefully let ${\Cal T}_n$ be the set of complete types
$p_n(x_0,\dotsc,x_{n-1})$ consistent with $T$, such that it is the 
type of a sequence
of length $n$ which is $m$-indiscernible for each $m \le n$.  The
order on ${\Cal T} = \cup\{{\Cal T}_n:n < \omega\}$ is inclusion, so
really ${\Cal T}_n$ is the $n$-th level.
We need DC$_{\aleph_0}$ to have
a rank function on this set which has power $\le {\Cal P}(|T|)$.  
We prove by induction on the ordinal $\gamma$ for each $n$,
that if $p \in {\Cal T}_n$ has rank $\gamma$ that no indiscernible 
$\bold I \subseteq M,M \in 
\text{ EC}(T,\Gamma)$ of cardinality $\ge \beth^{\bold L[T]}_{\omega
\gamma}$ exists.]
\tenpoint
\nl
9A) Of course, if EC$(T,\Gamma)$ has a model $M,|M| \subseteq \bold L$
of cardinality $\ge \beth'_\delta$ where $\delta := 
\theta({\Cal P}(|T|))$ then we do not need DC.
\nl
9B)  We can avoid ``$T$ has Skolem functions", see \cite{Sh:F701}, in
both parts (9) and (9A).  The point is that $T$ needs not be complete,
\wilog \, $T$ has elimination of quantifiers and we can define
$T^{\text{SK}}$ which is $T +$ the axioms of Skolem functions; now for
every $\alpha < \theta({\Cal P}(|T|)$, there is a model $M$ of $T$ of
cardinality $\ge \beth_\alpha$ it can be expanded to $M^+$, a model of
$T^{\text{SK}}$ and we can continue (with new function symbols).
\nl
10) Assume $T$ is complete uncountable.  Then all the proofs in \S2 +
\S3 holds except that we do not have the dichotomy OTOP/existence of
primes over stable amalgamation.  We intend to return to it in
\cite{Sh:F701}.
\nl
11) If $p(x_0,\dotsc,x_{n-1})$ is a set of $\Bbb L(\tau_T)$-formulas
    consistent with $T$ \ub{then} it is realized in some model $M$ in
    some universe $\bold L[T,Y]$ hence can be extended to a complete
    type realized in such $M,p$ hence $\in D_n(T)$ when $T$ is
    complete.
\nl
\eightpoint
[Why?  Work in $\bold L[T,p]$ O.K. as $p \subseteq \bold L$ as $T
\subseteq L$.]
\tenpoint
\endproclaim
\bigskip

\proclaim{\stag{4.10} Lemma}  \cite{Sh:c} can be done in {\rm ZF} $+ 
(\forall \alpha)([\alpha]^{\aleph_0}$ is well ordered), see
\cite{Sh:835} as long as
\mr
\item "{$(a)$}"  the theory $T$ is in a vocabulary which can be well
ordered
\sn 
\item "{$(b)$}"   we deal only with models whose power is a cardinal
\sn
\item "{$(c)$}"  all notions are in $\bold L[T,Y],Y \subseteq {
\text{\rm Ord\/}}$ large enough (so ${\frak C}$ is not constant it
depends on the universe)
\sn
\item "{$(d)$}"  in \cite[VIII]{Sh:c}, the case $\lambda > |T_1|$ regular
($\tau(T_1)$ well orderable, too) is clear
as using the well ordering
$[\lambda]^{\aleph_0}$ we can find $\langle C_\delta:\delta \in
S^\lambda_{\aleph_0} \rangle \in \bold L[T,Y]$ hence define a
partition $\langle S_\alpha:\alpha < \lambda \rangle$ of
$S^\lambda_{\aleph_0}$ such that $(\exists^\lambda \alpha)$
($S_\alpha$ stationary (in $\bold V$), so increasing $Y$ we are there
but
\sn
\item "{$(e)$}"  Ch VI on ultrapower should be considered separately.
\endroster
\endproclaim
\newpage

\head {\S7 Powers which are not cardinals} \endhead  \resetall \sectno=7
 \spuriousreset
\bigskip

We suggest to look at categoricity of countable theories in so-called
reasonable cardinals.  For them we have the completeness theorem in
\scite{8.8}.  We then uncharacteristically examine a classical
example: Ehrenfuecht example (with 3 models in $\aleph_0$, see
\scite{8.14}).  
\bn
We naturally ask
\nl
\ub{Question}:  Can an expansion of the theory of linear orders be
categorical in some uncountable power?

We then deal with criterion, i.e. sufficient conditions for categoricity.
We intend to continue this in \cite{Sh:F701}.
\bigskip

\demo{\stag{8.0.7} Convention}  $T$ not necessarily $\subseteq \bold L$.
\enddemo
\bn
We may consider
\definition{\stag{8.3} Definition}  1) For a class $\bold C$ of powers
we say $T_1 \le^{\text{ex}}_{\bold C} T_2$ when:
for every set $X$ of power $\in \bold C$ if $T_2$ has a model with
universe $X$ \ub{then} $T_1$ has a model with universe $X$.
\nl
2) For a class $\bold C$ of powers we say $T_1 \le^{\text{cat}}_{\bold
C} T_2$ when: for every set $X$ of power $\in \bold C$ if $T_2$ is
categorical in $|X|$, (i.e., has one and only one model with universe
$X$ up to isomorphism) then $T_1$ is categorical in $|X|$.
\nl
3) In both cases, if 
$\bold C$ is the class of all powers $\ge |T_2|$ we may omit it.
\enddefinition
\bigskip

\demo{\stag{8.4.1} Observation}  $\le^{\text{ex}}_{\bold C},
\le^{\text{cat}}_{\bold C}$ are partial orders.
\enddemo
\bn
We may also consider
\sn
\margintag{8.5}\ub{\stag{8.5} Question}:  1) For which countable theories 
$T$ is there a forcing extension $\bold V^{\Bbb P}$ of $\bold V$, 
model of ZF such that in $\bold V^{\Bbb P}$ the theory $T$ is 
categorical in some uncountable power?
\nl
2) As in (1) for reasonable powers, see below.
\bigskip

\definition{\stag{8.6} Definition}  We say that $X$ is a set of
reasonable power (or $|X|$ is a reasonable power) when:
\mr
\item "{$(a)$}"  there is a linear order of $X$
\sn
\item "{$(b)$}"   $|X| = |X \times X|$.
\endroster
\enddefinition
\bigskip

\proclaim{\stag{8.7} Claim}  If $T$ is countable theory and $X$ a set
of reasonable power \ub{then} $T$ has a model with universe $X$.
\endproclaim
\bigskip

\demo{Proof}  By \scite{8.8}.  \hfill$\square_{\scite{8.7}}$
\enddemo
\bigskip
 
\proclaim{\stag{8.8} Claim}  [ZF]  1) For some first order sentence
$\psi$ we have: for a set $X$ the following are equivalent:
\mr
\item "{$(a)$}"  $X$ is a set of reasonable power
\sn
\item "{$(b)$}"  if $T$ is a countable theory 
then $T$ has a model with universe $X$
\sn
\item "{$(c)$}"  $\psi$ has a model with universe $\psi$.
\ermn
2) If $T$ is categorical in $|X|$, a reasonable power 
then $T \cup \{(\exists^{\ge n}x)(x=x):n < \omega\}$ is a complete theory.
\endproclaim
\bigskip

\demo{Proof}  1) \ub{$(b) \Rightarrow (a)$}.

First apply clause (b) to $T_1 =$ (the theory of dense linear order with
neither first nor last elements), or just $T'_1 = \{\psi_1\} \subseteq
T_1$, where $\psi_1 \vdash T_1$ so it has a model $M=(X,<^M)$, so
$<^M$ linearly ordered $X$.

Second, apply clause (b) to $T_2 = \text{ Th}(\omega,F),F$ a one-to-one
function from $\omega \times \omega$ onto $\omega$, or just $T'_2 =
\{\psi_2\} \in T_2$ expresses this so there is a
model $M = (X,F^M)$ of $T$, so $F^M$ exemplifies $|X| = |X \times X|$.

Note that we have used (b) only for theories consisting of one
sentence.
\mn
\ub{$(a) \rightarrow (b)$}.

Use Ehrenfeuch-Mostoswki models.

That is it is enough to prove: using $I = (X,<)$ a linear order
\mr
\item "{$\boxplus$}"  if $T'$ is a countable complete theory with Skolem
functions, every term $\sigma(x_0,\dotsc,x_{n-1})$ is (by $T'$) equal
to a function symbol, $M' \models T$ and $\langle a_n:n <
\omega\rangle$ is an indiscernible sequence in $M',p_n = \text{
tp}_{\text{qf}}(\langle a_0,\dotsc,a_{n+1}\rangle,\emptyset,M)$ for $n
< \omega$ \ub{then} we can find $M,\langle a_t:t \in I\rangle$ such
that
{\roster
\itemitem{ $\circledast$ }  $(a) \quad M$ is a model of $T'$
\sn
\itemitem{ ${{}}$ }  $(b) \quad M^*$ has universe $X$
\sn
\itemitem{ ${{}}$ }  $(c) \quad \langle a_t:t \in I\rangle$ is an
indiscernible sequence in $M$
\sn
\itemitem{ ${{}}$ }  $(d) \quad \langle
a_{t_0},\dotsc,a_{t_{n-1}}\rangle$ realizes $p_n$ in $M$ when $t_0 <_I
\ldots <_I t_{n-1}$.
\endroster}
\ermn
Let $<^*$ be a well order $\tau(T)$.

Let $\langle(k_n,F_n):n < \alpha \le \omega\rangle$ list with no
repetition the pairs $(k,F)$ satisfying $(*)_{k,F}$ 
such that $k_0 = 1,M \models \forall x[F_0(x)=x]$ where
\mr
\item "{$(*)_{k,F}$}"  $(a) \quad F \in \tau(T)$ is a $k$-place function
symbol
\sn
\item "{${{}}$}"  $(b) \quad$ there is no $u \subset \{0,\dotsc,k-1\}$
such that $F^{M'}(a_0,a_1,\dotsc,a_{k-1}) \in$
\nl

\hskip25pt $\text{ Sk}_M(\{a_\ell:\ell \in u\})$
\sn
\item "{${{}}$}"  $(c) \quad$ there is no $k$-place function symbol
$F_1 \in \tau(T')$ such that
\nl

\hskip25pt $F_1 <^* F$ and 
$F^{M'}_1(a_0,\dotsc,a_{k_1}) = F^{M'}(a_0,\dotsc,a_{k-1})$. 
\ermn
Let $Y = \dbcu_{n <\omega} Y_n$ where $Y_n = \{(n,t_0,\dotsc,t_{k_n-1}):n <
\alpha$ and $t_0 <_I \ldots < t_{k_{n-1}}\}$.

Let $g:X \times X \rightarrow X$ be one to one onto.

Clearly there is a model as required with universe $Y$, hence it is
enough to prove $|Y| = |X|$.  Clearly $|X| \le |Y|$ as $\{(0,t):t \in I\}
\subseteq Y$.  Also $|Y_n| = |X|^{k_n}$ which is 1 if $k_n=0$ and is
$|X|$ if $k_n \ge 1$ as we can prove by induction on $n$.  Moreover,
we can choose $\langle f_n:n < \alpha,k_n \ge 1\rangle$ such that
$f_n$ is one-to-one from $Y_n$ onto $X$ as $f_n$ is gotten by
composition $k_n-1$ times of $g$.  This leads to $|Y| \le |X \times
\omega| + |\omega|$.  But trivially $\aleph_0 \le |X|$ by $g$ hence $|X|
\le |Y| \le |X| \times |X| + \aleph_0 = |X|$ hence we are 
done proving $(b) \Rightarrow (a)$.

Let $\psi$ say ``$<$ is a linear order and $F(x,y)$ is a one-to-one
function onto.

Now
\sn
\ub{$(c) \Rightarrow (a)$}:  as in the proof of $(b) \Rightarrow (a)$
\nl 
and also
\sn 
\ub{$(b) \Rightarrow (c)$}:  should be clear.
\nl
2) Easy, too.  \hfill$\square_{\scite{8.8}}$
\enddemo
\bn
\margintag{7f.16}\ub{\stag{7f.16} Discusion}   We can use 
an $\aleph_0$-saturated model $M$ of $T$ as a set of
urelements, i.e. we use a Fraenkel-Mostowski model for 
the triple ($M'$, a copy of $M$; finite support; finite partial 
automorphism of $M$).  Is $T$ categorical in
$|M'|$?  The problem is that maybe some $\psi \in \Bbb
L_{(2^{\aleph_0})^+,\omega}$ define in $M'$ with finitely many
parameters, a model $M''$ of $T$ with universe $|M'|$ such that
there is no permutation $f$ of $|M'|$ definable similarly
such that $f$ is an isomorphism from $M'$ onto $M''$.  But we may
consider $(D,\aleph_0)$-homomogeneous models of some extension of $T$
(in bigger vocabulary).  This seems related to \cite{Sh:199}, \cite{Sh:750}.
\bigskip

\definition{\stag{8.13} Definition}  $T_1$ is the theory of dense
linear order with neither first nor last element and $c_n < c_{n+1}$
for $n < \omega$ (so $\tau(T_1) = \{<\} \cup \{c_n:n < \omega\}$.
\enddefinition
\bigskip

\remark{Remark}  1) This is the Eherenfeucht example for $\dot I(\aleph_0,T)
= 3$.
\nl
2) We can replace $T_2$ by $T_{i,n}$ with 3 below replaced by $3+n$.
\endremark
\bigskip

\proclaim{\stag{8.14} Claim}  [ZF]  1) $T_1$ is a complete countable
first order which is not categorical in any infinite power.
\nl
2) In fact if $T_1$ has a model with universe $X$ \ub{then} $T_1$ has at
least three non-isomorphic models with this universe.
\nl
3) If in (2) the set $X$ is uncountable (i.e. $|X| \ne |\omega|$) 
\ub{then} $T_1$ has at least
$\aleph_0$ non-isomorphic models with this universe.
\endproclaim
\bn
\margintag{8.15}\ub{\stag{8.15} Question}:  1) Consistently (with ZF), in some uncountable
power, does Th$(\Bbb Q,<)$ has exactly 3 models. 
\bigskip

\demo{Proof}  1) Follows by (2).
\nl
2) Let $X$ be a set.  For $\ell=0,1,2,3$ let 

$$
\align
K_\ell = \{N:\, &(a) \quad N \text{ is a model of } T_1 
\text{ with universe } X; \\
  &(b) \quad \text{ if } \ell=1 \text{ then } N \text{ omit } p(x) =
  \{c_n < x:n < \omega\} \\
  &(c) \quad \text{ if } \ell=2 \text{ some } a \in N \text{ realizes
  } p(x) \text{ but no} \\
  &\,\,\qquad \quad a \in N \text{ is the first such element}; \\
  &(d) \quad \text{ if } \ell=3 \text{ some element } a \in N \text{
  realizes } p(x) \text{ and is the} \\
  &\,\,\qquad \quad \text{ first such element}\}.
\endalign
$$
\mn
Clearly
\mr
\item "{$\boxplus$}"  $(a) \quad K_0$ is the class of models of $T_1$
with universe $X$
\sn
\item "{${{}}$}"  $(b) \quad K_0$ is the disjoint union of
$K_1,K_2,K_3$.
\ermn
By $(*)_1,(*)_2,(*)_3$ below the result follows:
\mr
\item "{$(*)_1$}"  if $K_2 \ne \emptyset$ then $K_3 \ne \emptyset$.
\nl
[Why?  Let $M \in K_3$ and we define a $\tau(T_1)$-model $N$ as
follows:
{\roster
\itemitem{ $(i)$ }  the universe of $N$ is $X = |M|$
\sn
\itemitem{ $(ii)$ }  $c^N_n = c^M_{n+1}$ for $n < \omega$
\sn
\itemitem{ $(iii)$ }  $N \models a < b$ \ub{iff} $M \models ``a < b
\wedge a \ne c_0 \wedge b \ne c_0$ or $a = c^M_0 \wedge$ ($b$
realizes $p(x)$ in $M$) or $b = c^N_0 \wedge a \ne c_0 \wedge \dsize
\bigvee_{m < \omega}  b < c_m$"
\nl
Now check that $N \in K_3$ with $c^M_0$ being the $<^N$-first member
of $X$ realizing $p(x)]$
\endroster}
\item "{$(*)_2$}"  if $K_3 \ne \emptyset$ then $K_1 \ne \emptyset$.
\nl
[Why?  Let $M \in K_3$ and $c \in M$ realizes $p(x)$ be the first
such element.  We define a $\tau(T_1)$-model $N$ by
{\roster
\itemitem{ $(i)$ }  the universe of $N$ is $X = |M|$
\sn
\itemitem{ $(ii)$ }  $c^N_n = c^M_n$
\sn
\itemitem{ $(iii)$ }  $N \models a < b$ \ub{iff} 
$M \models ``a < b < c"$ or $M \models ``c < b <
a"$ or $M \models ``c \le a \wedge b < c"$.
\endroster}
Now check that $N \in K_1$.]
\sn
\item "{$(*)_3$}"  if $K_1 \ne \emptyset$ then $K_2 \ne \emptyset$.
\nl
[Why?  Let $M \in K_1$, let $Y = \{a \in X:M \models ``c_{2n+1} \le a <
c_{2n+2}"$ for some $n < \omega\}$ 
 and we define a $\tau(T_1)$-model $N$
{\roster
\itemitem{ $(i)$ }   the universe of $N$ is $X=|M|$
\sn
\itemitem{ $(ii)$ }  $c^N_n \equiv c^M_{2n}$ for $n < \omega$
\sn
\itemitem{ $(iii)$ }  $N \models a < b$ \ub{iff} $M \models ``a < b
\wedge (a \notin Y) \wedge (b \notin Y)"$ or
$M \models ``a < b \wedge a \in Y \wedge b \in Y"$ or 
$(b \in Y) \wedge (a \notin Y)$.
\endroster}
Now check that $N \in K_2$.]
\ermn
3) Let $M \in K_1$ have universe $X$ and stipulate $c_{-1} = -
\infty,X_n = \{a:M \models c_{n-1} < a \le c_n\}$ for $n < \omega$ so
$\langle X_n:n < \omega\rangle$ is a partition of $X$.
\nl
Let $S_* = \{n <\omega:X_n$ is uncountable$\}$.
\bn
\ub{Case 1}:  $S_*$ is infinite.

For any partition $\langle S_n:n < \omega\rangle$ of $\omega$ to
infinite sets we can define $N \in K_1$ with universe $X$ such that
$\{c^N_n:n < \omega\} = \{c^M_n:n \in S_0\}$, on this set $<^M,<^N$
agree, and the set $\{n < \omega:(c_n,c_{n+1})_N$ is uncountable$\}$ is
any infinite co-infinite set.
\bn
\ub{Case 2}:  $S_*$ is finite.

We can find $N \in K_0$ with universe $X$ such that
max$\{n:(c_n,c_{n+1})_N$ is uncountable$\}$ is any natural number.
\hfill$\square_{\scite{8.14}}$
\enddemo
\bigskip

\definition{\stag{8.19} Definition}  1) Let $N$ is $(\Bbb
L_{\infty,\kappa},\lambda)$-interpretable in $M$ means (without loss
of generality $\tau_N$ consist of predicates only): there is $\bar d \in
{}^{\lambda >} M$ and sequence $\langle \varphi_R(\bar x_R,\bar d):R
\in \tau_N \rangle$, including $R$ being equality such that

$$
\varphi_R(\bar x_R,\bar y) \in \Bbb L_{\infty,\kappa}
$$

$$
\ell g(\bar x_R) = \text{ arity}(R)
$$

$$
|N| = \{a \in M:M \models \varphi_= (a,a,\bar d)\}
$$

$$
R^N = \{\bar a \in {}^{\ell g(x_R)}|M|:M \models \varphi_R(\bar a,\bar
d)\}.
$$
\mn
2) We add ``fully" if $\varphi_= (\bar x_R) 
= (x_0=x_1)$ for $R$ being the equality.
\enddefinition
\bigskip

\proclaim{\stag{8h.26} Claim}  1) To prove the consistency of ``a first
order complete $T$ is categorical in some power $\ne \aleph_0$" it is
enough
\mr
\item "{$(*)$}"  find a model $N$ of $T$ and $\kappa > \aleph_0$ satisfying:
 if $M$ is a model of $T$ fully $\Bbb
L_{\infty,\kappa}(\tau_M)$-interpretable in $N$ \ub{then} $M \cong N$;
moreover there is a function which is 
$\Bbb L_{\infty,\kappa}(\tau_M)$-definable in $N$
(with $< \kappa$ parameters)  and is an isomorphism from $N$ onto $M$.
\ermn
2) We can replace $\Bbb L_{\infty,\kappa}(N)$ by: there is a set
${\Cal F}$ such that
\mr
\item "{$(a)$}"   ${\Cal F} \subseteq \{f:f$ a partial automorphism of
$N$ with domain of cardinality $< \kappa\}$
\sn
\item "{$(b)$}"   $(\forall A \subseteq N)(|A| < \kappa \Rightarrow
(\exists f \in {\Cal F})(A \subseteq \text{\rm Dom}(f))$
\sn
\item "{$(c)$}"   ${\Cal F}$ closed under inverse and composition
\sn
\item "{$(d)$}"   if $f \in {\Cal F},A \in N$ then $(\exists g \in
{\Cal F})(f \subseteq g \cap a \in \text{\rm Dom}(g))$.
\endroster
\endproclaim
\bigskip

\demo{Proof}  Straight.  
\enddemo
\bigskip

\remark{Remark}  So this categoricity does not imply ``not complicated".
\hfill$\square_{\scite{8h.26}}$ 
\endremark
\newpage
    
REFERENCES.  
\bibliographystyle{lit-plain}
\bibliography{lista,listb,listx,listf,liste}

\def\germ{\frak} \def\scr{\cal} \ifx\documentclass\undefinedcs
  \def\bf{\fam\bffam\tenbf}\def\rm{\fam0\tenrm}\fi 
  \def\defaultdefine#1#2{\expandafter\ifx\csname#1\endcsname\relax
  \expandafter\def\csname#1\endcsname{#2}\fi} \defaultdefine{Bbb}{\bf}
  \defaultdefine{frak}{\bf} \defaultdefine{=}{\B} 
  \defaultdefine{mathfrak}{\frak} \defaultdefine{mathbb}{\bf}
  \defaultdefine{mathcal}{\cal}
  \defaultdefine{beth}{BETH}\defaultdefine{cal}{\bf} \def\bbfI{{\Bbb I}}
  \def\mbox{\hbox} \def\text{\hbox} \def\om{\omega} \def\Cal#1{{\bf #1}}
  \def\pcf{pcf} \defaultdefine{cf}{cf} \defaultdefine{reals}{{\Bbb R}}
  \defaultdefine{real}{{\Bbb R}} \def\restriction{{|}} \def\club{CLUB}
  \def\w{\omega} \def\exist{\exists} \def\se{{\germ se}} \def\bb{{\bf b}}
  \def\equivalence{\equiv} \let\lt< \let\gt>
  \def\implies{\Rightarrow}\def\mathfrak{\bf}\def\germ{\frak} \def\scr{\cal}
  \ifx\documentclass\undefinedcs
  \def\bf{\fam\bffam\tenbf}\def\rm{\fam0\tenrm}\fi 
  \def\defaultdefine#1#2{\expandafter\ifx\csname#1\endcsname\relax
  \expandafter\def\csname#1\endcsname{#2}\fi} \defaultdefine{Bbb}{\bf}
  \defaultdefine{frak}{\bf} \defaultdefine{=}{\B} 
  \defaultdefine{mathfrak}{\frak} \defaultdefine{mathbb}{\bf}
  \defaultdefine{mathcal}{\cal}
  \defaultdefine{beth}{BETH}\defaultdefine{cal}{\bf} \def\bbfI{{\Bbb I}}
  \def\mbox{\hbox} \def\text{\hbox} \def\om{\omega} \def\Cal#1{{\bf #1}}
  \def\pcf{pcf} \defaultdefine{cf}{cf} \defaultdefine{reals}{{\Bbb R}}
  \defaultdefine{real}{{\Bbb R}} \def\restriction{{|}} \def\club{CLUB}
  \def\w{\omega} \def\exist{\exists} \def\se{{\germ se}} \def\bb{{\bf b}}
  \def\equivalence{\equiv} \let\lt< \let\gt>
\begin{thebibliography}{LwSh 518}
\makeatletter \renewcommand{\@biblabel}[1]{[#1]} \makeatother
\def\eprintfootnotetext{References of the form {\tt math.XX/$\cdots$}
 refer to {\tt arXiv.org} }
\ifx\documentstyle\undefinedcontrolsequence
   \def\anyfootnote{\footnote{*}}
   \else\def\anyfootnote{\footnote}\fi
\def\eprintfn{\ifEprint\anyfootnote{\eprintfootnotetext}\fi\Eprintfalse }
\newif\ifEprint  \Eprinttrue

\bibitem[BLSh 464]{BLSh:464}John~T. Baldwin, Michael~C. Laskowski, and Saharon
  Shelah.
\newblock {Forcing Isomorphism}.
\newblock {\em {Journal of Symbolic Logic}}, {\bf 58}:1291--1301, 1993.
\newblock math.LO/9301208.

\bibitem[Be84]{Be84}Steven Buechler.
\newblock {Kueker's conjecture for superstable theories}.
\newblock {\em Journal of Symbolic Logic}, {\bf 49}:930--934, 1984.

\bibitem[HHL00]{HHL00}Bradd Hart, Ehud Hrushovski, and Michael~C. Laskowski.
\newblock The uncountable spectra of countable theories.
\newblock {\em Annals of Mathematics}, {\bf 152}:207--257, 2000.

\bibitem[Hr89]{Hr89}Ehud Hrushovski.
\newblock Kueker's conjecture for stable theories.
\newblock {\em Journal of Symbolic Logic}, {\bf 54}:207--220, 1989.

\bibitem[Hr89d]{Hr89d}Ehud Hrushovski.
\newblock Unidimensional theories.
\newblock In {\em Logic Colloquium 88}. North--Holland, 1989.

\bibitem[Ke71a]{Ke71a}Jerome~H. Keisler.
\newblock {\em {On theories categorical in their own power}}, volume~36.
\newblock 1971.

\bibitem[Las88]{Las88}Michael~C. Laskowski.
\newblock {Uncountable theories that are categorical in a higher power}.
\newblock {\em The Journal of Symbolic Logic}, {\bf 53}:512--530, 1988.

\bibitem[LwSh 518]{LwSh:518}Michael~C. Laskowski and Saharon Shelah.
\newblock {Forcing Isomorphism II}.
\newblock {\em {Journal of Symbolic Logic}}, {\bf 61}:1305--1320, 1996.
\newblock math.LO/0011169.

\bibitem[Lv71]{Lv71}Richard Laver.
\newblock {On Fraiss\'e's order type conjecture}.
\newblock {\em Annals of Mathematics}, {\bf 93}:89--111, 1971.

\bibitem[Mo65]{Mo65}Michael Morley.
\newblock {Categoricity in power}.
\newblock {\em Transaction of the American Mathematical Society}, {\bf
  114}:514--538, 1965.

\bibitem[Sh 300f]{Sh:300f}Saharon Shelah.
\newblock {\em {Chapter VI}}.

\bibitem[Sh:e]{Sh:e}Saharon Shelah.
\newblock {\em {Non--structure theory}}, accepted.
\newblock {Oxford University Press}.

\bibitem[Sh 750]{Sh:750}Saharon Shelah.
\newblock {On Weak Beth for Cofinality Logic}.
\newblock {\em Preprint}.

\bibitem[Sh:F728]{Sh:F728}Saharon Shelah.
\newblock {PCF arithmetic with little choice}.

\bibitem[Sh 835]{Sh:835}Saharon Shelah.
\newblock {PCF without choice}.
\newblock {\em Archive for Mathematical Logic}, {\bf submitted}.
\newblock math.LO/0510229.

\bibitem[Sh 3]{Sh:3}Saharon Shelah.
\newblock {Finite diagrams stable in power}.
\newblock {\em {Annals of Mathematical Logic}}, {\bf 2}:69--118, 1970.

\bibitem[Sh 4]{Sh:4}Saharon Shelah.
\newblock {On theories $T$ categorical in $|T|$}.
\newblock {\em {Journal of Symbolic Logic}}, {\bf 35}:73--82, 1970.

\bibitem[Sh 12]{Sh:12}Saharon Shelah.
\newblock {The number of non-isomorphic models of an unstable first-order
  theory}.
\newblock {\em {Israel Journal of Mathematics}}, {\bf 9}:473--487, 1971.

\bibitem[Sh 31]{Sh:31}Saharon Shelah.
\newblock {Categoricity of uncountable theories}.
\newblock In {\em {Proceedings of the Tarski Symposium (Univ. of California,
  Berkeley, Calif., 1971)}}, volume XXV of {\em Proc. Sympos. Pure Math.},
  pages 187--203. {Amer. Math. Soc., Providence, R.I}, 1974.

\bibitem[Sh 52]{Sh:52}Saharon Shelah.
\newblock {A compactness theorem for singular cardinals, free algebras,
  Whitehead problem and transversals}.
\newblock {\em {Israel Journal of Mathematics}}, {\bf 21}:319--349, 1975.

\bibitem[Sh 54]{Sh:54}Saharon Shelah.
\newblock {The lazy model-theoretician's guide to stability}.
\newblock {\em {Logique et Analyse}}, {\bf 18}:241--308, 1975.

\bibitem[Sh:E18]{Sh:E18}Saharon Shelah.
\newblock {A combinatorial proof of the singular compactness theorem}.
\newblock {\em Mineograph notes and lecture in a mini-conference, Berlin,
  August'77}, 1977.

\bibitem[Sh 100]{Sh:100}Saharon Shelah.
\newblock {Independence results}.
\newblock {\em {The Journal of Symbolic Logic}}, {\bf 45}:563--573, 1980.

\bibitem[Sh 199]{Sh:199}Saharon Shelah.
\newblock {Remarks in abstract model theory}.
\newblock {\em {Annals of Pure and Applied Logic}}, {\bf 29}:255--288, 1985.

\bibitem[Sh:c]{Sh:c}Saharon Shelah.
\newblock {\em {Classification theory and the number of nonisomorphic models}},
  volume~92 of {\em {Studies in Logic and the Foundations of Mathematics}}.
\newblock {North-Holland Publishing Co., Amsterdam, xxxiv+705 pp}, 1990.

\bibitem[Sh 497]{Sh:497}Saharon Shelah.
\newblock {Set Theory without choice: not everything on cofinality is
  possible}.
\newblock {\em Archive for Mathematical Logic}, {\bf 36}:81--125, 1997.
\newblock A special volume dedicated to Prof. Azriel Levy. math.LO/9512227.

\bibitem[Sh:E38]{Sh:E38}{Shelah, Saharon}.
\newblock {Continuation of 497: Universes without Choice}.

\bibitem[Sh:F701]{Sh:F701}{Shelah, Saharon}.
\newblock More on model theory without choice.

\bibitem[WT05]{WT05}Agatha Walczak-Typke.
\newblock {The first-order structure of weakly Dedekind-finite sets}.
\newblock {\em Journal of Symbolic Logic}, {\bf 70}:1161--1170, 2005.

\bibitem[WT07]{WT07}Agatha Walczak-Typke.
\newblock A model-theoretic approach to structures in set theory without the
  axiom of choice.
\newblock In B.~Loewe, editor, {\em Algebra, Logic, Set Theory: Festschrift fur
  Ulrich Felgner zum 65 Geburtstag}, Studies in Logic. College Publications at
  Kings College London. to appear.

\end{thebibliography}

\enddocument